\documentclass[12pt,a4paper]{article}
\usepackage{amsmath,amssymb}
\usepackage{booktabs}
\usepackage{latexsym}
\usepackage{graphicx} 
\usepackage{epsfig}
\usepackage{epstopdf} 
\usepackage[title]{appendix}
\usepackage{color}
\usepackage[dvipsnames]{xcolor} 
\usepackage{multirow}
\usepackage{indentfirst}
\usepackage{float}                  
\usepackage{subfig}                 
\usepackage{overpic}                

\usepackage{mathrsfs}
\usepackage{amsmath}
\usepackage{siunitx} 
\usepackage{bm}
\usepackage{appendix}
\usepackage[pdftex]{hyperref}
\usepackage[numbers,sort&compress]{natbib}

\definecolor{yellowgreen}{RGB}{70, 210, 10}
\renewcommand{\theequation}{\thesection.\arabic{equation}}
\renewcommand{\thefigure}{\thesection.\arabic{figure}}
\renewcommand{\thetable}{\thesection.\arabic{table}}

\begin{document}
	\baselineskip=2pc
	
	\begin{center}
		{\large \bf  Uniformly High Order Discontinuous Galerkin Gas Kinetic Scheme for Compressible flows
}
	\end{center}

\centerline{
Mengqing Zhang%
\footnote{Institute of Applied Physics and Computational Mathematics and National Key laboratory of Computational physics, 
Beijing 100088, China. E-mail: zhangmengqing22@gscaep.ac.cn.},
Shiyi Li%
\footnote{Institute of Applied Physics and Computational Mathematics and National Key laboratory of Computational physics, 
Beijing 100088, China. E-mail: lishiyi14@tsinghua.org.cn.},
Dongmi Luo%
\footnote{Institute of Applied Physics and Computational Mathematics and National Key laboratory of Computational physics, 
Beijing 100088, China. E-mail: dongmiluo@stu.xmu.edu.cn.},
Jianxian Qiu%
\footnote{School of Mathematical Sciences and Fujian Provincial Key Laboratory of Mathematical Modeling and High-Performance Scientific Computing, Xiamen University, 
Xiamen, Fujian 361005, China. E-mail: jxqiu@xmu.edu.cn. },
Yibing Chen*%
\footnote{Institute of Applied Physics and Computational Mathematics and National Key laboratory of Computational physics, 
and Center for Applied Physics and Technology, Peking University,
Beijing 100088, China. E-mail: chen\_yibing@iapcm.ac.cn.},
}

\vspace{20pt}

\begin{abstract}

In this paper, a uniformly high-order discontinuous Galerkin gas kinetic scheme (DG-HGKS) is proposed to solve the Euler equations of compressible flows. 
The new scheme is an extension of the one-stage compact and efficient high-order GKS (CEHGKS, Li et al. , 2021. J. Comput. Phys. 447, 110661) in the finite volume framework. 
The main ideas of the new scheme consist of two parts. 
Firstly, starting from a fully discrete DG formulation, the numerical fluxes and volume integrals are expanded in time. 
Secondly, the time derivatives are replaced by spatial derivatives using the techniques in CEHGKS. 
To suppress the non-physical oscillations in the discontinuous regions while minimizing the number of ”troubled cells”, 
an effective limiter strategy compatible with the new scheme has been developed by combining the KXRCF indicator and the SHWENO reconstruction technique. 
The new scheme can achieve arbitrary high-order accuracy in both space and time, 
thereby breaking the previous limitation of no more than third-order accuracy in existing one-stage DG-HGKS schemes. 
Numerical tests in 1D and 2D have demonstrated the robustness and effectiveness of the scheme.
\end{abstract}

\pagenumbering{arabic}

\section{Introduction}
 
In the past few decades, high-order schemes for the Euler equations of compressible flows have been rapidly developed.
Many schemes have made great success, including the  Essentially Non-Oscillatory scheme (ENO) scheme \cite{ENO} , the Weighted Essentially Non-Oscillatory  (WENO) scheme \cite{WENO}, the Discontinuous Galerkin method (DG) \cite{DG}, the  Residual Distribution method (RD) \cite{RD} and the Spectral Volume method (SV) \cite{SV},
which can achieve arbitrary high-order accuracy in space.
Among these methods, 
the DG scheme is one of the most widely used numerical methods due to its compactness, as it requires no additional computational stencils even when accuracy increases. 
In smooth regions, the DG method relies solely on information from the local element,
while near discontinuities, it only uses data from adjacent elements.
Moreover, it allows the simultaneous update of both conservative variables and their high-order moments.
One of the primary challenges of the DG scheme is to avoid non-physical oscillations while maintaining high order accuracy and compactness when addressing discontinuous problems.
The predominant strategy involves using detectors to identify "troubled cells" that may necessitate the application of monotonicity limiters \cite{RKDG2,RKDG3,RKDG4,RKDG5,WL1,WL2}, 
followed by the employment of high-order, non-oscillatory, and compact reconstruction methods, 
including but not limited to Hermite WENO (HWENO) \cite{HWL1,HWL2} 
and Simple HWENO (SHWENO) \cite{SWENOL,SHWENOL,MRWENOL}.
The most widely used  discontinuity detecting methods to identify "troubled cells" include the minmod-based TVB limiter  \cite{HWL2} and a shock-detection technique KXRCF indicator \cite{KXRCF}.

In the context of time discretization, the multi-stage Runge-Kutta methods are commonly employed to ensure the stability of numerical schemes.
However, it has been established that a Total Variation Diminishing (TVD) Runge-Kutta method with positive coefficients cannot achieve an accuracy higher than fourth-order \cite{Timebarriers}.
Even if the TVD requirements are relaxed, a fifth-order or higher Runge-Kutta method would require a corresponding increase in the order of spatial discretization. 
Moreover, Runge-Kutta methods of order greater than third require more stages. 
In most high-order schemes, the high-order reconstruction process is the most computationally expensive component. 
Furthermore, at each stage of the Runge-Kutta method, a reconstruction process is necessitated.
Consequently, the efficiency of Runge-Kutta methods diminishes for orders beyond third \cite{Balsara}.
When employing the most common third-order TVD Runge-Kutta method, 
the overall accuracy is limited to third-order, which is known as the accuracy barrier \cite{ADER1,ADER}.

To overcome these drawbacks, 
{much attention has been devoted to the development of one-stage schemes that can achieve arbitrary high-order accuracy uniformly in both  spatial and  temporal dimensions.
These methods are known as uniformly  high-order schemes.}
The basic idea of one-stage schemes be traced back to the Lax-Wendroff procedure, and various techniques have been developed for nonlinear problems \cite{Qiu}, 
such as the generalized Riemann problem (GRP) \cite{GRP}, the scheme with arbitrary high order derivative (ADER) \cite{ADER} method,
arbitrary high order approach based on flux vector splitting (HFVS) \cite{HFVS},
 and the high order gas-kinetic scheme (HGKS) \cite{13,15,R1,31,32}. 
{ 
It is noteworthy that the HGKS leverages an underlying kinetic model,
which offers a natural physical mechanism for dissipation in discontinuous regions, 
thus ensuring inherent numerical stability \cite{GKS.GRP}.}

Some of these one-stage schemes have been extended to the DG framework \cite{Qiu,Dum1,Dum2}, with the aim of achieving uniformly high-order accuracy in both space and time.
Research on DG-HGKS is relatively limited. 
Xu \cite{Xu} first proposed a one-dimensional DG method  that employs a second-order BGK scheme for flux calculation, which overall accuracy cannot exceed second-order. 
Ni et al. \cite{NGX} further extended this scheme to a two-dimensional case, applying a WENO limiter to suppress oscillations. However, this approach destroy the original compactness of DG scheme. 
Ren et al. \cite{RXD} proposed a DG method based on a multi-dimensional third-order GKS for viscous flow calculations. 
They utilized a linear least squares method to reconstruct the equilibrium distribution at cell interfaces and applied a compact limiter strategy to address discontinuous problems. 
The development of DG-HGKS is primarily limited by two key factors: 
First, the HGKS formulation within the original finite volume framework is already highly complex, and extending it to the DG framework would further exacerbate this complexity, making it challenging to achieve third-order or higher accuracy. 
Second, advanced detectors and limiters have not yet been able to effectively adapt to the DG-HGKS.

In recent years, significant progress has been made in HGKS. 
The efficient high-order gas-kinetic scheme (EHGKS) proposed by Li et al. \cite{LSY} significantly reduces the complexity of the original HGKS by eliminating unnecessary high-order dissipation terms and incorporating the Lax-Wendroff procedure, thereby achieving uniformly arbitrary high-order accuracy in both space and time. 
Subsequently, Li et al. \cite{LSY2} combined EHGKS with a  SHWENO  reconstruction technique to construct a more compact and efficient high-order gas kinetic scheme (CEHGKS).

 
Inspired by the recent advancements in HGKS and compact reconstruction techniques SHWENO, 
this article aims to explore the design of a robust, one-stage uniformly arbitrary high-order DG-HGKS scheme utilizing compact stencils. 
The construction of the new scheme can be divided into two steps. 
The first step involves starting from the fully discrete DG weak form and performing Taylor expansion of the time-dependent numerical fluxes and volume integrals. 
The numerical flux terms are then directly solved using the CEHGKS algorithm, 
while the volume integrals are solved using the linearized ADER algorithm, 
thereby replacing the time derivatives with spatial derivatives. 
This results in the fully discrete form of the new one-stage DG-HGKS. 
The second step is to calculate the spatial derivative terms. 
In smooth regions, since the conservative variables and their high-order moments can be updated simultaneously within the DG framework, 
the required spatial derivatives can be obtained directly. 
This is a significant difference from the finite volume CEHGKS. 
For regions containing discontinuities, an effective limiter strategy compatible with the new scheme has been carefully developed by combining the KXRCF indicator and the SHWENO reconstruction technique. 
The new scheme has the following advantages: 
Firstly, it is the first one-stage DG-HGKS scheme to achieve uniformly arbitrary high-order accuracy in both space and time, 
breaking the limitation of existing  schemes,
which are constrained  to a maximum of third-order accuracy. 
Secondly, it avoids non-physical oscillations even when the detection range is significantly reduced, while maintaining robustness and better compactness. 
The number of "troubled cells" is also greatly reduced compared to RKDG schemes using similar strategies \cite{SHWENOL}.
It is particularly worth noting that this is also the main reason why the new DG-HGKS yields better computational results than the finite volume CEHGKS.
 

The paper is organized as follows. 
The one-stage uniformly high-order DG-HGKS scheme using the KXRCF indicator and SHWENO reconstruction technique is described in detail for one-dimensional and two-dimensional Euler systems in Sections 2 and 3, respectively. In Section 4, one- and two-dimensional numerical examples are presented to demonstrate the good performance of the new scheme. Finally, the conclusions are drawn in Section 5.

\section{The numerical method in one-dimensional case}   

In this section, we describe the one-stage uniformly high-order DG-HGKS scheme in detail for one-dimensional case.

\subsection{The BGK model and the Euler systems }
Let's introduce the relationship between the mesoscopic BGK model and the macroscopic Euler systems \cite{Relationship},
the   BGK model \cite{31} is given by
\begin{eqnarray}
\label{eq:bgk}
   \frac{\partial f}{\partial t}
+u \frac{\partial f}{\partial x}
 = \frac{g - f}{\tau } .
\end{eqnarray}
Here, $f$ represents the particle distribution function, 
and $g$ is the corresponding equilibrium distribution.
Both $f$ and $g$ are functions of $x,t,u,\bm{\xi}$,
where  ${ u}$ is the particle velocity,
$\bm{\xi}$ represents the internal variables.
{The  degrees of freedom $K$ of $\bm {\xi}$ is  equal to $(3-\gamma)/(\gamma-1)$,}
with   the specific  heat ratio $\gamma=1.4$.
$\tau=\mu/p$ is the mean collision time,
$\mu$ is the dynamical viscosity, 
and $p$ is the pressure.
{The equilibrium distribution function $g$ is the Maxwellian  
}
\begin{eqnarray}
\label{eq:bgkg}
g={\rho
\left( \frac{\lambda}{\pi } \right)^{\frac{K+1}{2}}}
e^{ - \lambda  \left( ({ u} - { U})^2+ {\bm {\xi}}^2 \right)} ,
\end{eqnarray}
{where $\rho,U$ represent density and macroscopic velocity, respectively.}
 $\lambda=1/(2RT)$,  where $T$ is the temperature and $R$  is the gas constant.
{In fact, $g$ can be determined by  the  macroscopic  conservative variables $ W=(\rho,\rho U,E)^{T}$.}
The total energy given by $E=\frac{p}{\gamma-1}+\frac{1}{2}\rho U^2$.
The time integral solution of the BGK model is given by \cite{32}
\begin{eqnarray}
f({ x},t,{ u},\bm {\xi}) =
\frac{1}{\tau} \int_0^t
{g({ x}',t',{ u},\bm {\xi})  e^{-(t - t')/\tau} { d t'}}
+ e^{-t/\tau} f({ x} - { u} t,0,{ u},\bm{ \xi}),
\label{eq:ints}
\end{eqnarray}	
{where ${ x}'= { x} - { u} (t-t')$ is the particle trajectory.
The solution $f$ depends on the equilibrium state and the initial gas distribution function at the beginning of each time step $t = 0$.}

The BGK model \eqref{eq:bgk} and its integral solution \eqref{eq:ints} provide a mesoscopic description of gas kinetics,
and the statistics of the particle distribution function $f$  provide a macroscopic description of the flow structure. 
The relationship between the particle distribution function $f$ and the macroscopic conserved quantity $W$   given by \cite{Relationship}
\begin{eqnarray}
W =  \int {f{\Psi}{d \Xi}} ,
\label{eq:rew}
\end{eqnarray}
{where $\Psi$ is the vector of moments  
\begin{equation}{\Psi} = ({1, { u}, ({ u}^2+{\bm \xi}^2 )/2})^{T},
\label{mom-psi}
\end{equation}}
and $d \Xi  = {d { u}}{d {\bm \xi}}$.
By taking the moments of the BGK model \eqref{eq:bgk} on $\Psi$,
the macroscopic conservative equations \eqref{eq1} can be derived 
\begin{eqnarray}
		\begin{aligned}
				\frac{\partial W}{\partial t}+\frac{\partial F(W)}{\partial x}=0 .
		\end{aligned}
		\label{eq1}
	\end{eqnarray}
The relationship between the particle distribution function $f$ and the macroscopic flux $F$   given by \cite{Relationship}
\begin{eqnarray}
F=\int { { u}f {\Psi}{d \Xi}} .
\label{eq:rew}
\end{eqnarray}
In particular, when $f=g$ in  \eqref{eq:bgk},
according to the moment of the Maxwellian equilibrium state \eqref{eq:bgkg},
the macroscopic equation \eqref{eq1}  represents the Euler systems
\begin{eqnarray}
\label{eq:euler}
\frac{\partial}{\partial t}\left(\begin{matrix}\rho\\\rho U\\E\\\end{matrix}\right)+\frac{\partial}{\partial x}\left(\begin{matrix}\rho U\\\rho U^2+p\ \ \\U(E+p)\\\end{matrix}\right)=0 .
\end{eqnarray}
 $F=\int { g{ u} {\Psi}{d \Xi}}=F(W)=(\rho U, \rho U^2+p, U(E+p))^{T}$ is the conservative flux.

Next, we will introduce the calculation process  of the fully discrete DG-HGKS scheme.

\subsection{The DG-HGKS scheme }

Given the space domain $[a, b]$ a uniform mesh division $a= x_{\frac{1}{2}} <  \cdots < x_{N+\frac{1}{2}}=b$. 
The cell $I_i=(x_{i-\frac{1}{2}}, x_{i+\frac{1}{2}})$, for $i=1, \cdots, N $, 
where $N$ is the number of cells, 
the cell center  $x_{i}=\frac{1}{2}(x_{i-\frac{1}{2}}+x_{i+\frac{1}{2}})$ 
and cell size $\Delta x=x_{i+\frac{1}{2}}-x_{i-\frac{1}{2}}$.
The  computational time $T$ is discretized into $0=t_0 < t_1 < \cdots < t_n < t_{n+1} < \cdots = T$.
 
First, we introduce the numerical solution under DG framework  briefly. 
The DG space is
$$V_h\equiv V_{h}^{\hspace{0.1em}k}=\left\{v\in L^2\left(I\right):v|_{I_i}\in P^k\left(I_i\right),\ i=1,\cdots,N\right\} ,$$ 
where $P^k\left(I_i\right)$ is the set of polynomials of  degree at most $k$ in  cell $I_i$.
The numerical solution is written as 
\begin{equation}\label{so1d}
W_h(x,t)=\sum_{m=0}^{k}{w_{i,\hspace{0.1em}m}(t){v_{i,\hspace{0.1em}m}(x)}}=\sum_{m=0}^{k}{w_{i,\hspace{0.1em}m}(t)P_m({\zeta})} ,
\end{equation}
where $P_m({\zeta})$ is an orthogonal Legendre polynomial of degree $m$ on the standard interval $\left[-1,1\right]$, 
after interval transformation by ${\zeta}=\frac{2(x-x_i)}{ \Delta x}$, 
{we obtain the test function
$v_{i,\hspace{0.1em}m}(x)$, $m=0,1,\cdots,k$.} 
{By orthogonality,
denote the diagonal matrix of mass 
\begin{equation}
M = \begin{pmatrix}
1 & 0 & \cdots & 0 \\
0 & \frac{1}{3} & \cdots & 0 \\
 \vdots & \vdots & \ddots & \vdots \\
0 & 0 & \cdots & \frac{1}{2k + 1} 
\end{pmatrix} ,
\label{DiagM}
\end{equation}}
{and the degrees of freedom $w_{i,\hspace{0.1em}m}(t)$ are the moments defined by
\begin{equation}\label{moment1d}
w_{i,\hspace{0.1em}m}(t)=\frac{2m+1}{\Delta x}\int_{I_i}W_h(x,t)v_{i,\hspace{0.1em}m}(x)dx,m=0,...,k.
\end{equation}
}
{For convenience, the subscript $i$ is omitted in the following, 
which are $v_m(x)$ and $w_m(x)$.}

The Euler systems \eqref{eq1} is multiplied by a test function $v_m(x)$ and integrated over each element $I_i \times [t_n,~t_{n+1}]$.
After integration by parts, we  can obtain the fully discrete DG scheme
	\begin{eqnarray}
		\begin{aligned}
     {\int_{t_{n}}^{t_{n + 1}}{{\int_{x_{i - \frac{1}{2}}}^{x_{i + \frac{1}{2}}}{\left( W_{h} \right)_{t}v_m(x)}}dx}}dt 
 - {\int_{t_{n}}^{t_{n + 1}}{{\int_{x_{i - \frac{1}{2}}}^{x_{i + \frac{1}{2}}}F ( W_{h}  ) ( v_m(x)  )_{x}}dx}}dt \\
  +{\int_{t_{n}}^{t_{n + 1}}{F ( x_{i + \frac{1}{2}}, t  )v_m ( x_{i + \frac{1}{2}}^{-}  )}}dt 
 - {\int_{t_{n}}^{t_{n + 1}}{F ( x_{i - \frac{1}{2}}, t  )v_m ( x_{i - \frac{1}{2}}^{+}  )}}dt =0 .
		\end{aligned}
		\label{eq2}
	\end{eqnarray}
Then transform the cell $I_i$ into the reference interval $\left[-1,1\right]$ and use the orthogonal property of  the Legendre polynomial to obtain
\begin{eqnarray}
		\begin{aligned}
{\int_{t_{n}}^{t_{n + 1}}{{\int_{x_{i - \frac{1}{2}}}^{x_{i + \frac{1}{2}}}{\left( W_{h} \right)_{t}v_m(x)}}dx}}dt 
&={\int_{t_{n}}^{t_{n + 1}}{{\int_{-1}^{1}{\left( \sum_{l=0}^{k}{w_l(t)P_l({\zeta})} \right)_{t}	P_m({\zeta})}}\frac{\Delta x}{2}d{\zeta}}}dt \\
&=\frac{\Delta x}{2m+1}\left( w_m(t_{n+1})-w_m(t_n)\right), m=0,...,k ,
             \end{aligned}
		\label{1-1}
\end{eqnarray}
and \eqref{eq2} equivalent to
 \begin{eqnarray}
		\begin{aligned}
  w_m(t_{n+1}) &= w_m(t_n) + \frac{2m+1}{\Delta x} {\int_{t_{n}}^{t_{n + 1}}{{\int_{x_{i - \frac{1}{2}}}^{x_{i + \frac{1}{2}}}{F( W_{h} )(v_m(x) )}_{x}}dx}}dt \\
&-\frac{2m+1}{\Delta x}\left( {\int_{t_{n}}^{t_{n + 1}}{F( x_{i + \frac{1}{2}}, t)v_m( x_{i + \frac{1}{2}}^{-} )}}dt - {\int_{t_{n}}^{t_{n + 1}}{F( x_{i - \frac{1}{2}}, t )v_m( x_{i - \frac{1}{2}}^{+} )}}dt \right) .
         \end{aligned}
		\label{eq2-3}
\end{eqnarray}

We need to solve the moments  ${w_m(t_{n+1})}, m=0,...,k$.
The following  will describe the  calculation of numerical flux terms  ${\int_{t_{n}}^{t_{n + 1}}{F( x_{i + \frac{1}{2}}, t)v_m( x_{i + \frac{1}{2}}^{-} )}}dt $ along the cell interfaces,
and the volume integral terms ${\int_{t_{n}}^{t_{n + 1}}{{\int_{x_{i - \frac{1}{2}}}^{x_{i + \frac{1}{2}}}{F( W_{h} )(v_m(x) )}_{x}}dx}}dt $  inside the cells. 

In order to preserve the conservative, the calculation of the numerical  flux $F( x_{i + \frac{1}{2}}, t )$ and $F( x_{i - \frac{1}{2}}, t )$ are similar, here we only introduce $F_{i + \frac{1}{2}}(t)=F( x_{i + \frac{1}{2}}, t )$.
Since the basis function  depended only on the spatial variable $x$, we have

\begin{eqnarray}
\int_{t_{n}}^{t_{n + 1}}{F( x_{i + \frac{1}{2}}, t)v_m( x_{i + \frac{1}{2}}^{-} )}dt=v_m( x_{i + \frac{1}{2}}^{-} )\int_{t_{n}}^{t_{n + 1}}{F_{i+\frac{1}{2}}(t)}dt.
\end{eqnarray}
		
The  numerical flux is  directly calculated  from the existing EHGKS flux, 
which is also used in the CEHGKS scheme. 
We will briefly describe this process,
more details refer to \cite{LSY,LSY2}.

The numerical flux is  a combination of equilibrium and non-equilibrium states
\begin{eqnarray}
		\begin{aligned}
F_{i + \frac{1}{2}}(t) = \left( 1 - e^{- {\Delta}t/\tau} \right)\mathcal{P}_{r}F^{e} + e^{-{\Delta}t/\tau}F^{k} ,
\end{aligned}
		\label{EHGKS}
	\end{eqnarray}
{here 
$
\tau  =\epsilon _1 \Delta t
+ \epsilon_2 \frac{{\left|p^{\rm L}-p^{\rm R} \right|}}{p^{\rm L}+p^{\rm R} }\Delta t
$,}
{where  the constant parameters $\epsilon _1=0.02$ and $\epsilon _2=2 $\cite{LSY},}
and $\mathcal{P}_{r} = {\sum\limits_{l = 0}^{r}{\frac{t^{l}}{l!}\frac{\partial^{l}}{\partial t^{l}}}}$,
{$r$ depends on the specific  order, for example, when the numerical scheme is of the fifth-order,
$r = 4$.}
$F^{k}$ is the second-order KFVS flux \cite{31}.
 {In smooth regions where $\tau$ is small, $F_{i + \frac{1}{2}}(t) \xrightarrow{}{\mathcal{P}_{r}} F^{e}$.
 While near discontinuities, it  approaches to $F_{i + \frac{1}{2}}(t) \xrightarrow{}F^{k}$.}

In the  high-order equilibrium state  $\mathcal{P}_{r}F^{e} $ term,
bypass the calculations involving the derivatives of   $g^{e}$,
{In fact, after taking the moments in $F^{e}$, }
$F^{e}$ can be converted into the macroscopic flux	 
\begin{eqnarray}
		\begin{aligned}
F^{e} = F_{Eu}\left( W^{e} \right) = \begin{pmatrix}
{\rho^{e}U^{e}} \\
{\rho^{e}U^{e2} + p^{e}} \\
{\rho^{e}{E^{e}U}^{e} + p^{e}U^{e}} 
\end{pmatrix} ,
         \end{aligned}
		\label{eq6}
\end{eqnarray}
{and consequently
 $
 \mathcal{P}_{r}F^{e} = \mathcal{P}_{r}F_{Eu}\left( W^{e} \right)
 $,}
the time integrated flux is evaluated as
\begin{eqnarray}
		\begin{aligned}
\mathbb{F}_{i + \frac{1}{2}}\left( {{\Delta}t} \right) = {\int_{0}^{{\Delta}t}{F_{i + \frac{1}{2}}(t)}}dt = \left( 1 - e^{- {\Delta}t/\tau} \right)\mathbb{F}_{i + \frac{1}{2}}^{e}\left( {{\Delta}t} \right) + e^{- {\Delta}t/\tau}\mathbb{F}_{i + \frac{1}{2}}^{k}\left( {{\Delta}t} \right) .
         \end{aligned}
		\label{eq8}
\end{eqnarray}
$\mathbb{F}_{i + \frac{1}{2}}^{k}\left({{\Delta}t}\right)$  is the time integral of the KFVS flux.
{Where $\mathbb{F}_{i+\frac12}^e(\Delta t)$ is the time integral of $\mathcal{P}_{r}F_{Eu}(W^e)$,
$$
\mathbb{F}_{i+\frac{1}{2}}^e(\Delta t)=\int_0^{\Delta t}\mathcal{P}_r\mathbf{F}_{Eu}({W}^e)\mathrm{d}t.
$$
By the Gaussian rule and} 
based on  primitive variables ${Q}^e=(\rho^e, U^e, p^e)^\mathrm{T}$ to approximate $\mathbb{F}_{i+\frac12}^e(\Delta t)$, we have
\begin{eqnarray}
		\begin{aligned}
\mathbb{F}_{i + \frac{1}{2}}^{e}\left( { {\Delta}t} \right) 
\approx {\sum\limits_{\alpha = 0}^{K}{F_{Eu}\left( {\mathcal{Q}^{e}\left( \kappa_{\alpha} {\Delta}t \right)} \right)\omega_{\alpha}}} , 
         \end{aligned}
		\label{flux1}
\end{eqnarray}
here
\begin{eqnarray}
		\begin{aligned}
\mathcal{Q}^{e}(t) &= {\sum\limits_{l = 0}^{r}{\frac{t^{l}}{l!}\frac{\partial^{l}Q^{e}_{i+\frac{1}{2}}}{\partial t^{l}}}} ,
         \end{aligned}
		\label{flux2}
\end{eqnarray}
where 
$\kappa_{\alpha}, \omega_{\alpha}$ are the Gaussian nodes and weights. 
The time related derivatives of $\mathcal{Q}^{e}$ are obtained from the spatial derivatives according to 
{
$$
\left.\frac{\partial^{m+q}}{\partial x^mt^q}\left(\frac{\partial\mathcal{Q}^e}{\partial t}+\mathbf{A}^e\cdot\frac{\partial\mathcal{Q}^e}{\partial x}\right)=\mathbf{0},\quad\mathbf{A}=\left[
\begin{array}
{ccc}U & \rho & 0 \\
0 & U & 1/\rho \\
0 & \gamma p & U
\end{array}\right.\right],
$$
this is also  equivalent to the Lax-Wendroff procedure,} which has been widely used in ADER \cite{ADER}.
As for each order of spatial derivatives can be calculated by{\cite{LSY,LSY2}}
 \begin{eqnarray}
		\begin{aligned}
\frac{\partial^{l}W^{e}}{\partial x^{l}} = \omega^{e}\frac{\partial^{l}W^{L}}{\partial x^{l}} + \left( 1 - \omega^{e} \right)\frac{\partial^{l}W^{R}}{\partial x^{l}} ,
         \end{aligned}
		\label{eq12}
\end{eqnarray}
where
$
\omega^{e} = {\int_{u \geq 0}{g^{e}~d\Xi}} = erfc\left( {- \sqrt{\lambda^{e}}U^{e}} \right)/2 .
$
{
The spatial derivatives of conservative variables $W^e$ and primitive variables $\mathcal{Q}^e$ can be transformed into each other.}
Because the DG method can naturally update conservative variables and their higher-order moments, 
the spatial derivatives  $\frac{\partial^{l}W^{L}}{\partial x^{l}} , \frac{\partial^{l}W^{R}}{\partial x^{l}} $ at the interface in smooth regions can be directly obtained.
This is a notable difference compared to the original CEHGKS method. 
Since the volume integral term only involves the computation of sufficiently smooth macroscopic quantities within the elements, 
inspired by the calculation of the high-order term in the numerical flux, 
it can  directly adopt the   equilibrium state. 

Initially, a reference interval transformation is performed, then the Gaussian rule  is applied to the spatial integral, yielding
\begin{eqnarray}
		\begin{aligned}
{\int_{t_{n}}^{t_{n + 1}}{{\int_{x_{i - \frac{1}{2}}}^{x_{i + \frac{1}{2}}}{F( W_{h} )( {v_{m}(x)} )}_{x}}dx}}dt &= {\int_{0}^{\Delta t}{{\int_{- 1}^{1}{F( {W_{h}( \frac{h}{2}s + x_{i}, t )} )( {P_{m}(s))_s}}}ds}}dt\\
&\approx{\int_{0}^{\Delta t}{\sum\limits_{g = 0}^{K_{g}}{\omega_{g}{F( {W_{h}(\frac{h}{2}{s_x}_g + x_{i}, t )} ){P_{m}^{'}\left( {s_x}_g \right)} }}}}dt\\
& = \sum\limits_{g = 0}^{K_{g}}\omega_{g}{P_{m}^{'}\left( {s_x}_g \right)}{\int_{0}^{\Delta t}{{F\left( {W_{h}\left(x_g, t \right)} \right)}}}dt ,
%
         \end{aligned}
		\label{3-x}
\end{eqnarray}
where $s_{x_g}$ represent the standard Gaussian nodes,
$x_g$ and $\omega_{g}$  represent the Gaussian nodes and weights in the spatial  directions, respectively.
For the time integrals  at the Gauss points $x_g$ within each element are obtained by
\begin{eqnarray}
		\begin{aligned}
{\int_{0}^{\Delta t}{{F\left( {W_{h}\left(x_g, t \right)} \right)}}}dt
\approx{\sum\limits_{\alpha=0}^{K}{{{F\left( {\mathcal{W}^{g}\left( \kappa_{\alpha} {\Delta}t \right)} \right) } \omega_{\alpha} }}} ,
  \end{aligned}
		\label{vx}
\end{eqnarray}
where
\begin{eqnarray}
		\begin{aligned}
\mathcal{W}^{g}(t)=  {\sum\limits_{l = 0}^{r}{\frac{{t}^l}{l!}\frac{\partial^{l}W_{h}\left( {x_{g}, t_n} \right)}{\partial t^{l}}}} ,
         \end{aligned}
		\label{wht}
\end{eqnarray}
here directly adopt the high-order equilibrium state of the numerical flux \eqref{flux1}-\eqref{flux2}.
The temporal derivatives of various orders in  \eqref{wht} are still converted into combinations of spatial derivatives through the Lax-Wendroff procedure in ADER. 
Due to the smoothness of the numerical solution within the element, 
the spatial derivatives of various orders are directly obtained by 

 \begin{equation}\label{wh1dx}
\frac{\partial^{l}W_{h}\left( {x_{g}, t_n} \right)}{\partial x^{l}}=\sum_{m=0}^{k}{w_{i,\hspace{0.1em}m}(t_n)v^{\prime}_m(x_g).}
\end{equation}

 
The above is the computational process for the fully discrete scheme \eqref{eq2-3} in smooth regions. 
In order to suppress oscillations in discontinuous regions, 
it is necessary to use a monotonic limiter. 
Next, we will describe the limiter strategy applied in this paper.

\subsection{The limiter strategy}
Let's assume that   starting from  the time step $n$, that is, 
each order of moments $w_m(t_n), m=0,1,...,k$  are known.
Before advancing  to the next time step $n+1$, 
we need to apply limiter to the "troubled cells" to obtain the $w_m(t_n)^{new}$,
 which are the degrees of freedom after the application of the limiter,
 that is, the new DG moments of each order.
 

 \textbf{(1) Motivation:}
In order to fully leverage the computational advantages of the DG method, 
the limiter strategy aims to minimize the identification of "troubled cells",
while utilizing  a more compact stencil for the reconstruction of these "troubled cells".

  \textbf{(2) Strategy:}
The limiter strategy adopts the idea of combining a detector and a limiter.
Firstly, a {KXRCF\cite{KXRCF}} detector is used to identify "troubled cells".
Subsequently, a {SHWENO\cite{LSY2}} reconstruction technique is applied to limit these "troubled cells" that may contain discontinuities.

Assuming that the "troubled cell" $I_i$ has been identified, 
let us first introduce how to apply limiter to   $I_i$. 
The main idea is to apply SHWENO reconstruction technique  for the "troubled cell" $I_i$ to suppress numerical oscillations. 
Furthermore, SHWENO polynomials are utilized to reconstruct the moments $w_m(t_n)^{new}$ within this cell $I_i$.
We mainly introduce the SHWENO reconstruction under scalar equations.
For simplicity,  in the following, we denote $W_h(x,t)=W(x)$  and omit the index $i$.

The compact stencil of SHWENO is $\{I_{i-1}, I_i, I_{i+1}\}$.
Based on the  cell-average  ${\bar{W}_l}$ and cell-average slope ${\bar{W}_l}^\prime$ of conservative variables, 
{we can define a  quartic polynomial $p_1(x)\in span\left\{1,x,x^2,x^3,x^4\right\}$ and two linear polynomials $p_2(x), p_3(x)\in span\left\{1,x\right\}$,
which satisfy
\begin{equation}
\begin{aligned}
\frac{1}{\Delta x}\int_{I_l}p_1\left(x\right)\mathrm{d}x={\bar{W}_l}, l=i-1,i,i+1,\\
\frac{1}{\Delta x}\int_{I_l}\frac{\mathrm{d}p_1\left(x\right)}{\mathrm{d}x}\mathrm{d}x={\bar{W}_l}^\prime, l=i-1,i+1.
\end{aligned}
\end{equation}
and
\begin{equation}
\begin{aligned}
\frac{1}{\Delta x}\int_{I_l}p_2\left(x\right)\mathrm{d}x={\bar{W}_l}, l=i-1,i,\\
\frac{1}{\Delta x}\int_{I_l}p_3\left(x\right)\mathrm{d}x={\bar{W}_l}, l=i,i+1\\ 
\end{aligned}
\end{equation}
}respectively,
The final SHWENO reconstruction polynomial over the $i$-th cell is given by
\begin{eqnarray}
 p(x)=\alpha_1\left(  \frac{1}{\gamma_1}p_1(x) -  \frac{\gamma_2}{\gamma_1}p_2(x)  -  \frac{\gamma_3}{\gamma_1}p_3(x) \right) + \alpha_2p_2(x) + \alpha_3p_3(x) ,
           \label{shwenop}
\end{eqnarray}
where  $\gamma_1, \gamma_2, \gamma_3$ are the positive linear weights, 
{select $\gamma_1=0.998, \gamma_2=\gamma_3=0.001$,}
and $\alpha_1, \alpha_2, \alpha_3$ are the nonlinear weights.
 For the specific reconstruction steps, see \cite{LSY2}.

In CEHGKS, the calculation of the  cell-average  slopes ${\bar{W}}_l^\prime$  use the equilibrium and non-equilibrium states for time advancing,
which becomes significantly more complex in two-dimensional case.
In contrast, the DG   method allows for the updating of a complete $p^k$   polynomial numerical solution $W(x)$ \eqref{so1d} after time evolution, 
which provides the direct computation of the cell-average ${\bar{W}}_l$ and cell-average slopes ${\bar{W}}_l^\prime$, namely
\begin{eqnarray}
{\bar{W}}_l=&\frac{1}{\Delta x}\int_{I_l}{W\left(x\right)}dx ,\\
{\bar{W}}_l^\prime=&\frac{1}{\Delta x}\int_{I_l}\frac{dW\left(x\right)}{dx}dx , 
   \label{shwenow}
\end{eqnarray}
where $l=i-1, i, i+1$.
This is much simpler and marks a significant difference from the original CEHGKS method.
 
By the idea of projection
\begin{eqnarray}
 \int_{I_i}{W^{new}(x)v_m\left(x\right)}dx=\int_{I_i}{p(x)v_m\left(x\right)}dx, m=0,1...k,
\label{shwenown}
\end{eqnarray}
and the moments 
$w_m^{new}$ satisfy
\begin{eqnarray}
w_m^{new}=\frac{1}{M_m}\int_{I_i}{p(x)v_m\left(x\right)}dx, m=0,1...k,
\label{dgwn} 
\end{eqnarray}
which can be calculated by the Gaussian rule
\begin{eqnarray}
w_m^{new}\approx\frac{\Delta x}{M_m}\sum_{g=1}^{G}\omega_g p\left(x_g\right)v_m\left(x_g\right), 
\label{dgwnG}
\end{eqnarray}
where $x_g, \omega_g$ are the Gaussian nodes and weights.

{In order to better suppress the oscillation, the SHWNO   under the Euler systems is used with a local characteristic field decomposition \cite{MRWENOL}. 

Next, let's revisit how to detect "troubled cells". 
For the detector,   we use the KXRCF indicator \cite{KXRCF},
which is commonly used as a shock  detection for the DG method and can efficiently identify potential discontinuous regions \cite{WL2}.
As shown in \cite{SHWENOL,MRWENOL}, 
if the correlation quantity $\mathcal{K}$ in the cell $I_i$
\begin{eqnarray}
\mathcal{K}=\frac{ |\int_{\partial I_i^-}{(u_i (x)-u_{n_i} (x))ds|}}{h^R |\partial I_i^- |\| u_i (x)\| }>C_k ,
\label{kxrcf}
\end{eqnarray}
where $C_k=1$, then $I_i$ is identified as a "troubled cell". 
Here $\partial I_i^-$ is the inflow boundary (if $\vec{v}\cdot\vec{n}<0$, $\vec{v}$ is the velocity,  and $\vec{n}$ is the outward normal vector of  the boundary $\partial I_i^-$),
$u_i(x)$ is the numerical solution on $I_i$,
and $u_{n_i} (x)$ is the numerical solution on the neighboring cell of $I_i$ in the side $\partial I_i^-$.
$h=\Delta x$,
and $R=1$ for $k=1$, $R=1.5$ for $k>1$. 
{The  $\|u_i(x)\|$ is defined as the maximum value of $|u_i(x)|$.}

In practical applications, when using the standard computational parameters of the KXRCF indicator,
similar to RKDG, an excessive number of "troubled cells" are identified at higher orders, 
particularly evident in fourth- and fifth-order schemes. 
This prevents the high-order schemes from demonstrating significant computational advantages compared to lower-order schemes.

To tackle the above challenge, we have made  improvements to the KXRCF indicator: 
When assessing discontinuities between elements, we attempt to increase the value of the index $\mathcal{K}$ \eqref{kxrcf},
with the aim of reducing the detection range and enhancing the numerical results.
Specifically, by choosing $C_k>1$, which amplifies the difference between the solution inside the element and   at the boundary, thereby increasing the threshold for discontinuity detection.
In our method, we can choose $C_k=5$ in one-dimensional and $C_k=2$ in two-dimensional case.

\textbf{Remark 1:}{
\textit {The feasibility of the adjustment is mainly attributed to
the stability and  robustness demonstrated by the SHWENO reconstruction technique, 
as well as the  remarkable robustness exhibited by the HGKS method \cite{GKS.GRP}.}
}

This adjustment to  the conventional KXRCF indicator \cite{KXRCF}
exhibits high compatible with our DG-HGKS scheme that utilizes the SHWENO reconstruction, 
enabling effective identification of discontinuities and
obtaining a more sharply detection range.
Numerical results will demonstrate this adjustment significantly reduces the number of "troubled cells",
while maintaining the stability and accuracy of the scheme.
\section{The numerical method in two-dimensional case }
\setcounter{equation}{0}
 In this section, we describe the DG-HGKS scheme and the limiter strategy in two-dimensional case.

\subsection{The DG-HGKS scheme in two-dimensional case }
 Considering the 2D Euler  systems
	\begin{eqnarray}
		\begin{aligned}
\frac{\partial W}{\partial t} + \frac{\partial F(W)}{\partial x} + \frac{\partial G(W)}{\partial y} = 0 ,
		\end{aligned}
		\label{eq:euler2d}
	\end{eqnarray}
where the  conservative variables 	
$W = \left( {\rho, \rho U, \rho V, E} \right)^{T}$,
and the  conservative flux
\begin{eqnarray}
	\begin{aligned}
	F(W) &= \left( {\rho U, \rho U^{2} + p, \rho UV, U\left( {E + p} \right)} \right)^{T} ,\\
	G(W) &= \left( {\rho V, \rho UV, \rho V^{2} + p, V(E + p)} \right)^{T} ,
	\end{aligned}
		\label{eqf2d}
\end{eqnarray}
here $U, V$ represent the  velocity in the $x$ and $y$ directions, respectively.
The  energy $E=\frac{p}{\gamma-1}+\frac{1}{2}\rho \left(U^{2}+V^{2}\right)$, 
with $\gamma=1.4$.

Given calculation domain of two-dimensional space $\Omega \triangleq  \left\lbrack {x_{a}, x_{b}} \right\rbrack \times \left\lbrack y_{a}, y_{b} \right\rbrack$.
The uniform cell  $I_{ij}=(x_{i-\frac{1}{2}}, x_{i+\frac{1}{2}})\times(y_{j-\frac{1}{2}}, y_{j+\frac{1}{2}})$,~$ i=1, \cdots, N,~j=1, \cdots, M $,
 where $N, M$ represent the number of cells in the $x$ and $y$ directions respectively,
 the  cell size  $ \Delta x=x_{i+\frac{1}{2}}-x_{i-\frac{1}{2}},~\Delta y=y_{j+\frac{1}{2}}-y_{j-\frac{1}{2}}$.
First introduce the numerical solution of the DG method in 2D, the DG space is 
$$V_h\equiv V_{h}^{\hspace{0.1em}k}=\left\{v(x, y): v|_{I_{ij}}\in P^k\left(I_{ij}\right), i=1, \cdots, N,~j=1, \cdots, M\right\} ,$$
where $P^k\left(I_{ij}\right)$ represents a set of 2D polynomials of degree no more than $k$,
and the numerical solution is 
\begin{equation}\label{only2}
W_h(x, y, t)=\sum_{l=0}^{K}{w_{ij,\hspace{0.1em}l}(t)v_{ij,\hspace{0.1em}l}(x, y)} ,
\end{equation}
here $K=\frac{(k+1)(k+2)}{2}-1$, 
the specific form of the 2D basis functions in fifth-order($p^4$) selected in this paper are as follows:  
$$
v_{ij,\hspace{0.1em}0}(x, y)=1,
v_{ij,\hspace{0.1em}1}(x, y)=v_{i,\hspace{0.1em}1}(x),
v_{ij,\hspace{0.1em}2}(x, y)=v_{j,\hspace{0.1em}1}(y),
$$
$$
v_{ij,\hspace{0.1em}3}(x, y)=v_{i,\hspace{0.1em}2}(x),
v_{ij,\hspace{0.1em}4}(x, y)=v_{i,\hspace{0.1em}1}(x)v_{j,\hspace{0.1em}1}(y),
v_{ij,\hspace{0.1em}5}(x, y)=v_{j,\hspace{0.1em}2}(y),
$$
$$
v_{ij,\hspace{0.1em}6}(x, y)=v_{i,\hspace{0.1em}3}(x),
v_{ij,\hspace{0.1em}7}(x, y)=v_{i,\hspace{0.1em}2}(x)v_{j,\hspace{0.1em}1}(y),
v_{ij,\hspace{0.1em}8}(x, y)=v_{i,\hspace{0.1em}1}(x)v_{j,\hspace{0.1em}2}(y),
$$
\begin{equation}\label{basis2d}
v_{ij,\hspace{0.1em}9}(x, y)=v_{j,\hspace{0.1em}3}(y),
v_{ij,\hspace{0.1em}10}(x, y)=v_{i,\hspace{0.1em}4}(x),
v_{ij,\hspace{0.1em}11}(x, y)=v_{i,\hspace{0.1em}3}(x)v_{j,\hspace{0.1em}1}(y),
\end{equation}
$$
v_{ij,\hspace{0.1em}12}(x, y)=v_{i,\hspace{0.1em}2}(x)v_{j,\hspace{0.1em}2}(y),
v_{ij,\hspace{0.1em}13}(x, y)=v_{i,\hspace{0.1em}1}(x)v_{j,\hspace{0.1em}3}(y),
v_{ij,\hspace{0.1em}14}(x, y)=v_{j,\hspace{0.1em}4}(y) ,
$$
where $v_{i,\hspace{0.1em}l}(x)$ is the basis functions under the one-dimensional DG framework mentioned in \eqref{so1d}. For simplicity, the following omits $i, j$.

Multiply  the Euler systems \eqref{eq:euler2d} by the test function  $v_m(x, y), m=0,1,...K$, and then  integrating on the element  $ I_{ij} \times (t_n, t_{n+1})$, we can obtain the two-dimensional fully discrete DG scheme

\begin{eqnarray}
		\begin{aligned}	
	{\int_{t_{n}}^{t_{n + 1}}{\int{{\int_{I_{ij}}{\left( W_{h} \right)_{t}v_m(x, y)}}dxdy}}}dt &+ {\int_{y_{j - \frac{1}{2}}}^{y_{j + \frac{1}{2}}}{{\int_{t_{n}}^{t_{n + 1}}{F( x_{i + \frac{1}{2}}, y, t )v_m( x_{i + \frac{1}{2}}^{-}, y )}}dtdy}}\\
	& - {\int_{y_{j - \frac{1}{2}}}^{y_{j + \frac{1}{2}}}{{\int_{t_{n}}^{t_{n + 1}}{F( {x_{i - \frac{1}{2}}, y, t})v_m( {x_{i - \frac{1}{2}}^{+}, y} )}}dtdy}} \\
	&+ {\int_{x_{i - \frac{1}{2}}}^{x_{i + \frac{1}{2}}}{{\int_{t_{n}}^{t_{n + 1}}{G( {x, y_{j + \frac{1}{2}}, t} )v_m( {x, y_{j + \frac{1}{2}}^{-}} )}}dtdx}} \\
	&- {\int_{x_{i - \frac{1}{2}}}^{x_{i + \frac{1}{2}}}{{\int_{t_{n}}^{t_{n + 1}}{G( x, y_{j - \frac{1}{2}}, t )v_m( x, y_{j - \frac{1}{2}}^{+} )}}dtdx}} \\
	&- {\int_{t_{n}}^{t_{n + 1}}{\int{{\int_{I_{ij}}{F( W_{h} )( {v_m(x, y)})}_{x}}dxdy}}}dt\\
	& - {\int_{t_{n}}^{t_{n + 1}}{\int{{\int_{I_{ij}}{G( W_{h} )( {v_m(x, y)} )}_{y}}dxdy}}}dt = 0 .
\end{aligned}
		\label{fully2d}
	\end{eqnarray}
Transforming the cell $I_{ij}$ into the reference interval, 
and using the orthogonal property of Legendre polynomials, we obtain
	\begin{eqnarray}
		\begin{aligned}
{\int_{t_{n}}^{t_{n + 1}}{\int{{\int_{I_{ij}}{\left( W_{h} \right)_{t}v_m(x, y)}}dxdy}}}dt 
             &={\int_{t_{n}}^{t_{n + 1}}\!{{\int{\int_{I_{ij}}({\sum_{l=0}^{K}{w_{l}(t)v_{l}(x, y)}})_t v_{m}(x, y) dxdy}}}}dt \\
&=\frac{\Delta x\Delta y}{M_{m}}\left( w_m(t_{n+1})-w_m(t_n)\right), m=0,...,K.
		\label{2d1}
		\end{aligned}
\end{eqnarray}
The 2D diagonal mass matrix is
$$
M=diag(1, \frac{1}{3}, \frac{1}{3}, \frac{1}{5}, \frac{1}{9}, \frac{1}{5}, \frac{1}{7}, \frac{1}{15}, \frac{1}{15}, \frac{1}{7}, \frac{1}{9}, \frac{1}{21}, \frac{1}{25}, \frac{1}{21}, \frac{1}{9}) ,
$$
where $M_m$ is the $m$-th diagonal element.

As for the numerical flux at the cell boundary $x_{i+\frac{1}{2}}$, based on the form of basis function \eqref{basis2d}, let $v_m(x, y)=v_l(x)·v_l(y)$, we have
\begin{eqnarray}
		\begin{aligned}
{\int_{y_{j - \frac{1}{2}}}^{y_{j + \frac{1}{2}}}{{\int_{t_{n}}^{t_{n + 1}}{F( x_{i + \frac{1}{2}}, y, t )v_m( x_{i + \frac{1}{2}}^{-}, y)}}dydt}}= &v_{l}( x_{i + \frac{1}{2}}^{-} )\mathbb{F}_{i + \frac{1}{2}, j}\left( {{\Delta}t} \right) ,
             \end{aligned}
		\label{2d2}
\end{eqnarray}
the time integral of numerical flux  is evaluated by
\begin{eqnarray}
		\begin{aligned}
   \mathbb{F}_{i + \frac{1}{2}, j}\left( {{\Delta}t} \right) = &{\int_{y_{j - \frac{1}{2}}}^{y_{j + \frac{1}{2}}}{{\int_{t_{n}}^{t_{n + 1}}{F( x_{i + \frac{1}{2}}, y, t )v_{l}( y )}}dydt}} \\
   =& \left( 1 - e^{- {\Delta}t/\tau} \right)\mathbb{F}_{i + \frac{1}{2}, j}^{e}\left( {{\Delta}t} \right) 
   + e^{- {\Delta}t/\tau}\mathbb{F}_{i + \frac{1}{2}, j}^{k}\left( {{\Delta}t} \right) ,
         \end{aligned}
		\label{2DFT}
\end{eqnarray}
where
$\mathbb{F}_{i + \frac{1}{2}, j}^{k}\left({{\Delta}t}\right)$ depends on the time integral of the second-order KFVS flux.
The high-order term is approximated by the Gaussian rule  
\begin{eqnarray} \label{2DFlux-gauss}
\mathbb{F}_{i+\frac12, j}^e(\Delta t)\thickapprox\sum_{\alpha=0}^{\hat{K}}\sum_{\beta=0}^{\hat{K}}\mathcal{F}_{Eu}\left(\mathcal{Q}^e\left(y_{j_\alpha}, \kappa_\beta\Delta t\right)\right)v_{l}(y_{j_\alpha})\omega_\alpha\omega_\beta,
\end{eqnarray}
here 
\begin{equation} \label{qe}
    \begin{aligned}    
    \mathcal{Q}^e\left(y, t\right)=\sum_{l=0}^r\frac1{l!}\left(t\frac\partial{\partial t}+y\frac\partial{\partial y}\right)^l\mathbf{Q}_{i+\frac12}^e(y_j, t_n) ,
    \end{aligned}
\end{equation}
where $\mathbf{Q}=(\rho, U, V, p)^\mathrm{T}$.
The fifth-order DG scheme takes $r=4, \hat{K}=3$. 
{$y_{j_\alpha}, \kappa_\beta, \omega_\alpha, \omega_\beta$ is the Gaussian points  and weights,}  
which can be obtained after  the reference interval transformation.
%

As for the   volume integral,
firstly substituting the basis function and applying the Gaussian rule in the $x$-spatial integral, we can obtain 
\begin{eqnarray}
\begin{aligned}
{\int_{t_{n}}^{t_{n + 1}}{\int{{\int_{I_{ij}}{F\left( W_{h} \right)\left( {v_m(x, y)} \right)}_{x}}dxdy}}}dt
=&{\int_{t_{n}}^{t_{n + 1}}\!{{\int_{x_{i-\frac{1}{2}}}^{x_{i+\frac{1}{2}}}{\int_{y_{j-\frac{1}{2}}}^{y_{j+\frac{1}{2}}}{F\left( {W_{h}} \right)\left( {v_l(x)}_{x}{v_l(y)} \right)}}}dxdy}}dt 
\\
\approx&\sum_{\iota=0}^{\hat{K}}\omega_\iota v_l^{'}(x_{\iota})\mathbb{F}_{i_\iota, j}(\Delta t) ,	
\end{aligned}
	\label{2DFX}
\end{eqnarray}
where $\omega_{\iota},~x_{\iota}$
is the Gaussian weights and nodes  in $x$ direction. 
$\mathbb{F}_{i_\iota, j}(\Delta t)$ represents the time integral  term at  each point  $x_\iota$ within elements,
the specific form is 
\begin{eqnarray}
\mathbb{F}_{i_\iota, j}(\Delta t)
={\int_{y_{j - \frac{1}{2}}}^{y_{j + \frac{1}{2}}}{{\int_{t_{n}}^{t_{n + 1}}{F\left(W_h(x_{\iota}, y, t)\right) v_l(y)}}dydt}} ,
\end{eqnarray}
following as the high-order terms of the numerical flux \eqref{2DFlux-gauss},
it can be calculated by 
\begin{eqnarray} \label{2Dvolume-gauss}
\mathbb{F}_{i_\iota, j}(\Delta t)
\thickapprox\sum_{\alpha=0}^{\hat{K}}\sum_{\beta=0}^{\hat{K}}\mathcal{F}\left(\mathcal{W}^{\iota}\left(y_{j_\alpha},\kappa_\beta\Delta t\right)\right)v_{l}(y_{j_\alpha})\omega_\alpha\omega_\beta,
\end{eqnarray}
here 
\begin{equation} \label{qe}
    \begin{aligned}    
    \mathcal{W}^{\iota}\left(y, t\right)=\sum_{l=0}^r\frac1{l!}\left(t\frac\partial{\partial t}+y\frac\partial{\partial y}\right)^l{W}_{h}(x_{\iota}, y_j, t_n) .
    \end{aligned}
\end{equation}
The time related derivatives of $W_h$ still be converted into combinations of spatial derivatives using the Lax-Wendroff  in ADER.
The conservative variables $W_h$ and its high-order spatial derivatives can be directly obtained using numerical solution \eqref{only2}.

The calculation in the spatial $y$ direction can be obtained in a similar way after the local coordinate transformation.
For more details  refer to \cite{LSY,LSY2}.

\subsection{The limiter strategy in two-dimensional case}

For the limiter strategy in  two-dimensional case, we continue to adopt the combination of a detector and a limiter.
The detector remains the KXRCF indicator, with specific details provided in \cite{MRWENOL}.
Additionally, the improvement  for the KXRCF indicator discussed in the one-dimensional  is   applicable to the two-dimensional case.
As for the limiter, we utilize a dimension-by-dimension  SHWENO reconstruction technique.
Assuming that a "troubled cell" $I_{ij}$ has been identified, the following describes how to apply the SHWENO reconstruction in the cell $I_{ij}$ briefly.
For the specific  process, please refer to the relevant literature \cite{LSY2,LDM} for details.

Similar to the one-dimensional case, 
based on the idea of projection 
\begin{eqnarray}
 \int_{I_{ij}}{W^{new}(x, y)v_m\left(x, y\right)}dxdy=\int_{I_{ij}}{\hat{w}(x, y)v_m\left(x, y\right)}dxdy\ ,m=0,1...K,
\label{shwenown2D}
\end{eqnarray}
the moments of the new polynomial $W^{new}(x, y)$ satisfy
 \begin{eqnarray}
w_m^{new}=\frac{\Delta x\Delta y}{M_{\text{diag}}}\sum_{gx=1}^{G}{\sum_{gy=1}^{G}{\omega_{gx}\omega_{gy}}\hat{w}\left(x_{gx},y_{gy}\right)v_m\left(x_{gx}, y_{gy}\right)},m=0,1,...K. 
\label{dgwnG2D}
\end{eqnarray}
where $x_{gx},  y_{gy}$ are Gaussian  points  in the $x$ and $y$ directions, respectively, 
and $\omega_{gx}, \omega_{gy}$ are Gaussian weights.
{The value at each Gauss point $\hat{w}\left(x_{gx}, y_{gy}\right)$ is obtained by a dimension-by-dimension  SHWENO reconstruction technique.
Let's give a brief introduction.

Define the cell-average   ${\bar{W}}_{i,j}$,
and cell-average slopes  ${\bar{W}}_{i,j}^{(x)}, {\bar{W}}_{i,j}^{(y)}, {\bar{W}}_{i,j}^{(xy)}$,
yielding
 \begin{eqnarray}
{\bar{W}}_{i,j}=&\frac{1}{\Delta x \Delta y}\int_{I_{i,j}}{W\left(x, y\right)}dxdy,  ~~
{\bar{W}}_{i,j}^{(xy)}=&\frac{1}{\Delta x \Delta y}\int_{I_{i,j}}\frac{\partial W\left(x, y\right)}{\partial xy}dxdy,\\
{\bar{W}}_{i,j}^{(x)}=&\frac{1}{\Delta x \Delta y}\int_{I_{i,j}}\frac{\partial W\left(x, y\right)}{\partial x}dxdy,   ~~
{\bar{W}}_{i,j}^{(y)}=&\frac{1}{\Delta x\Delta y}\int_{I_{i,j}}\frac{\partial W\left(x, y\right)}{\partial y}dxdy.
   \label{shwenow}
\end{eqnarray}
First, we perform two $y$-direction reconstructions
\begin{eqnarray}
\begin{aligned}
\{\bar{W}_{mn},\bar{W}^{(y)}_{mn}\}\to\bar{W}_{i+l,j}(y_{g_{y}})\approx\frac{1}{\Delta x}\int_{I_{i+l},j}W(x,y_{g_{y}})dx, \\
\{\bar{W}^{(x)}_{mn},\bar{W}^{(xy)}_{mn}\}\to\bar{W}^{(x)}_{i+l,j}(y_{g_{y}})\approx\frac{1}{\Delta x}\int_{I_{i+l},j}W_{x}(x,y_{g_{y}})dx. 
\end{aligned}
\end{eqnarray}
Then we use $\bar{W}_{i+l,j}(y_{g_y})$ and $\bar{W}^{(x)}_{i+l,j}(y_{g_y})$ to perform $x$-direction reconstruction to get an approximation to $W(x_{g_x}, y_{g_y})$, i.e.,
\[
\{\bar{W}_{mn}(y_{g_y}), \bar{W}_{x,mn}(y_{g_y})\} \rightarrow \hat{w}(x_{g_x}, y_{g_y}).
\]
To maintain symmetry, the same approach is first applied for reconstruction in the $x$-direction and then in the $y$-direction. 
The final result is obtained by taking the average of the two.

}

It is important to emphasize that in the CEHGKS \cite{LSY2}, 
the cell-average slopes are advanced in time using both non-equilibrium and equilibrium states, 
this approach results in higher computational complexity, particularly in two-dimensional case,
as it specifically requires calculations related to six spatial Gauss points within each cell.
However, in  DG method   can naturally update the  cell-average slopes
by directly calculating from the numerical solution $W(x, y)$ \eqref{only2}.
This is   very different and simpler compared to the CEHGKS.



\section{Numerical examples}
\setcounter{equation}{0}
In this section, numerical examples in both one- and two-dimensional are given to test the accuracy,
effectiveness, and robustness of the DG-HGKS scheme.
Specifically,
Example 4.1 tests the 1D accuracy, 
while Examples 4.6 and 4.7 focus on 2D accuracy.
To verify the ability of the DG-HGKS scheme  in capturing discontinuities, 
Examples 4.2 through 4.5 offer 1D   numerical results, 
and Examples 4.8 and 4.9 provide corresponding 2D  numerical results.

The time step $\Delta t$ is determined by the CFL condition
$$\Delta t=CFL\times {\rm{min}}\{ \frac{\Delta x}{\lvert U \rvert+c }, \frac{\Delta y}{\lvert V \rvert+c}\} ,$$
where  $c=\sqrt{\gamma RT}$ is the speed of sound. 

\textbf{Remark 2:}{\textit{
The $CFL$ condition numbers required for the one-stage DG-HGKS scheme to achieve the corresponding orders without  limiter are lower than those for the RKDG method, 
with specific values of 0.26, 0.13, 0.08, 0.06 for the second-order($p^1$),  third-order($p^2$),  fourth-order($p^3$), and fifth-order($p^4$), respectively. 
When the limiter is applied to all cells, the $CFL$ numbers  can be raised to the level used in RKDG method.
The specific numerical results are provided in the Appendix A.
}}

The numerical implementation of boundary conditions is handled through ghost cells, with specific treatments applied to different boundary types. 
Periodic boundary conditions mimic a periodically repeated domain, with ghost cell values taken from the interior cells near the opposite boundary. 
Inflow/outflow boundary conditions use zero-gradient extrapolation for density, velocity and pressure, with more sophisticated schemes for high-order accuracy or non-steady conditions. 
Non-reflecting boundary conditions minimize numerical reflections using characteristic-based methods to extrapolate outgoing waves. 
Reflective boundary conditions (solid wall) set the normal velocity component in ghost cells to the negative of the adjacent interior cell (no-slip or slip condition), while tangential velocity components, density, and pressure are set equal to their interior values.

\subsection{One-dimensional examples}

{\bf Example 4.1}{ 1D linear advection of the density perturbation}

Here we test the accuracy  of 1D DG-HGKS scheme when the solution is smooth.
The initial condition is given by
\begin{eqnarray}
			\left( {\rho, U, p } \right) = \left( {1 + 0.2{\sin{\pi x}}, 0.7, 1} \right),
 \label{1dorder}
\end{eqnarray}
under the periodic boundary condition, the analytic solution is
$$
			\left( {\rho, U, p} \right) = \left( {1 + 0.2{\sin{\pi (x-0.7t)}}, 0.7, 1} \right).
		$$
The computational domain  $[ {0,2 } ]$ is divided into $N$ uniform cells.
The output time is $T=1.0$.
The numerical errors and   orders  are listed in the {Table \ref {tab:ac1d} whether or not a limiter is applied},
which demonstrate that  the DG-HGKS scheme  can achieve the  second-order($p^1$), third-order($p^2$), fourth-order($p^3$), and fifth-order($p^4$) in both space and time.

 \begin{table}[h]
\centering
\begin{tabular}{|c|c|cccc|cccc|}
\hline
 &   & \multicolumn{4}{c|}{ \textbf{DG without limiter}} & \multicolumn{4}{c|} {\textbf{DG with limiter}} \\ \hline
 &{ { N }}   & { {$L^1$ error}} & { {order}} & { {$L^\infty$ error}} & { {order}} & { {$L^1$ error}} & { {order}} & { {$L^\infty$ error }}& { {order}} \\ \hline
\multirow{4}{*}{{ {$\boldsymbol{p^1}$}}} 
& 20   & 1.64e-03 & --      & 4.13e-03 & --     & 9.74e-04 & --      & 3.46e-03 & --\\
& 40   & 4.22e-04 & 1.96 & 1.07e-03 & 1.95 & 2.30e-04 & 2.08 & 9.00e-04 & 1.94 \\
& 80   & 1.06e-04 & 1.99 & 2.77e-04 & 1.95 & 5.61e-05 & 2.03 & 2.28e-04 & 1.98\\
& 160 & 2.66e-05 & 1.99 & 6.95e-05 & 1.99 & 1.40e-05 & 2.00 & 5.74e-05 & 1.99\\ \hline
\multirow{4}{*}{{ {$\boldsymbol{p^2}$}}} 
& 20   & 2.73e-05 & --      & 1.31e-05 & --     & 3.89e-05 & --      & 3.07e-05 & --\\
& 40   & 3.38e-06 & 3.01 & 1.65e-06 & 2.99 & 4.10e-06 & 3.25 & 3.29e-06 & 3.22 \\
& 80   & 4.23e-07 & 3.00 & 2.12e-07 & 2.95 & 4.60e-07 & 3.16 & 3.70e-07 & 3.15\\
& 160 & 5.26e-08 & 3.01 & 2.64e-08 & 3.01 & 5.42e-08 & 3.08 & 4.38e-08 & 3.08\\ \hline
\multirow{4}{*}{{ $\boldsymbol{p^3}$}} 
& 20   & 8.32e-07 & --      & 6.32e-07 & --     & 7.03e-06 & --      & 6.25e-06 & --\\
& 40   & 5.30e-08 & 3.97 & 3.90e-08 & 4.02 & 3.21e-07 & 4.45 & 2.84e-07 & 4.46 \\
& 80   & 3.33e-09 & 3.99 & 2.41e-09 & 4.02 & 1.62e-08 & 4.31 & 1.39e-08 & 4.35\\
& 160 & 2.08e-10 & 4.00 & 1.50e-10 & 4.01 & 8.90e-10 & 4.18 & 7.48e-10 & 4.22\\ \hline
\multirow{4}{*}{{ $\boldsymbol{p^4}$}} 
& 20   & 4.31e-09 & --      & 3.95e-09 & --     & 4.56e-06 & --      & 3.71e-06 & --\\
& 40   & 1.36e-10 & 4.98 & 1.34e-10 & 4.87 & 1.30e-07 & 5.14 & 1.16e-07 & 5.00 \\
& 80   & 4.45e-12 & 4.93 & 4.73e-12 & 4.87 & 4.05e-09 & 5.00 & 3.63e-09 & 4.99\\
& 160 & 1.39e-13 & 4.99 & 1.58e-13 & 4.90 & 1.28e-10 & 4.98 & 1.13e-10 & 5.00\\ \hline
\end{tabular}
\caption{{\label{tab:ac1d} Accuracy test  for 1D of Example 4.1 \eqref{1dorder}.}}
\end{table}

Next, we will verify the ability of the DG-HGKS scheme to capture discontinuities.
In all 1D examples,  
numerical   results  for the second-order($p^1$), third-order($p^2$), fourth-order($p^3$), and fifth-order($p^4$) schemes are represented by squares of different colors,  
whereas the exact or reference solution is represented by the black line.

{\bf Example 4.2} { 1D Riemann problems}

The typical  1D Riemann problems including the Sod problem and Lax problem.
The initial condition of the Sod problem is given by
\begin{eqnarray}
	\begin{aligned}
\left( {\rho, U, p} \right) = \left\{ \begin{matrix}
{\left( {1, 0, 1} \right) ,~ - 5 \leq x \leq 0,} \\
{\left( {0.125, 0, 0.1} \right) ,~0 < x \leq 5.}
\end{matrix} \right.
\end{aligned}
	\label{sod}
\end{eqnarray}
The Lax problem with the initial condition
\begin{eqnarray}
	\begin{aligned}
\left( {\rho, U, p} \right) = \left\{ \begin{matrix}
{\left( {0.445, 0.698, 3.528} \right) ,~ - 5 \leq x \leq 0,} \\
{\left( {0.5, 0, 0.571} \right) ,~0 < x \leq 5,}
\end{matrix} \right.
\end{aligned}
	\label{lax}
\end{eqnarray}
with the inflow/outflow boundary condition.
200 uniform cells are used in the case. 
The output time is $T=2.0$ for the Sod problem and $T=1.3$ for the Lax problem.

The density distributions for each order are given   by Fig. \ref{fig:SOD} and Fig. \ref{fig:LAX},  
and the time history of the troubled cells are shown in Fig. \ref{fig:SODKX} and Fig. \ref{fig:LAXKX}.
The numerical results show that the DG-HGKS scheme constructed in this paper has a good ability to solve the 1D Riemann problems, and the KXRCF indicator can  capture the troubled cells accurately.

\begin{figure}[h!]
	\centering
	\subfloat{
	\begin{minipage}[t]{0.5\textwidth}  
	\centering
	 \includegraphics[width=\textwidth]{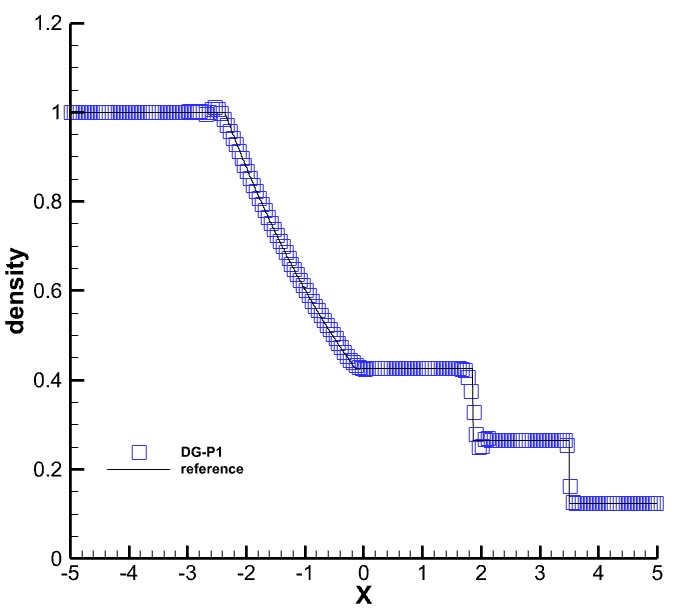}
  \end{minipage}
   \begin{minipage}[t]{0.5\textwidth}  
	\centering
	 \includegraphics[width=\textwidth]{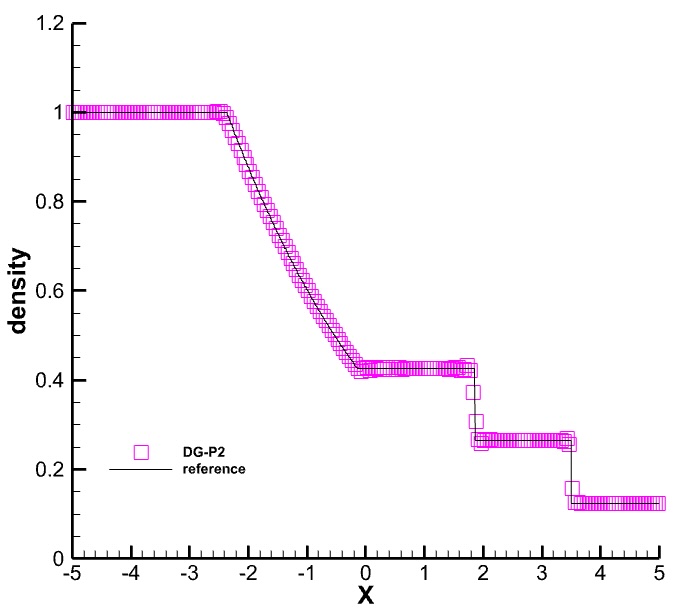}
  \end{minipage}
 }\\
   \subfloat{
  \begin{minipage}[t]{0.5\textwidth}  
	\centering
	 \includegraphics[width=\textwidth]{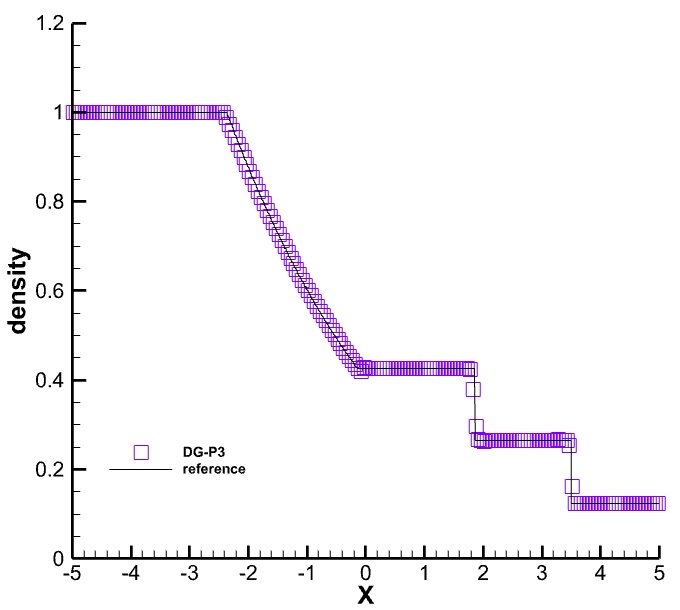}
  \end{minipage}
  \begin{minipage}[t]{0.5\textwidth}  
	\centering
	 \includegraphics[width=\textwidth]{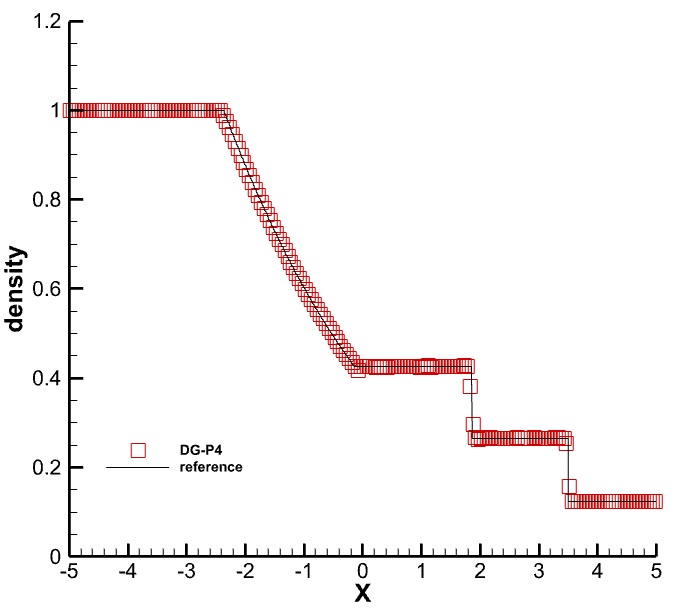}
  \end{minipage}
  }
  \caption{The density distribution of Example 4.2 \eqref{sod} with $N=200.$
  Left top: second-order($p^1$), right top: third-order($p^2$), left bottom: fourth-order($p^3$), right bottom: fifth-order($p^4$).}
  \label{fig:SOD}
\end{figure}
\begin{figure}[h!]
	\centering
	\subfloat{
	\begin{minipage}[t]{0.24\textwidth}  
	\centering
	 \includegraphics[width=\textwidth]{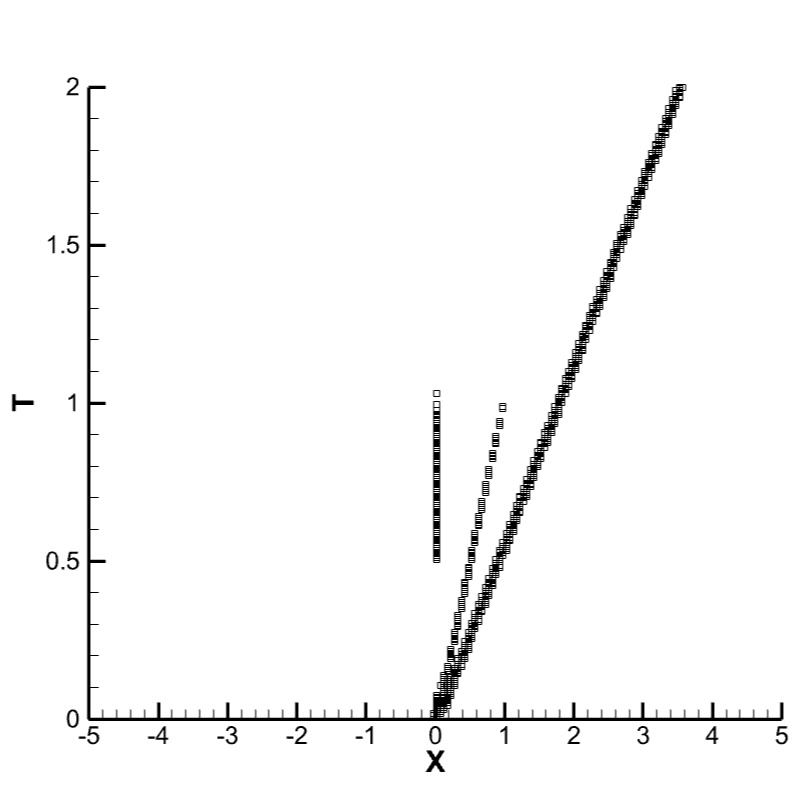}
  \end{minipage}
   \begin{minipage}[t]{0.24\textwidth}  
	\centering
	 \includegraphics[width=\textwidth]{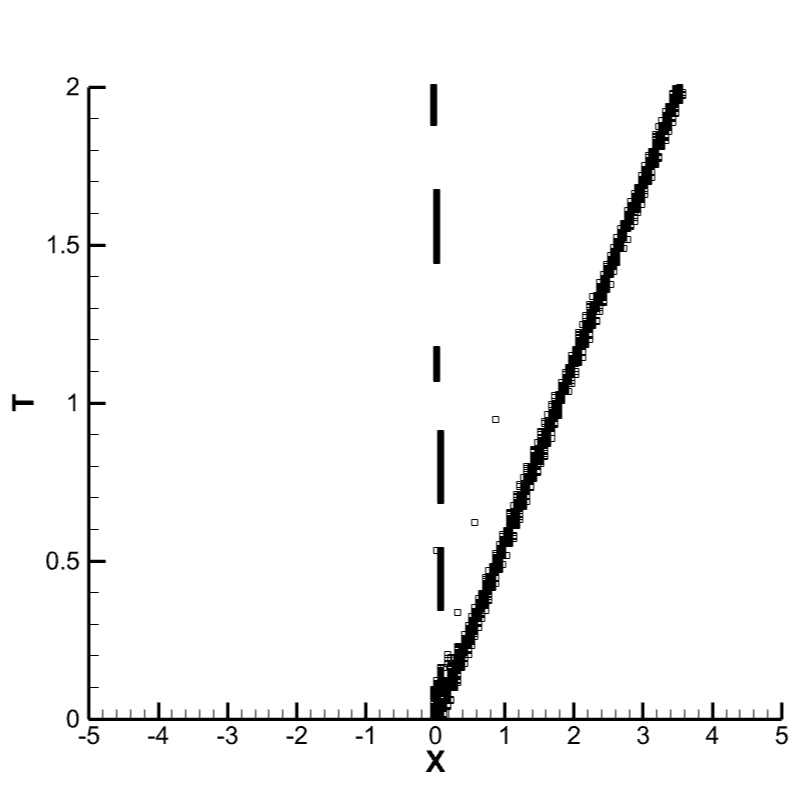}
  \end{minipage}
  \begin{minipage}[t]{0.24\textwidth}  
	\centering
	 \includegraphics[width=\textwidth]{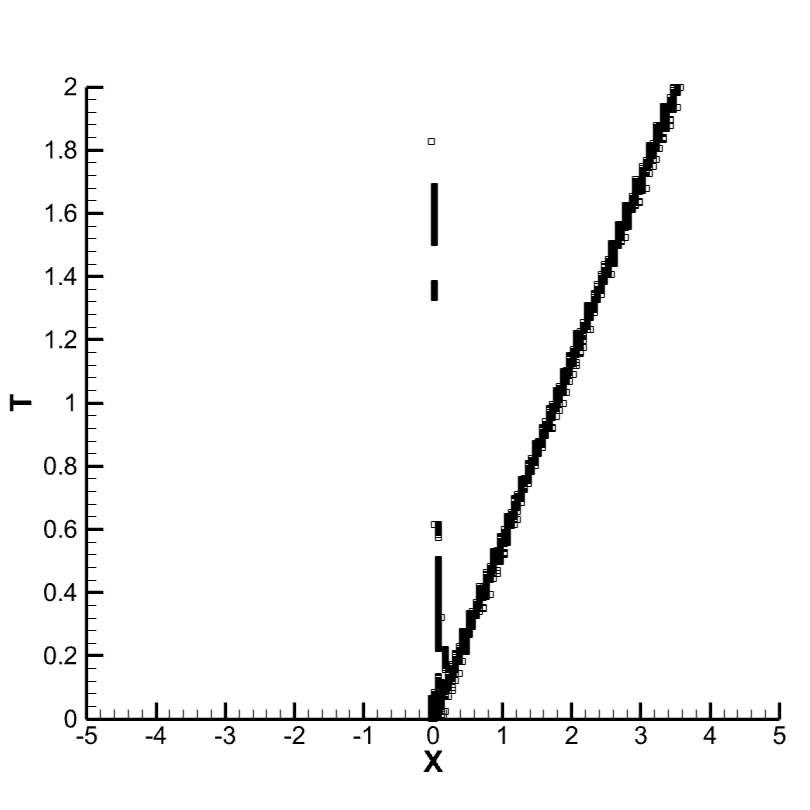}
  \end{minipage}
  \begin{minipage}[t]{0.24\textwidth}  
	\centering
	 \includegraphics[width=\textwidth]{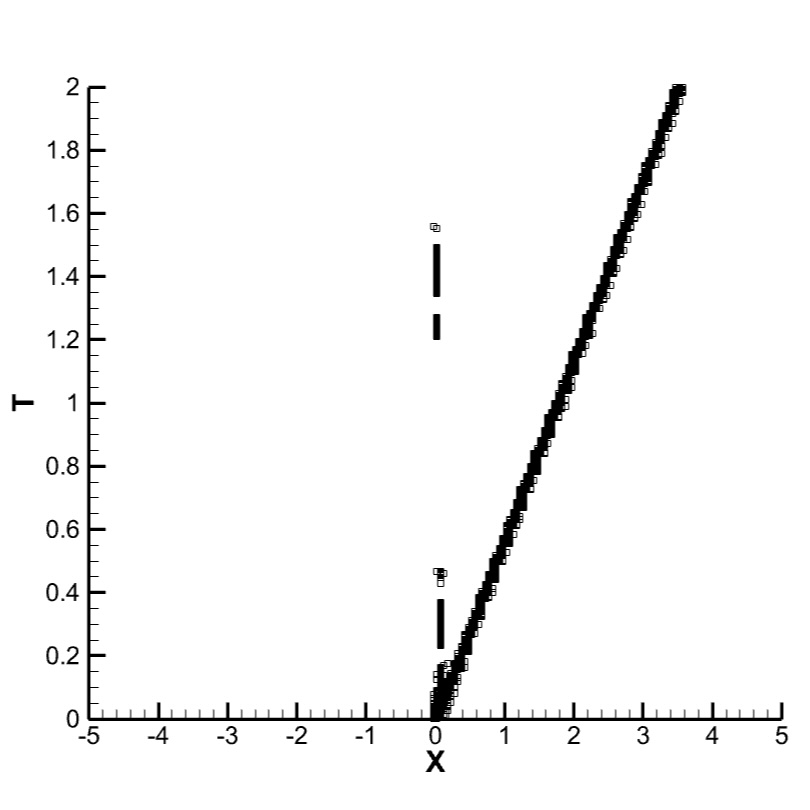}
  \end{minipage}
}
  \caption{The time history of the troubled cells for Example 4.2 \eqref{sod}   with  $N=200.$
  DG-HGKS with SHWENO limiter. 
  From left to right: second-order($p^1$), third-order($p^2$), fourth-order($p^3$), fifth-order($p^4$).}
  \label{fig:SODKX}
\end{figure}

\begin{figure}[h!]
	\centering
	\subfloat{
	\begin{minipage}[t]{0.5\textwidth}  
	\centering
	 \includegraphics[width=\textwidth]{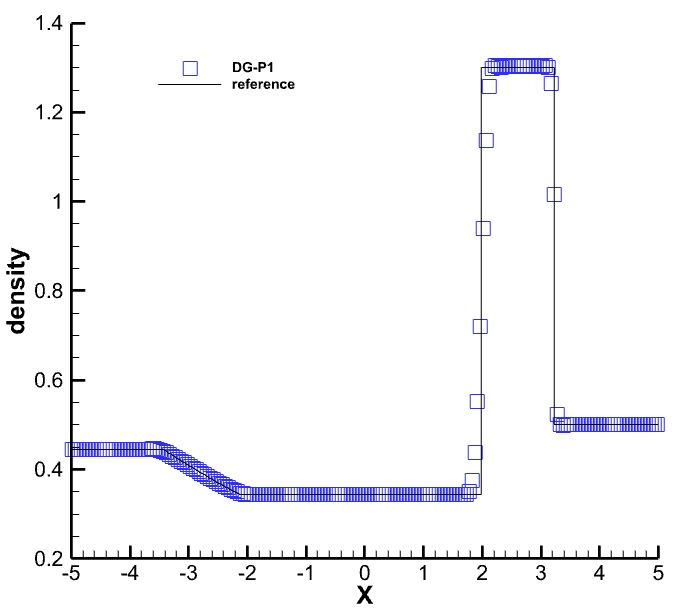}
  \end{minipage}
   \begin{minipage}[t]{0.5\textwidth}  
	\centering
	 \includegraphics[width=\textwidth]{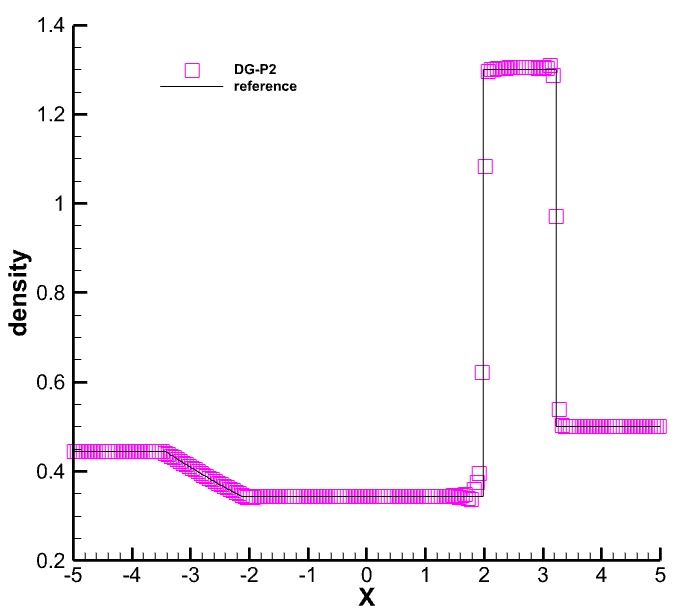}
  \end{minipage}
 }\\
  \subfloat{
  \begin{minipage}[t]{0.5\textwidth}  
	\centering
	 \includegraphics[width=\textwidth]{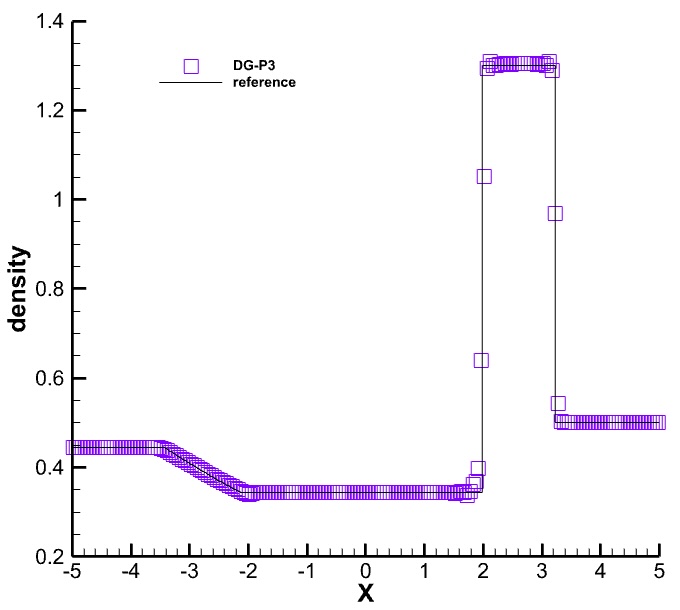}
  \end{minipage}
  \begin{minipage}[t]{0.5\textwidth}  
	\centering
	 \includegraphics[width=\textwidth]{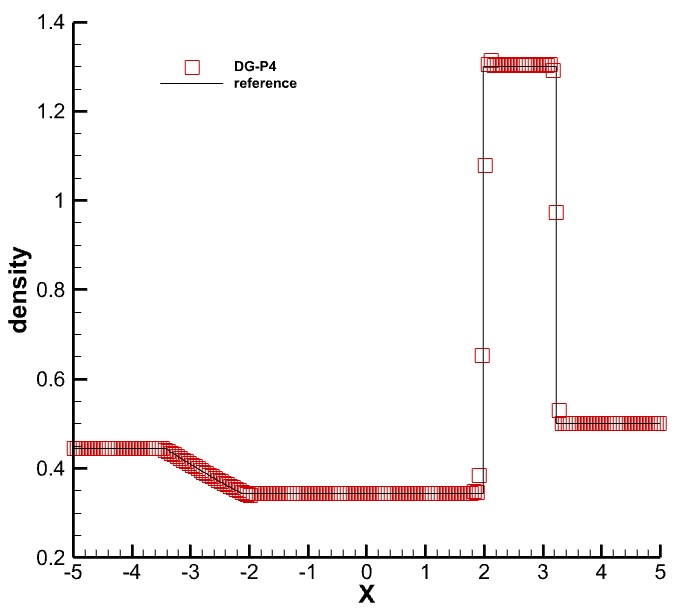}
  \end{minipage}
  }
  \caption{The density distribution of  Example 4.2 \eqref{lax}  with  $N=200.$
  Left top: second-order($p^1$), right top: third-order($p^2$), left bottom: fourth-order($p^3$), right bottom: fifth-order($p^4$).}
  \label{fig:LAX}
\end{figure}

\begin{figure}[h!]
	\centering
	\subfloat{
	\begin{minipage}[t]{0.24\textwidth}  
	\centering
	 \includegraphics[width=\textwidth]{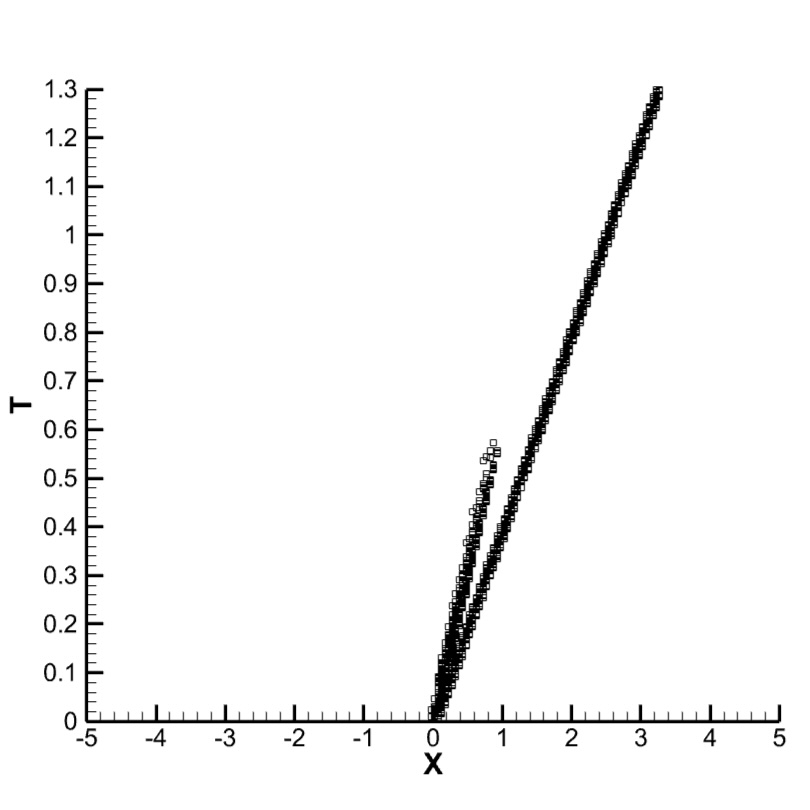}
  \end{minipage}
   \begin{minipage}[t]{0.24\textwidth}  
	\centering
	 \includegraphics[width=\textwidth]{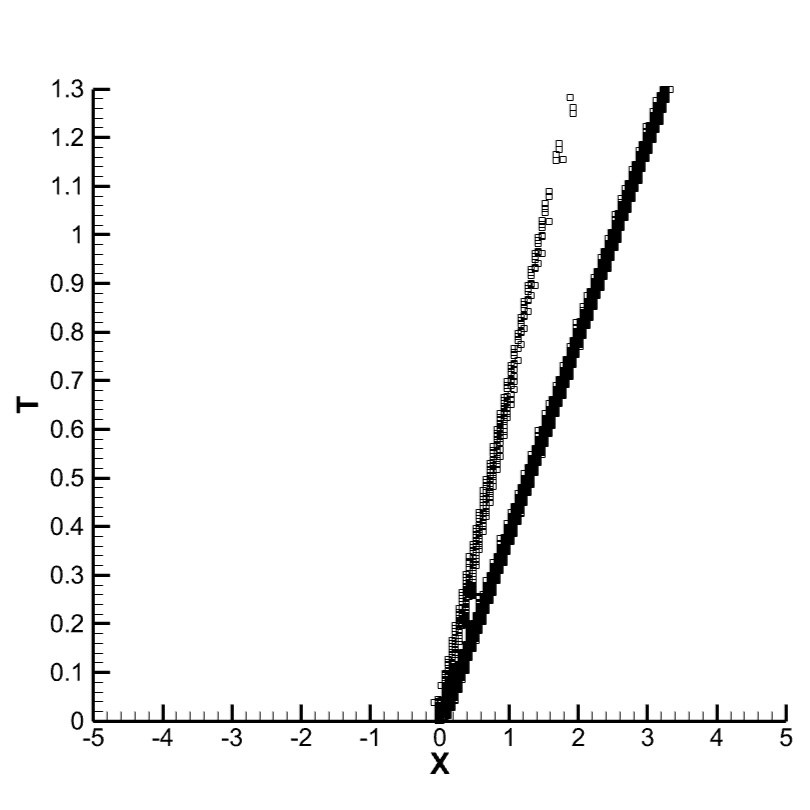}
  \end{minipage}
  \begin{minipage}[t]{0.24\textwidth}  
	\centering
	 \includegraphics[width=\textwidth]{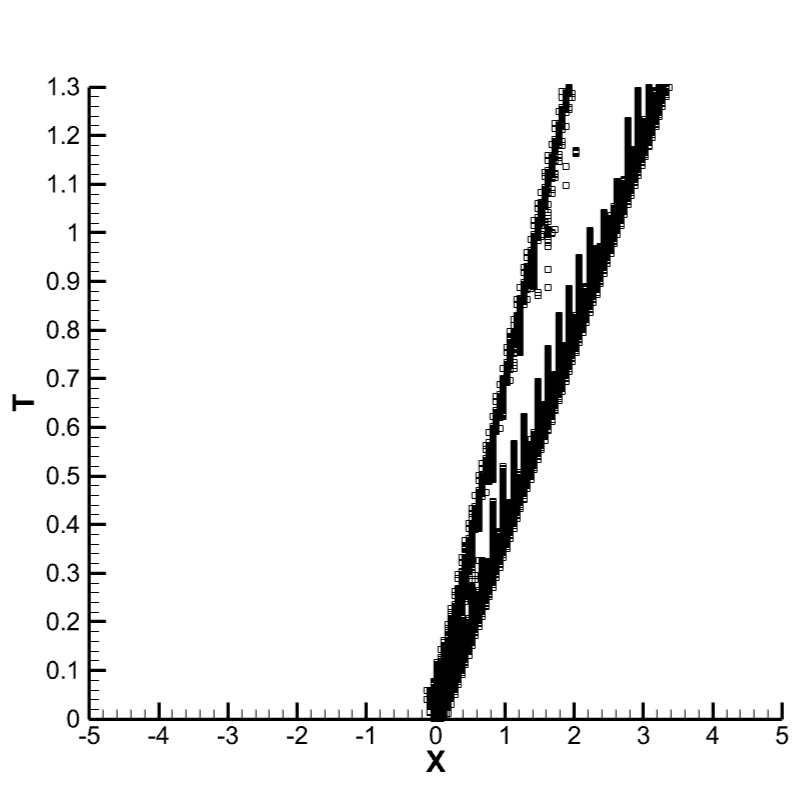}
  \end{minipage}
  \begin{minipage}[t]{0.24\textwidth}  
	\centering
	 \includegraphics[width=\textwidth]{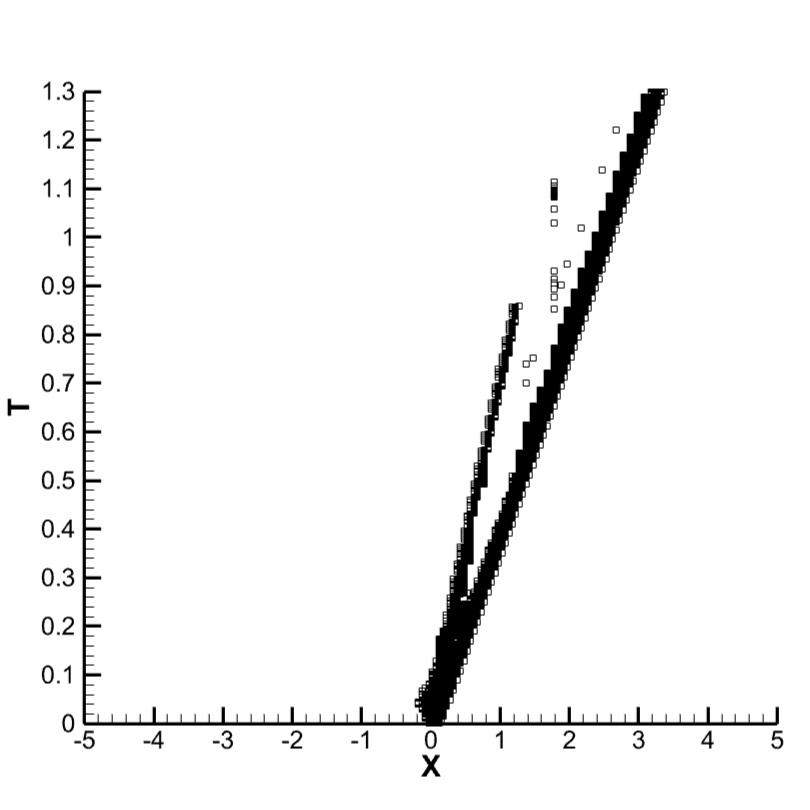}
  \end{minipage}
}
  \caption{The time history of the troubled cells for Example 4.2 \eqref{lax}   with  $N=200.$
  DG-HGKS with SHWENO limiter. 
  From left to right: second-order($p^1$), third-order($p^2$), fourth-order($p^3$), fifth-order($p^4$).}
  \label{fig:LAXKX}
\end{figure}

{\bf Example 4.3} { {Shu-Osher} shock acoustic wave interaction}

The Shu-Osher  problem  is considered in this example,
which contains both small-scale perturbations and shock waves \cite{ADER}.
This case is to assess whether the high-order schemes can capture the small-scale information of the flows exactly.
The initial condition is given by

\begin{eqnarray}
	\begin{aligned}
\left( {\rho, U, p} \right) = \left\{ \begin{matrix}
{\left( {3.857134, 2.629369, 10.33333} \right) ,~x \leq - 4.0,} \\
{\left( {\left. 1 + 0.2{\sin\left( 5\pi x \right.} \right), 0, 1} \right) ,~x > - 4.0,}
\end{matrix} \right.
\end{aligned}
	\label{so}
\end{eqnarray}
with the inflow/outflow boundary condition.
The computational domain is $[-5, 5]$, 
the output time is $T=1.8$.
Fig. \ref{fig:SO} shows the density distributions   under the cells $N=200$,
and the reference solution obtained by a fifth-order finite volume WENO scheme  with 10,000 uniform points.
Fig. \ref{fig:SOKX} shows the time history of the troubled cells.
The results show that the higher-order schemes exhibits better convergence than the lower-order schemes.
 
\begin{figure}[h!]
	\centering
	\subfloat{
	\begin{minipage}[t]{0.5\textwidth}  
	\centering
	 \includegraphics[width=\textwidth]{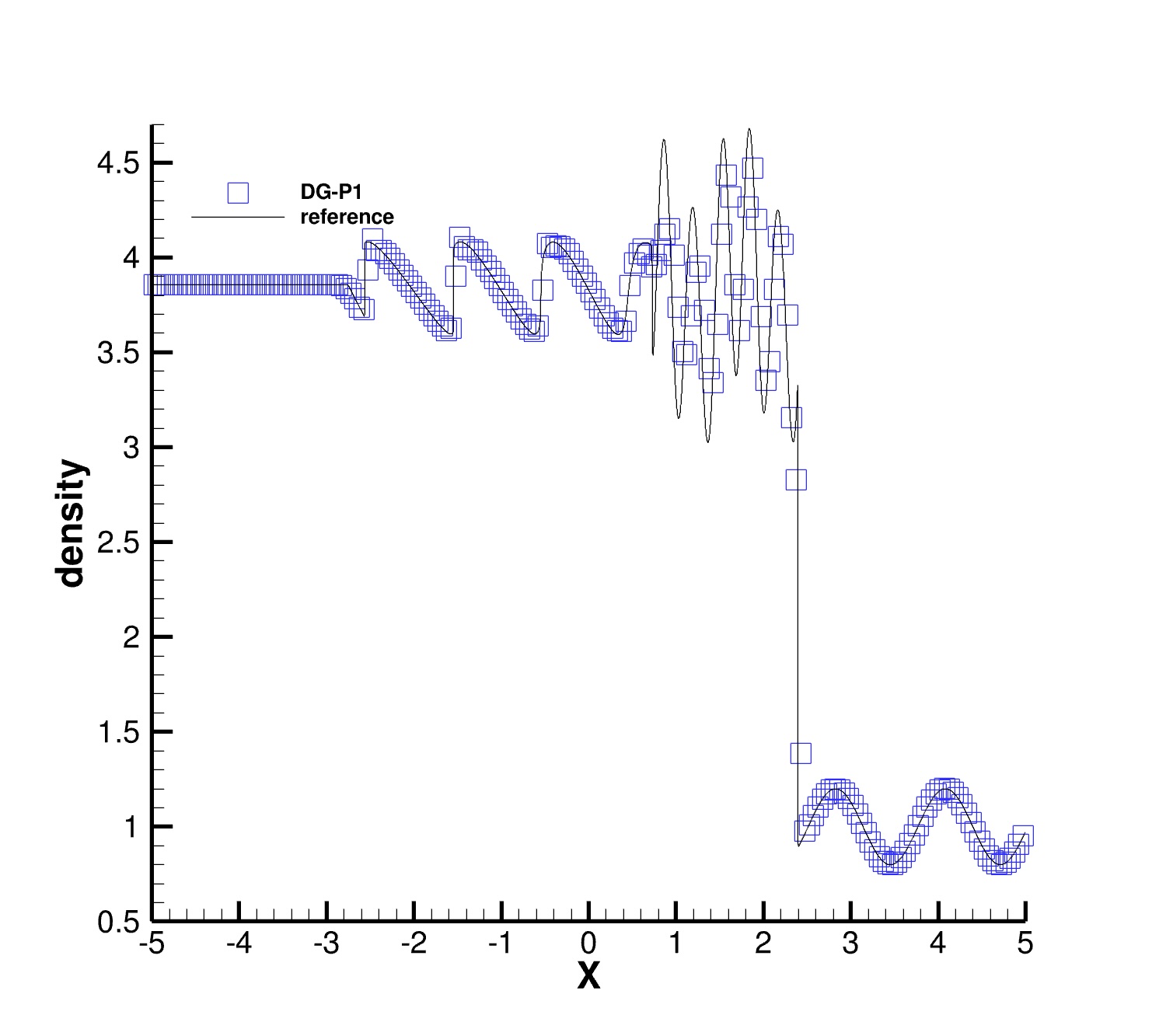}
  \end{minipage}
   \begin{minipage}[t]{0.5\textwidth}  
	\centering
	 \includegraphics[width=\textwidth]{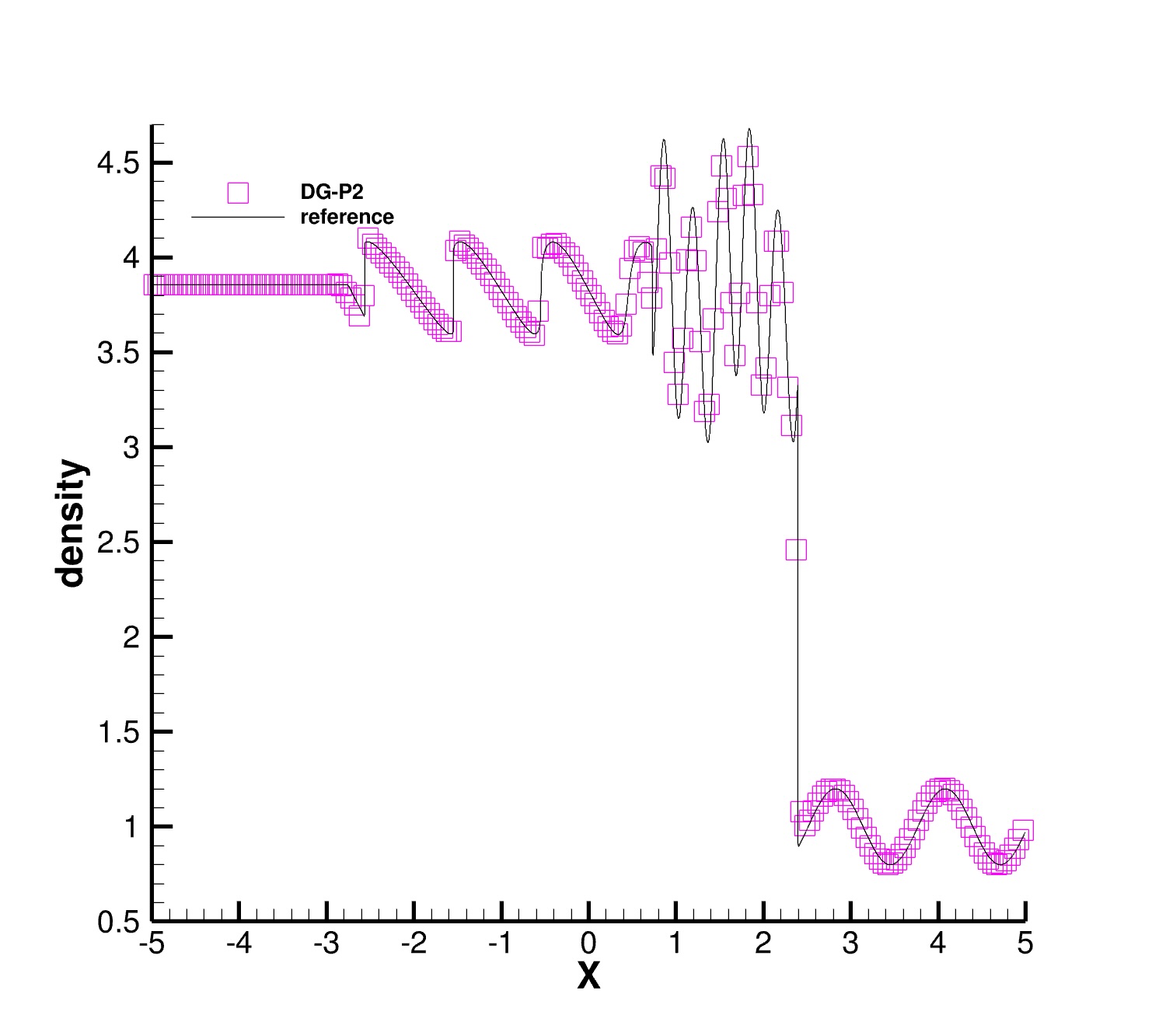 }
  \end{minipage}
  }\\
  \subfloat{
  \begin{minipage}[t]{0.5\textwidth}  
	\centering
	 \includegraphics[width=\textwidth]{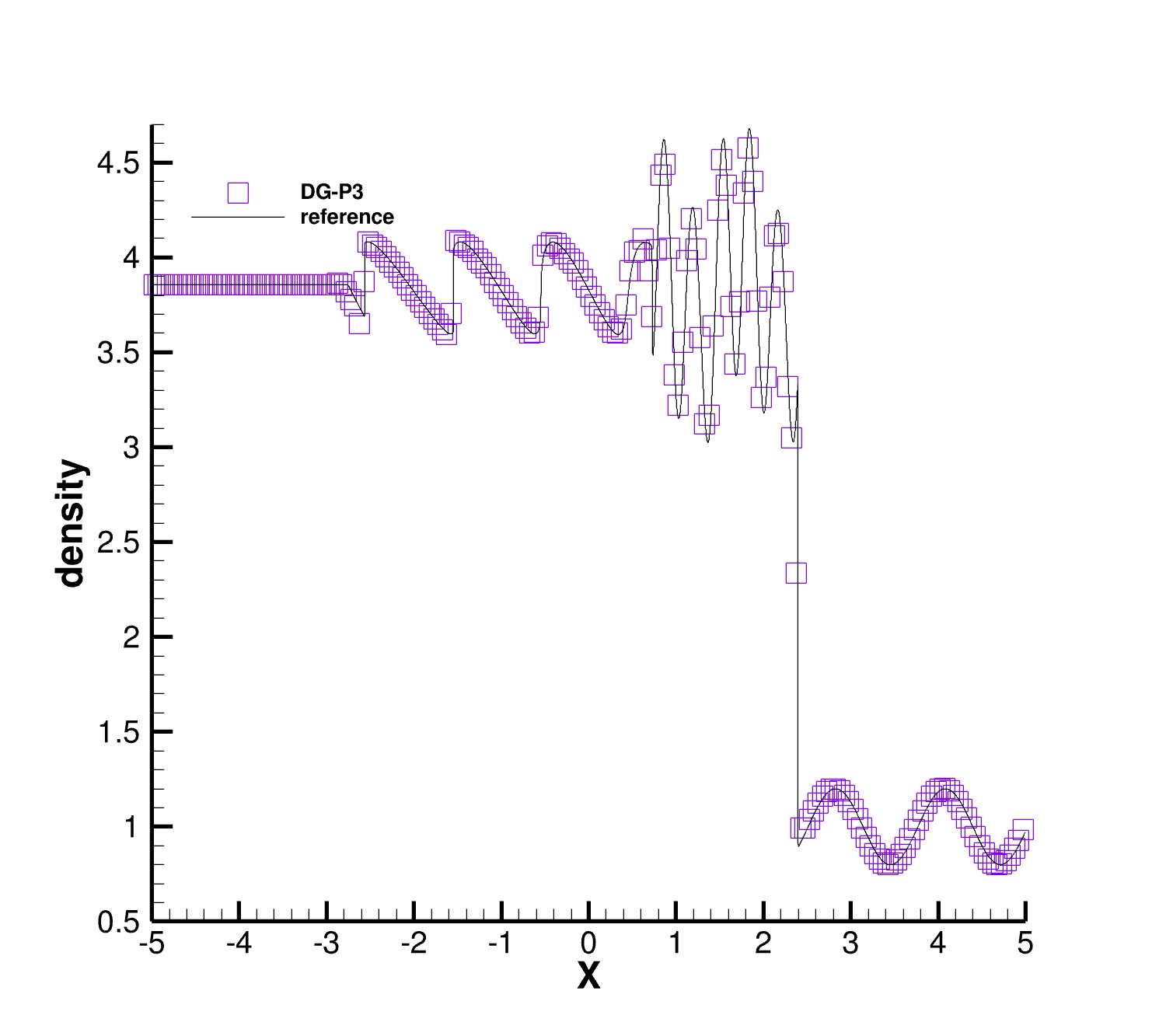 }
  \end{minipage}
  \begin{minipage}[t]{0.5\textwidth}  
	\centering
	 \includegraphics[width=\textwidth]{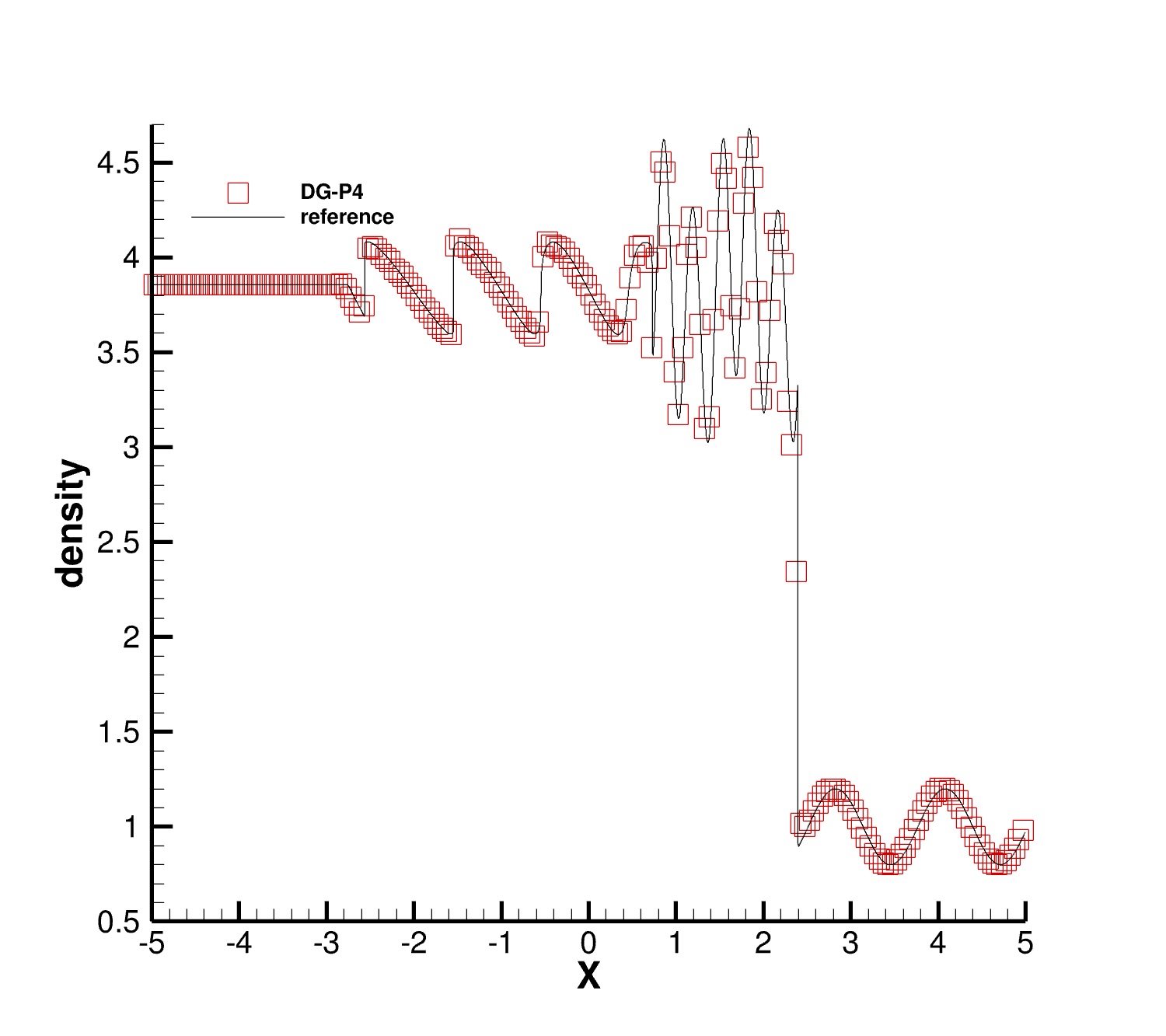 }
  \end{minipage}
  }
  \caption{The density distribution of  Example 4.3 \eqref{so}   with  $N=200.$
   Left top: second-order($p^1$), right top: third-order($p^2$), left bottom: fourth-order($p^3$), right bottom: fifth-order($p^4$).}
  \label{fig:SO}
\end{figure}
\begin{figure}[h!]
	\centering
	\subfloat{
	\begin{minipage}[t]{0.24\textwidth}  
	\centering
	 \includegraphics[width=\textwidth]{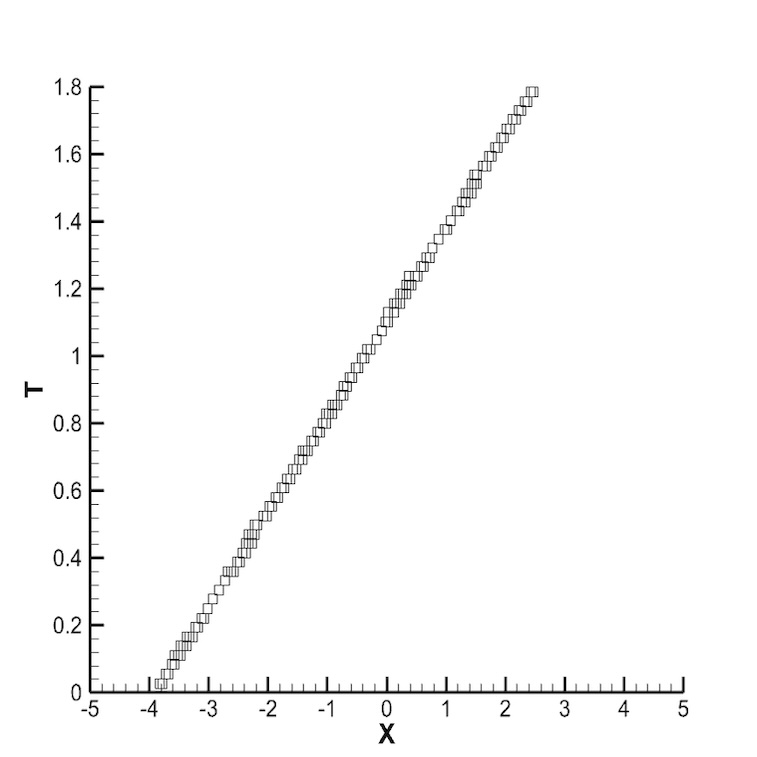 }
  \end{minipage}
   \begin{minipage}[t]{0.24\textwidth}  
	\centering
	 \includegraphics[width=\textwidth]{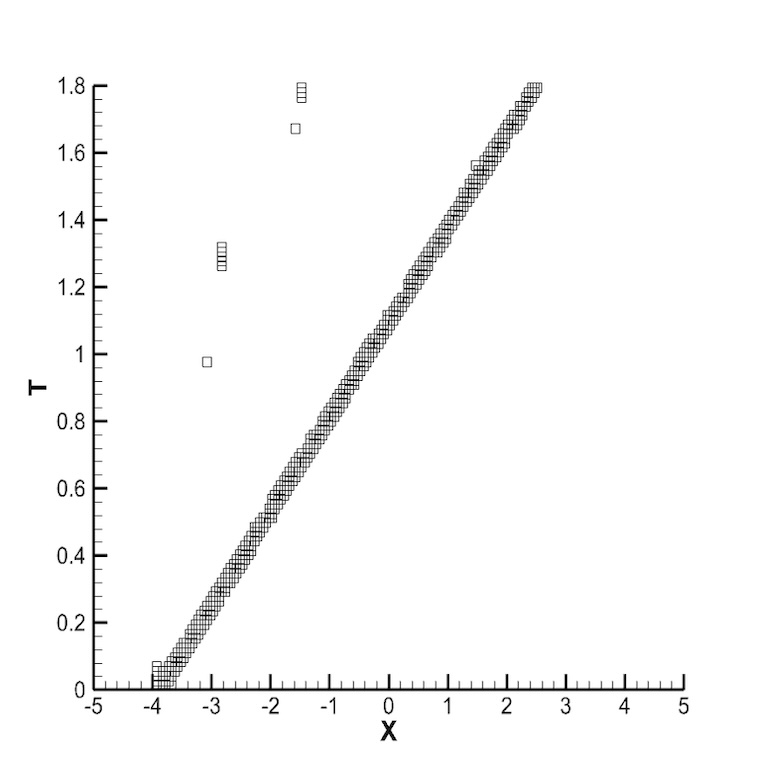 }
  \end{minipage}
  \begin{minipage}[t]{0.24\textwidth}  
	\centering
	 \includegraphics[width=\textwidth]{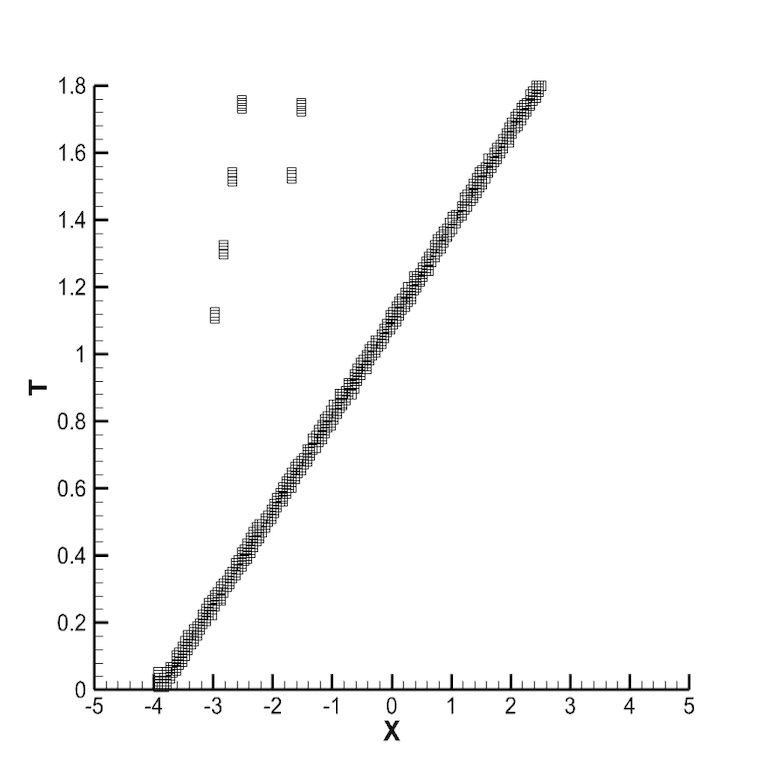 }
  \end{minipage}
  \begin{minipage}[t]{0.24\textwidth}  
	\centering
	 \includegraphics[width=\textwidth]{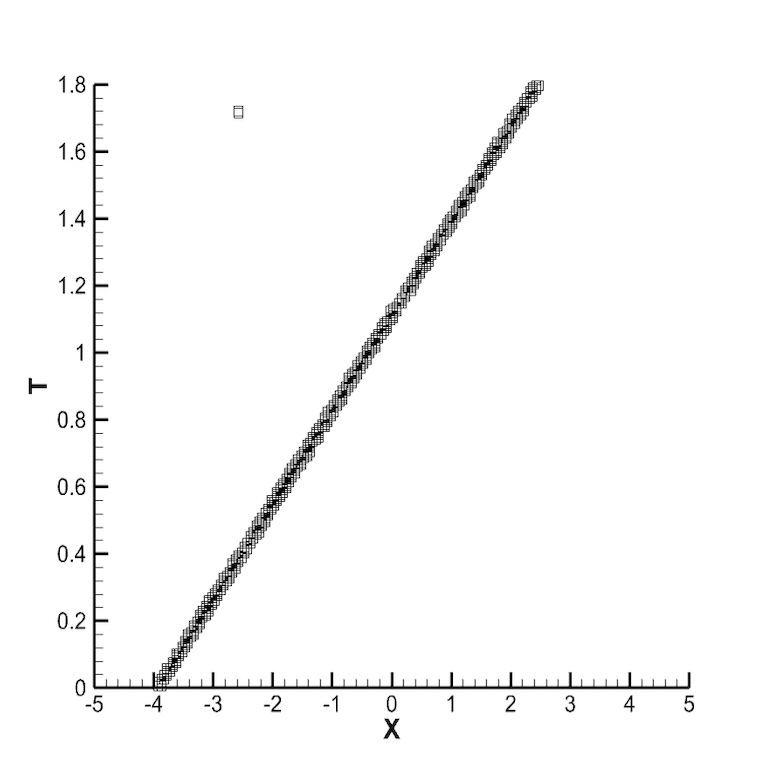 }
  \end{minipage}
}
  \caption{The time history of the troubled cells for Example 4.3 \eqref{so} with  $N=200.$
  DG-HGKS with SHWENO limiter. 
 From left to right: second-order($p^1$), third-order($p^2$), fourth-order($p^3$), fifth-order($p^4$).}
  \label{fig:SOKX}
\end{figure}

{\bf Example 4.4} {Woodward-Colella blast wave}

This case contains the interactions between strong shock waves and contact discontinuities, 
which is a very challenging problem to assess the robustness
of a numerical scheme \cite{DM}. 
The initial condition is given by
\begin{eqnarray}
	\begin{aligned}
\left( {\rho, U, p} \right) = \left\{ \begin{matrix}
{\left( {1, 0, 1000} \right) ,~0 \leq x \leq 0.1,} \\
{\left( {1, 0, 0.01} \right) ,~0.1 < x \leq 0.9,} \\
{\left( {1, 0, 100} \right) ,~0.9 < x \leq 1,}
\end{matrix} \right.	
\end{aligned}
	\label{bw}
\end{eqnarray}

 \begin{figure}[h!]
	\centering
	\subfloat{
	\begin{minipage}[t]{0.5\textwidth}  
	\centering
	 \includegraphics[width=\textwidth]{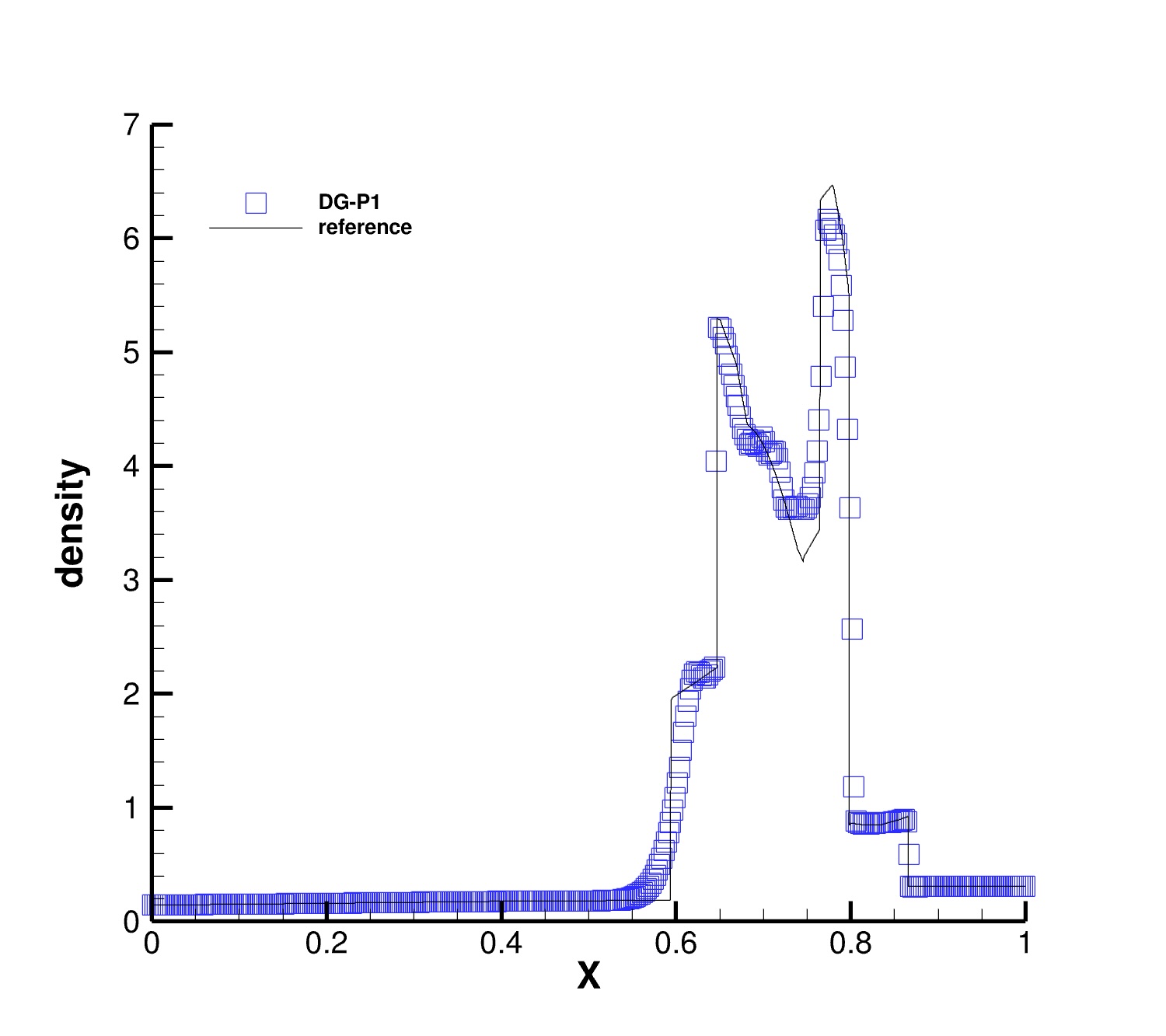 }
  \end{minipage}
   \begin{minipage}[t]{0.5\textwidth}  
	\centering
	 \includegraphics[width=\textwidth]{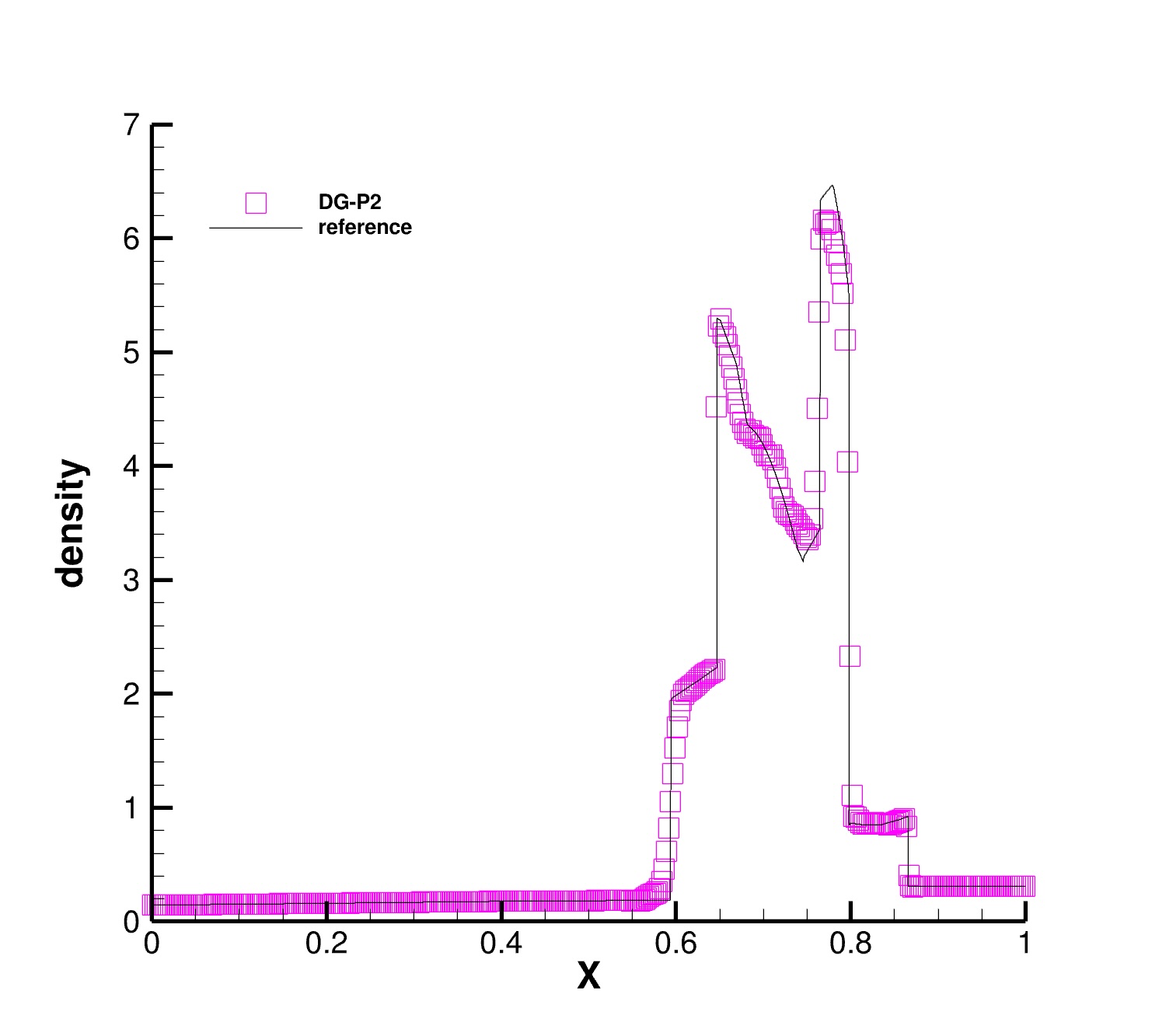 }
  \end{minipage}
  }\\
  \subfloat{
  \begin{minipage}[t]{0.5\textwidth}  
	\centering
	 \includegraphics[width=\textwidth]{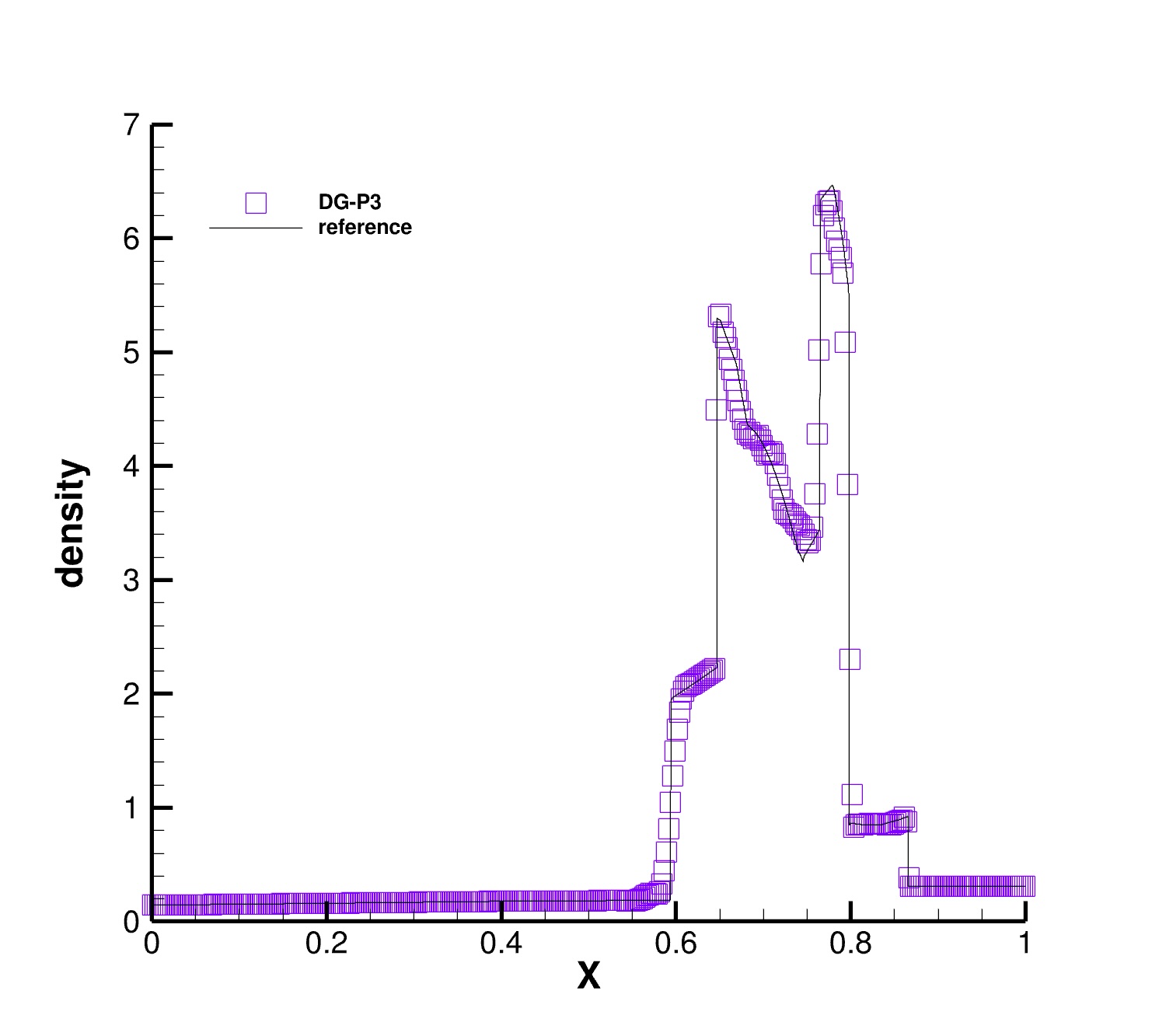 }
  \end{minipage}
  \begin{minipage}[t]{0.5\textwidth}  
	\centering
	 \includegraphics[width=\textwidth]{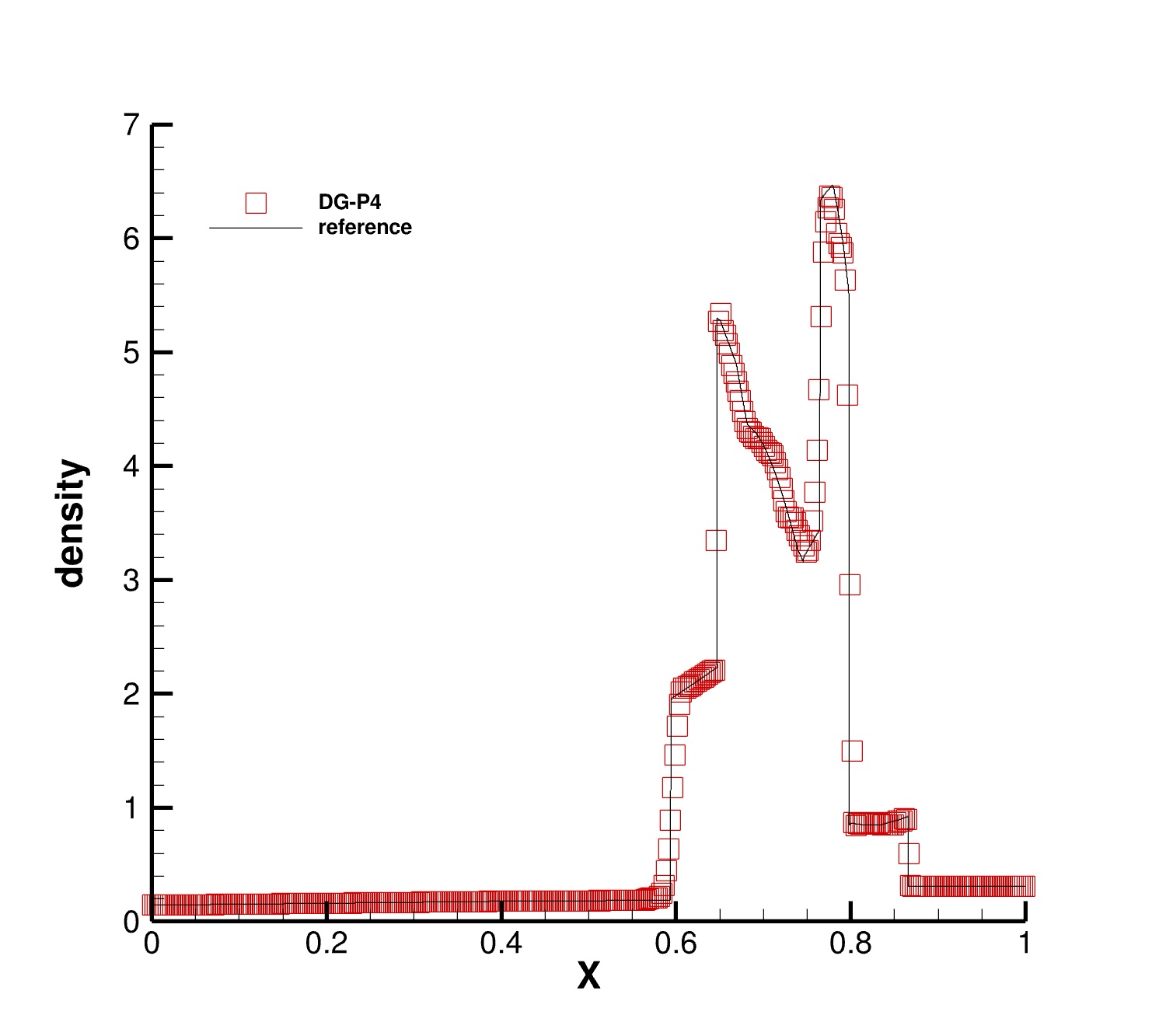 }
  \end{minipage}
  }
  \caption{The density distribution of Example 4.4 \eqref{bw}   with  $N=400.$
  Left top: second-order($p^1$), right top: third-order($p^2$), left bottom: fourth-order($p^3$), right bottom: fifth-order($p^4$).}
  \label{fig:BW}
\end{figure}

\begin{figure}[h!]
	\centering
	\subfloat{
	\begin{minipage}[t]{0.24\textwidth}  
	\centering
	 \includegraphics[width=\textwidth]{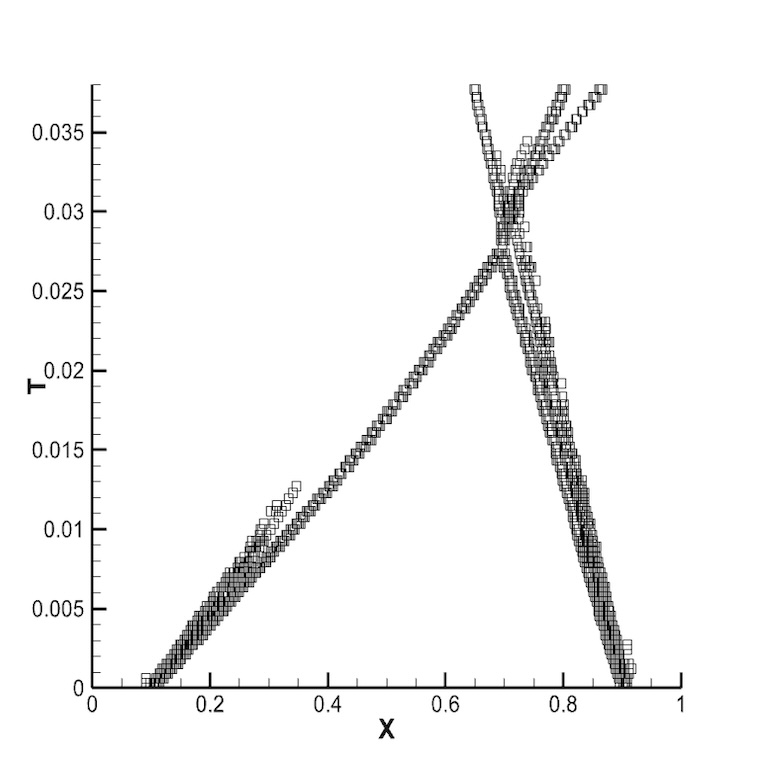 }
  \end{minipage}
   \begin{minipage}[t]{0.24\textwidth}  
	\centering
	 \includegraphics[width=\textwidth]{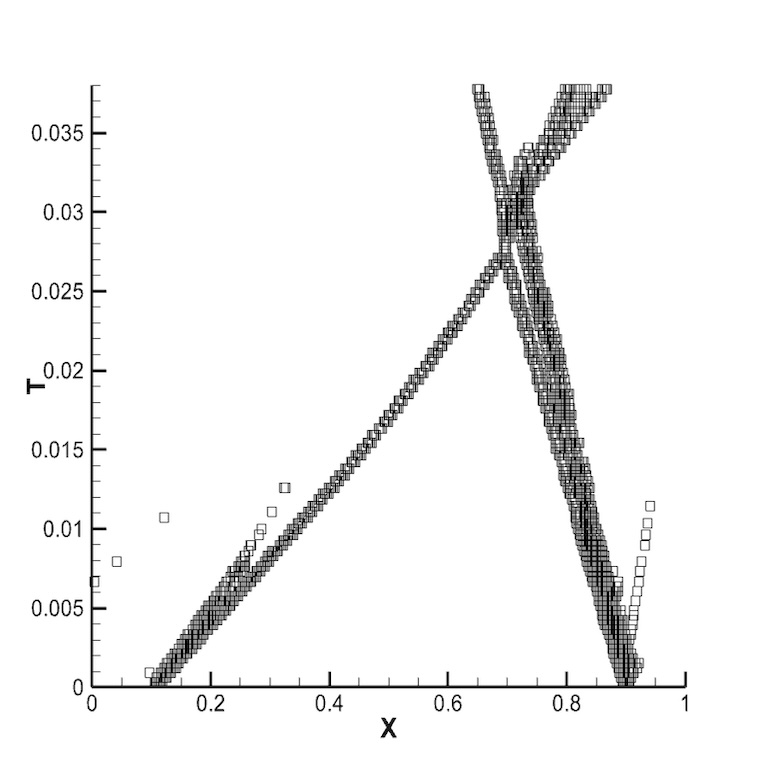 }
  \end{minipage}
  \begin{minipage}[t]{0.24\textwidth}  
	\centering
	 \includegraphics[width=\textwidth]{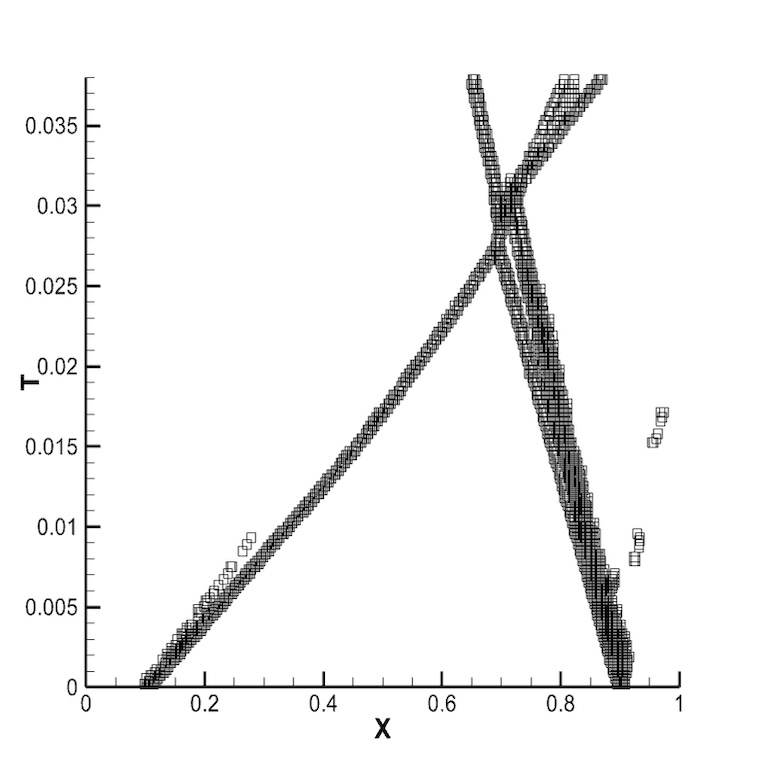 }
  \end{minipage}
  \begin{minipage}[t]{0.24\textwidth}  
	\centering
	 \includegraphics[width=\textwidth]{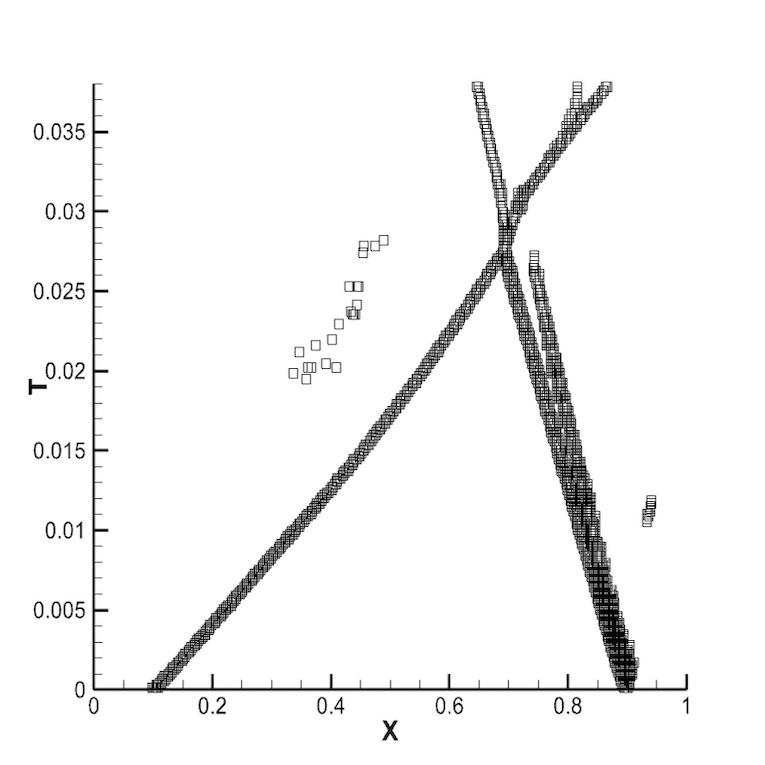 }
  \end{minipage}
}
  \caption{The time history of the troubled cells for Example 4.4 \eqref{bw} with  $N=400.$
  DG-HGKS with SHWENO limiter. 
  From left to right: second-order($p^1$), third-order($p^2$), fourth-order($p^3$), fifth-order($p^4$).}
  \label{fig:BWKX}
\end{figure}
\noindent with the reflective boundary conditions on both sides of the computational domain $[0, 1]$.
The output time is $T=0.038$.  
Fig. \ref{fig:BW} shows the density distributions   under the cells $N=400$,
this reference solution computed by a  fifth-order finite difference WENO scheme \cite{WENO} with 81,920 uniform mesh points.
Fig. \ref{fig:BWKX} is the  corresponding troubled cells  detection domain.
The numerical results show that the DG-HGKS scheme  is effective in dealing with strong shock waves, which is robustness.
 Additionally, it exhibits strong capabilities in solving discontinuous problems, 
and the KXRCF indicator sharply identifies the “troubled cells”.

{\bf Example 4.5} {The modified shock/turbulence interaction }

The computational domain is $\left\lbrack {-5, 5} \right\rbrack$. The initial flow field is \cite{LSY2}
\begin{eqnarray}
	\begin{aligned}
\left( {\rho, U, p} \right) = \left\{ \begin{matrix}
{\left( {1.515695, 0.523346, 1.80500} \right) ,~x \leq - 4.5,} \\
{\left( {\left. 1 + 0.1{\sin\left( 20\pi x \right.} \right), 0, 1} \right) ,~x > - 4.5,}
\end{matrix} \right.	
\end{aligned}
	\label{gp}
\end{eqnarray}
with the inflow/outflow boundary condition.
The output time is $T=5.0$. 
Fig. \ref {fig:GP} show the  density distributions  divided by  $1000$ uniform cells,
compared against a reference solution obtained by the  fifth-order DG-HGKS scheme  with 10,000 uniform mesh points.
The conclusion is that  with the improvement of the accuracy,
the numerical solution of the fifth-order scheme is more effective.

\begin{figure}[h!]
	\centering
	\subfloat{
	\begin{minipage}[t]{0.5\textwidth}  
	\centering
	 \includegraphics[width=\textwidth]{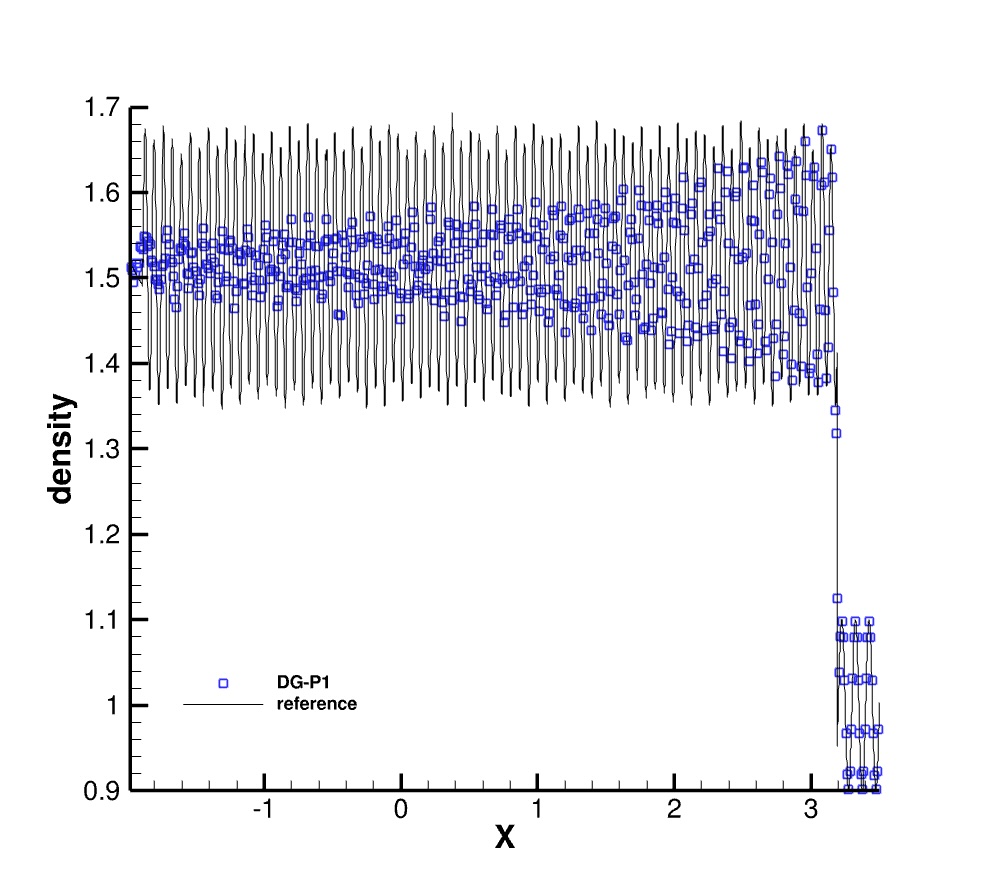 }
  \end{minipage}
   \begin{minipage}[t]{0.5\textwidth}  
	\centering
	 \includegraphics[width=\textwidth]{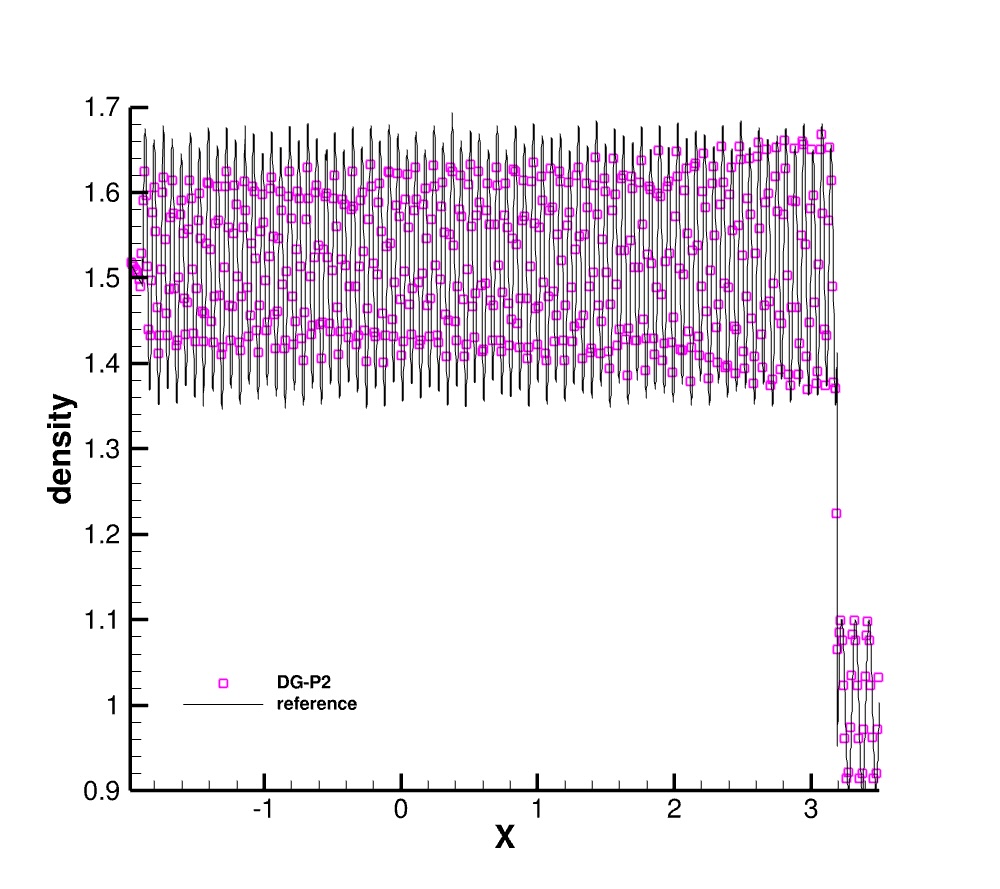 }
  \end{minipage}
  }\\
  \subfloat{
  \begin{minipage}[t]{0.5\textwidth}  
	\centering
	 \includegraphics[width=\textwidth]{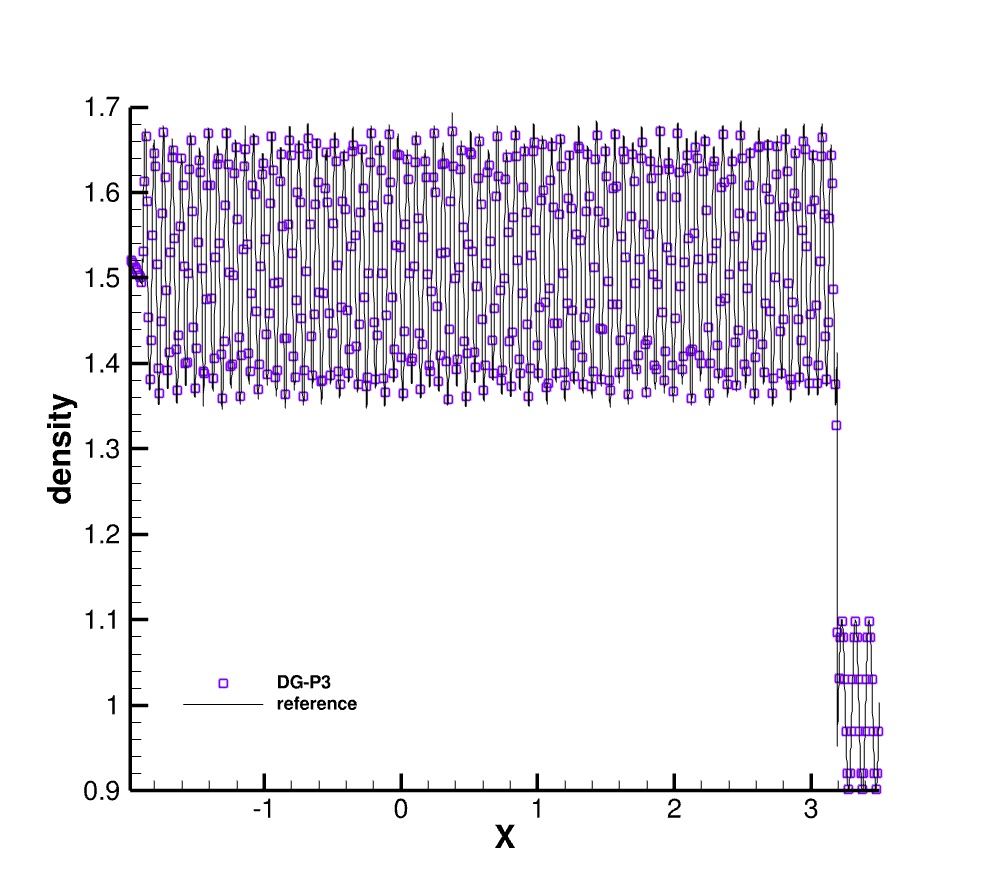 }
  \end{minipage}
  \begin{minipage}[t]{0.5\textwidth}  
	\centering
	 \includegraphics[width=\textwidth]{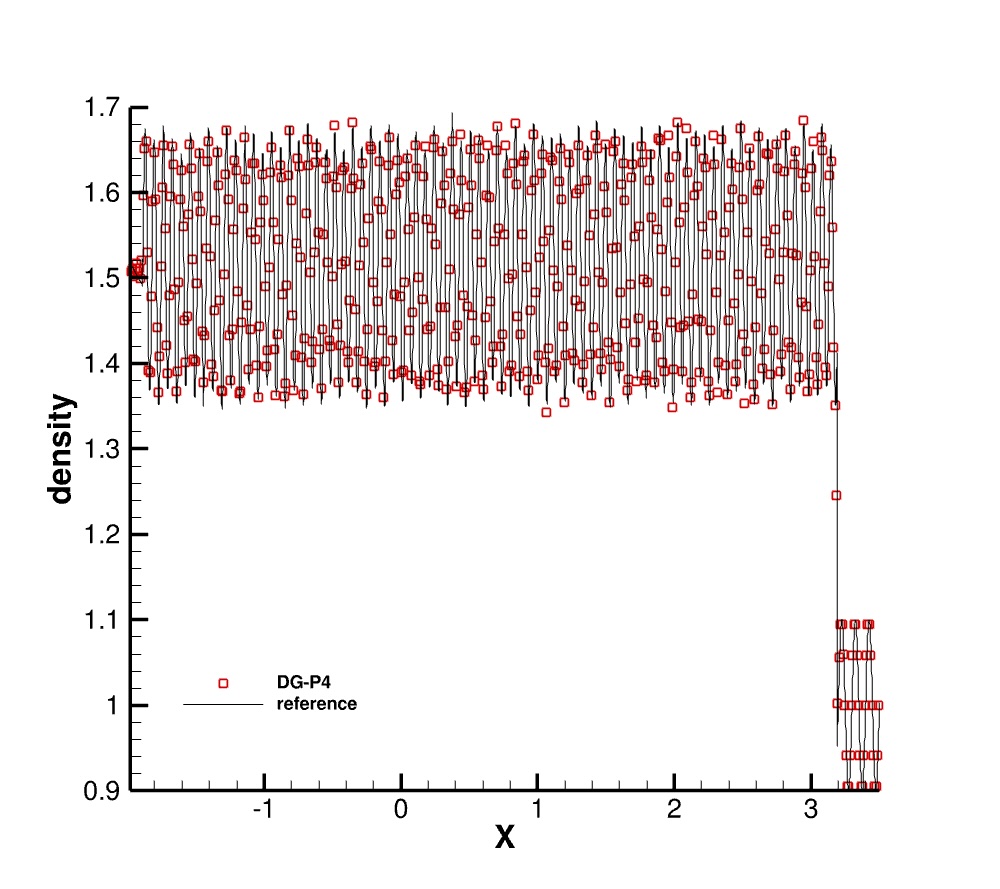 }
  \end{minipage}
  }
  \caption{The density distribution of Example 4.5 \eqref{gp}     with $N=1000.$
  Left top: second-order($p^1$), right top: third-order($p^2$), left bottom: fourth-order($p^3$), right bottom: fifth-order($p^4$).}
  \label{fig:GP}
\end{figure}

 \subsection{Two-dimensional examples}
 
{\bf Example 4.6} {2D linear advection of the density perturbation}

Here we assess the accuracy  of the 2D problem when the solution is linear and smooth.
The initial condition is given by
\begin{eqnarray}
	\left( {\rho, U, V, p} \right) = \left( {1 + 0.2{\sin{\pi (x+y)}}, 0.7, 0.3, 1}\right),
\label{2dorder}
\end{eqnarray}
and under the periodic boundary condition, the analytic solution is
$$
      \left( {\rho, U, V, p} \right) = \left( {1 + 0.2{\sin{\pi (x+y-t)}}, 0.7, 0.3, 1}\right).
$$
The computational domain  $[ {0, 1 } ] \times [ {0, 1 } ]$ is divided into $N \times M$ uniform cells.
The output time is $T=1.0$. 
The numerical errors and  accuracy orders  are shown in the Table \ref{tab:ac2d}.
It is observed that the DG-HGKS scheme  can achieve the uniformly second-order($p^1$), third-order($p^2$), fourth-order($p^3$), and fifth-order($p^4$) in 2D linear problem.

 \begin{table}[h]
\centering
\resizebox{\textwidth}{!}{
\begin{tabular}{|c|c|cccc|cccc|}
\hline
 &   & \multicolumn{4}{c|}{ \textbf{DG without limiter}} & \multicolumn{4}{c|} {\textbf{DG with limiter}} \\ \hline
 &{ { N=M }}   & { {$L^1$ error}} & { {order}} & { {$L^\infty$ error}} & { {order}} & { {$L^1$ error}} & { {order}} & { {$L^\infty$ error }}& { {order}} \\ \hline
\multirow{4}{*}{{ {$\boldsymbol{p^1}$}}} 
& 20   & 4.74e-03 & --      & 1.68e-03 & --     & 3.76e-03 & --      & 1.70e-03 & --\\
& 40   & 1.27e-03 & 1.90 & 4.16e-04 & 2.01 & 9.34e-04 & 2.01 & 4.22e-04 & 2.01 \\
& 80   & 3.38e-04 & 1.90 & 1.03e-04 & 2.01 & 2.32e-04 & 2.01 & 1.05e-04 & 2.00\\
& 160& 8.71e-05 & 1.96 & 2.55e-05 & 2.01 & 5.85e-05 & 1.99 & 2.63e-05 & 2.00\\ \hline
\multirow{4}{*}{{ {$\boldsymbol{p^2}$}}} 
& 20   & 3.42e-04 & --      & 1.93e-05 & --     & 2.17e-04 & --      & 8.78e-06 & --\\
& 40   & 4.53e-05 & 2.92 & 3.65e-06 & 2.41 & 2.82e-05 & 2.94 & 1.50e-06 & 2.55 \\
& 80   & 5.58e-06 & 3.02 & 2.38e-07 & 3.94 & 3.70e-06 & 2.93 & 2.12e-07 & 2.82\\
& 160 & 6.96e-07 & 3.00 & 3.26e-08 & 2.87 & 5.03e-07 & 2.88 & 2.95e-08 & 2.84\\ \hline
\multirow{4}{*}{{ $\boldsymbol{p^3}$}} 
& 20   & 1.17e-05 & --      & 4.76e-06 & --     & 8.80e-06 & --      & 5.19e-06 & --\\
& 40   & 7.69e-07 & 3.93 & 2.84e-07 & 4.07 & 5.02e-07 & 4.13 & 2.92e-07 & 4.15 \\
& 80   & 4.67e-08 & 4.03 & 1.91e-08 & 3.89 & 3.03e-08 & 4.05 & 1.75e-08 & 4.07\\
& 160 & 2.89e-09 & 4.01 & 1.06e-09 & 4.17 & 1.88e-09 & 4.01 & 1.07e-09 & 4.03\\ \hline
\multirow{4}{*}{{ $\boldsymbol{p^4}$}} 
& 20   & 2.98e-07 & --      & 7.79e-08 & --     & 5.13e-06 & --      & 9.16e-07 & --\\
& 40   & 8.40e-09 & 5.15 & 2.18e-09 & 5.16 & 1.49e-07 & 5.11 & 2.39e-08 & 5.26 \\
& 80   & 2.63e-10 & 5.00 & 6.67e-11 & 5.03 & 4.84e-09 & 4.94 & 7.56e-10 & 4.99\\
& 160 & 8.24e-12 & 5.00 & 2.43e-12 & 4.78 & 1.47e-10 & 5.04 & 2.36e-11 & 5.00\\ \hline
\end{tabular}
}
\caption{{\label{tab:ac2d} Accuracy test  for 2D of Example 4.6 \eqref{2dorder}.}}
\end{table}
   
{\bf Example 4.7} {2D isotropic vortex  problem}

The 2D  isotropic vortex  problem  has a nonlinear but smooth solution \cite{Vortex}.
The  initial flow is given by
\begin{eqnarray}
\left\{\begin{matrix}	
{\rho = \left( {1 - \frac{25(\gamma - 1)}{8\gamma\pi^{2}}e^{1 - r^{2}}} \right)^{\frac{1}{\gamma - 1}}} ,\\
{U = 1 - \frac{5}{2\pi}e^{\frac{1 - r^{2}}{2}}\left( {y - 5} \right)} ,\\
{V = 1 + \frac{5}{2\pi}e^{\frac{1 - r^{2}}{2}}\left( {x - 5} \right)} ,\\
{p = \rho^{\gamma} ,	}
\end{matrix}\right.
\label{vor}
\end{eqnarray}
where $r^2=(x-5)^2+(y-5)^2$. The periodic boundary condition is employed in both directions. The exact solution is the vortex along the upper right direction with velocities $(U, V)=(1, 1)$. We compute the numerical solution at the output time $T= 1.0$. 
The accuracy results are listed in Table \ref{tab:ac2d2}, 
which  can  observe that the new scheme can still achieve the uniformly fifth-order in the 2D isotropic vortex  problem.

 \begin{table}[h]
\centering
\resizebox{\textwidth}{!}{
\begin{tabular}{|c|c|cccc|cccc|}
\hline
 &   & \multicolumn{4}{c|}{ \textbf{DG without limiter}} & \multicolumn{4}{c|} {\textbf{DG with limiter}} \\ \hline
 &{ { N=M }}   & { {$L^1$ error}} & { {order}} & { {$L^\infty$ error}} & { {order}} & { {$L^1$ error}} & { {order}} & { {$L^\infty$ error }}& { {order}} \\ \hline
\multirow{4}{*}{{ {$\boldsymbol{p^1}$}}} 
& 20   & 1.46e-01 & --      & 1.49e-02 & --     & 2.11e-01 & --      & 2.81e-02 & --\\
& 40   & 3.69e-02 & 1.99 & 3.91e-03 & 1.81 & 6.34e-02 & 1.73 & 9.40e-03 & 1.58 \\
& 80   & 9.35e-03 & 1.98 & 1.00e-03 & 1.92 & 1.69e-02 & 1.90 & 2.71e-03 & 1.80\\
& 160 & 2.46e-03 & 1.93 & 2.65e-04 & 1.92 & 4.33e-03 & 1.97 & 7.29e-04 & 1.89\\ \hline
\multirow{4}{*}{{ {$\boldsymbol{p^2}$}}} 
& 20   & 6.77e-02 & --      & 1.59e-03& --     & 5.45e-02 & --       & 7.35e-03 & --\\
& 40   & 1.65e-03 & 3.25 & 1.31e-04 & 3.61 & 2.86e-03 & 4.25 & 5.82e-04 & 3.66 \\
& 80   & 1.83e-04 & 3.17 & 1.39e-05 & 3.23 & 3.91e-04 & 2.87 & 1.01e-04 & 2.52\\
& 160 & 2.18e-05 & 3.07 & 1.93e-06 & 2.85 & 5.00e-05 & 2.97 & 1.32e-05 & 2.94\\ \hline
\multirow{4}{*}{{ $\boldsymbol{p^3}$}} 
& 20   & 1.80e-03 & --      & 3.44e-04 & --     & 5.51e-02 & --      & 2.24e-02 & --\\
& 40   & 1.24e-04 & 3.86 & 3.08e-05 & 3.48 & 7.23e-04 & 6.25 & 2.33e-04 & 6.59 \\
& 80   & 8.27e-06 & 3.90 & 2.04e-06 & 3.92 & 3.85e-05 & 4.23 & 1.27e-05 & 4.20\\
& 160 & 6.22e-07 & 3.73 & 1.21e-07 & 4.07 & 2.28e-06 & 4.08 & 6.91e-07 & 4.19\\ \hline
\multirow{4}{*}{{ $\boldsymbol{p^4}$}} 
& 20   & 6.77e-05 & --      & 2.20e-05 & --     & 7.73e-03 & --      & 3.23e-04 & --\\
& 40   & 2.18e-06 & 4.95 & 7.15e-07 & 4.94 & 1.87e-04 & 5.37 & 4.95e-06 & 6.03 \\
& 80   & 6.86e-08 & 4.99 & 2.25e-08 & 4.99 & 6.36e-06 & 4.88 & 1.39e-07 & 5.16\\
& 160 & 2.14e-09 & 5.00 & 7.06e-10 & 5.00 & 2.15e-07 & 4.89 & 3.98e-09 & 5.12\\ \hline
\end{tabular}
}
\caption{{\label{tab:ac2d2} Accuracy test  for 2D of Example 4.7 \eqref{vor}.}}
\end{table}

{\bf Example 4.8} {2D Riemann problems}

To evaluate the capability of DG-HGKS in solving two-dimensional problems, 
we tested two cases of 2D Riemann problems.
The computational domain is $\left\lbrack {0, 1} \right\rbrack \times \lbrack {0, 1}\rbrack$,
and the non-reflecting boundary condition is used at all boundaries. 

The first problem is the  interaction between shocks and contact continuities \cite{R1}.
The initial condition is given by
\begin{eqnarray}
\begin{aligned}
\left( {\rho, U, V, p} \right) = \left\{ \begin{matrix}
{\left( {1.5, 0.0, 0.0, 1.5} \right) ,             x \geq 0.7, y \geq 0.7 ,} \\
{\left( {0.5323, 1.206, 0.0, 0.3} \right),      x \le 0.7, y\geq 0.7 ,} \\
{\left( {0.138, 1.206, 1.206, 0.029} \right),  x\le 0.7, y \le0.7 ,} \\
{\left( {0.5065, 0.0, 0.8939, 0.35} \right),   x \geq 0.7, y \le 0.7 .}
 \end{matrix}\right.
\end{aligned}	
\label{2dr3}
\end{eqnarray}
The output time is $T=0.6$.
The domain is divided into  300$\times$300 cells uniformly.
The density distributions   are shown in  Fig. \ref{fig:R3KX}, and Fig. \ref{fig:R3KXD} is the troubled cells detection domain at the end of time $T=0.6$. 
\begin{figure}[h!]
	\centering
	\subfloat{
	\begin{minipage}[t]{0.5\textwidth}  
	\centering
	 \includegraphics[width=\textwidth]{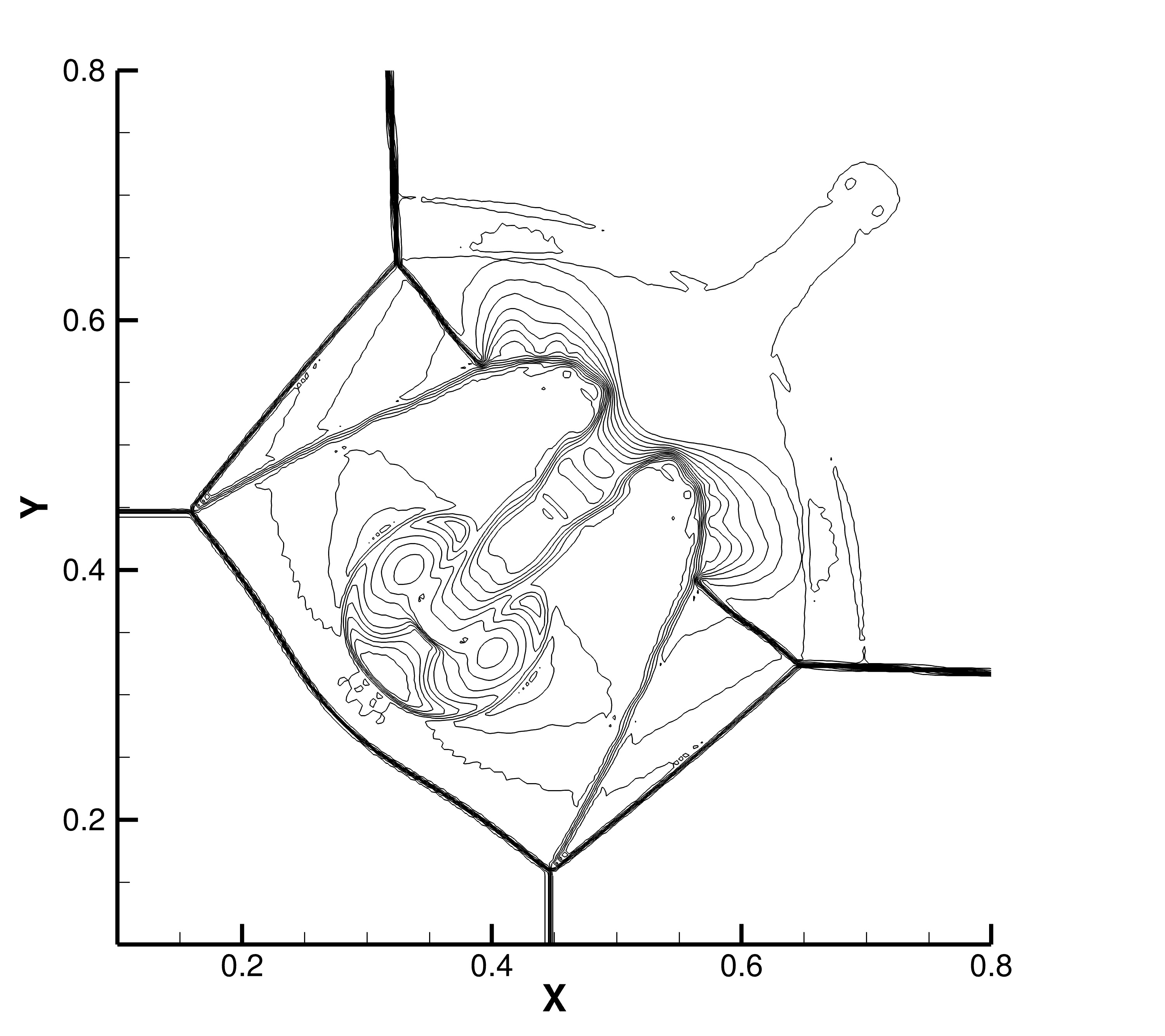}
  \end{minipage}
   \begin{minipage}[t]{0.5\textwidth}  
	\centering
	 \includegraphics[width=\textwidth]{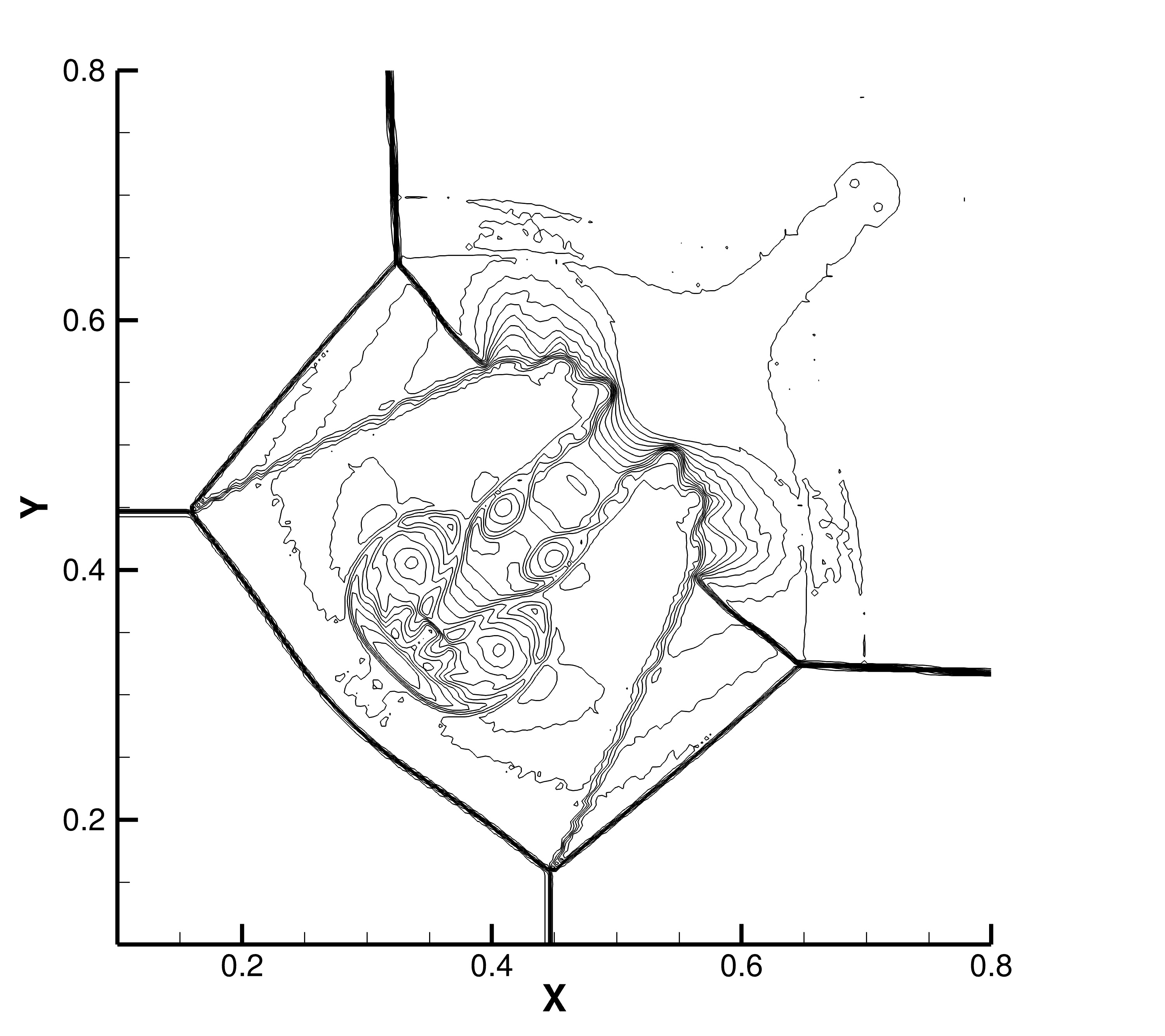}
  \end{minipage}
  }\\
  \subfloat{
  \begin{minipage}[t]{0.5\textwidth}  
	\centering
	 \includegraphics[width=\textwidth]{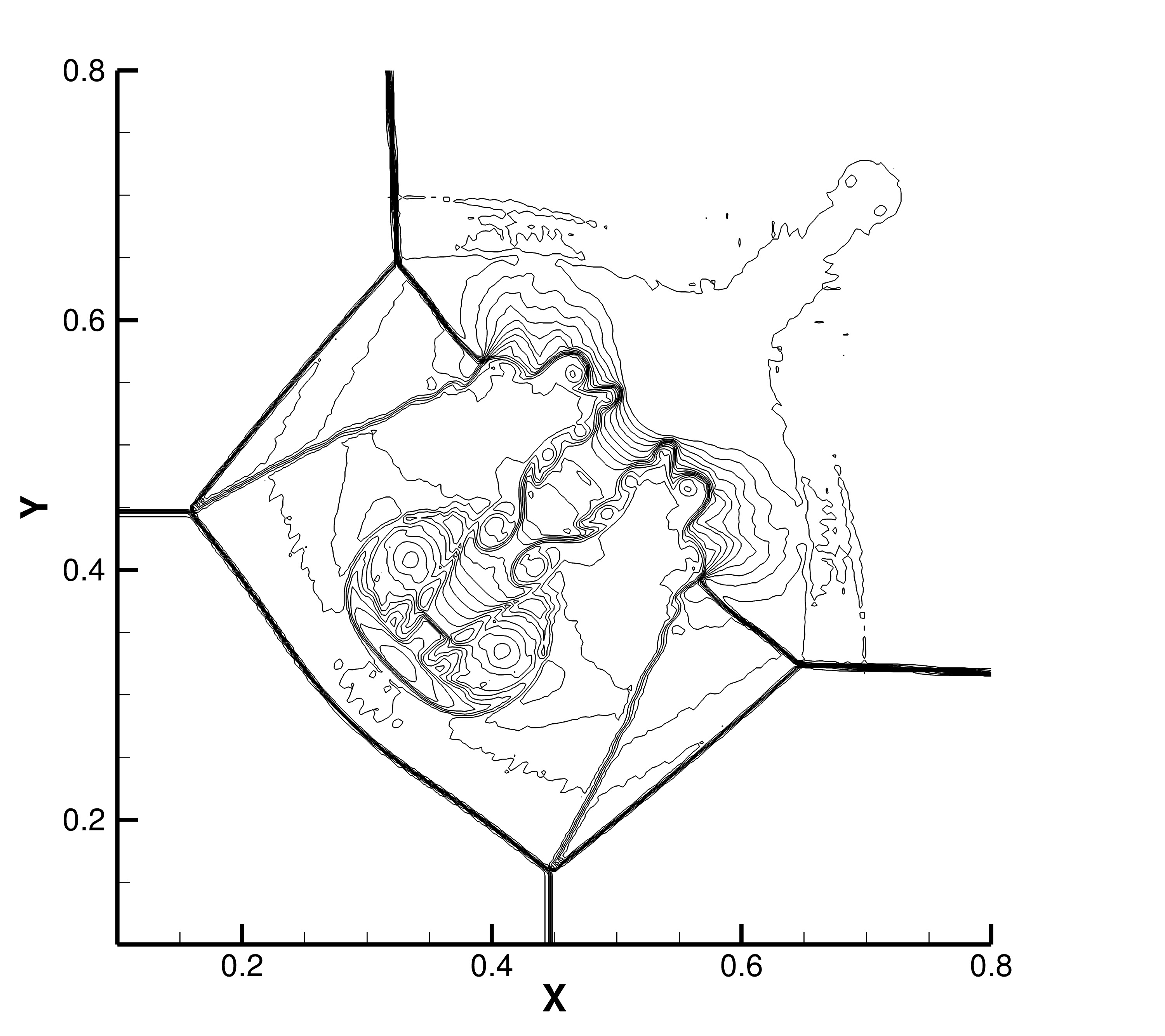}
  \end{minipage}
  \begin{minipage}[t]{0.5\textwidth}  
	\centering
	 \includegraphics[width=\textwidth]{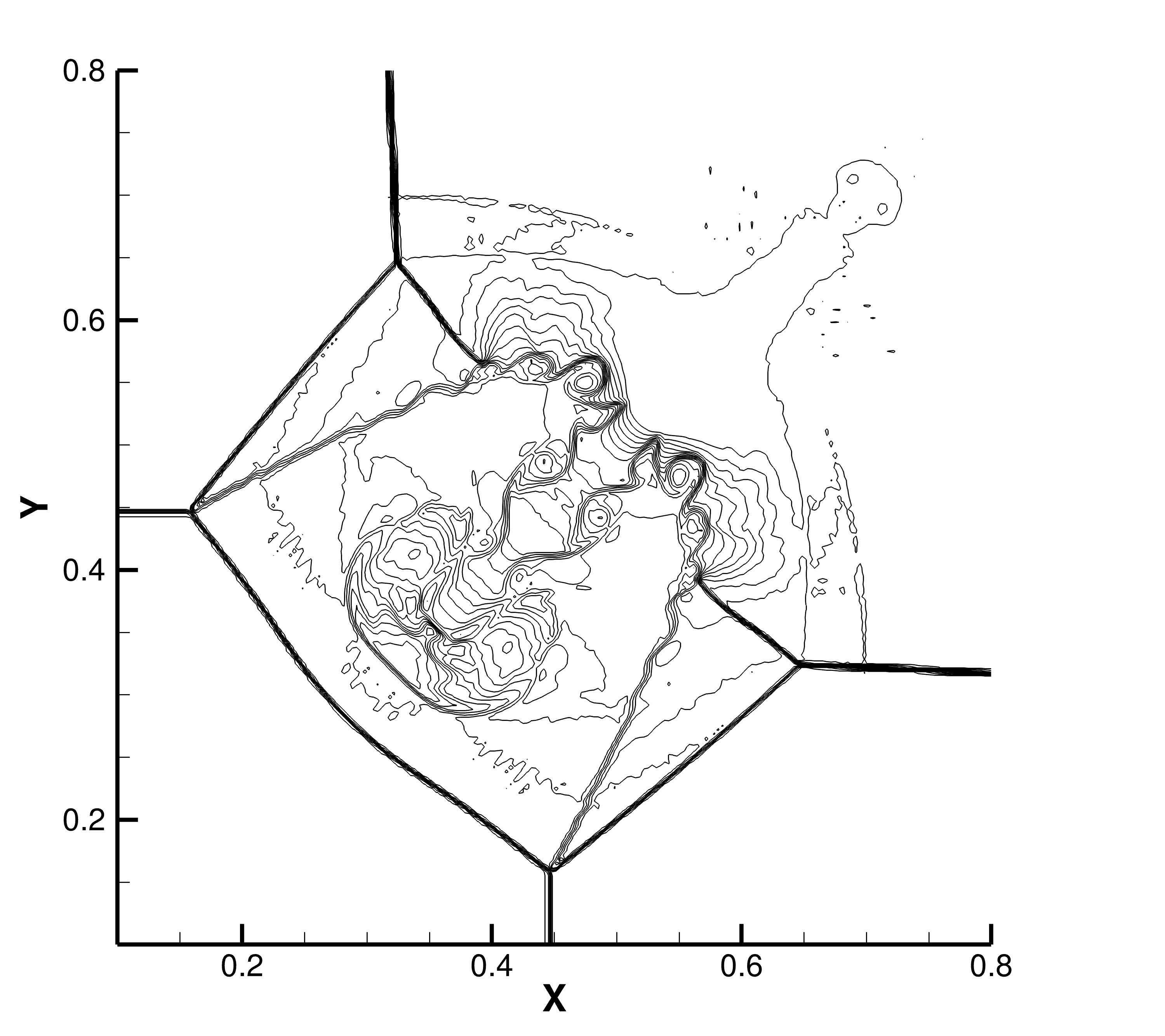}
  \end{minipage}
  }
  \caption{The  density distribution of Example 4.8 \eqref{2dr3}  with 300$\times$300.
  Left top: second-order($p^1$), right top: third-order($p^2$), left bottom: fourth-order($p^3$), right bottom: fifth-order($p^4$).
  20 contours are drawn from 0.15 to 1.5.}
  \label{fig:R3KX}
\end{figure}  
\begin{figure}[h!]
	\centering
	\subfloat{
	\begin{minipage}[t]{0.24\textwidth}  
	\centering
	 \includegraphics[width=\textwidth]{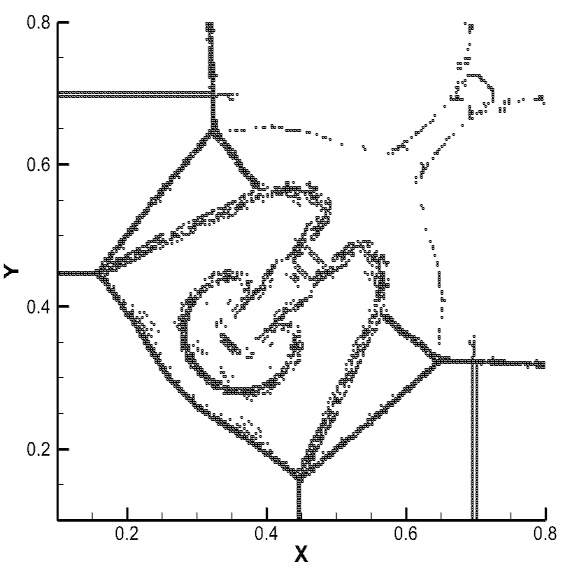}
  \end{minipage}
   \begin{minipage}[t]{0.24\textwidth}  
	\centering
	 \includegraphics[width=\textwidth]{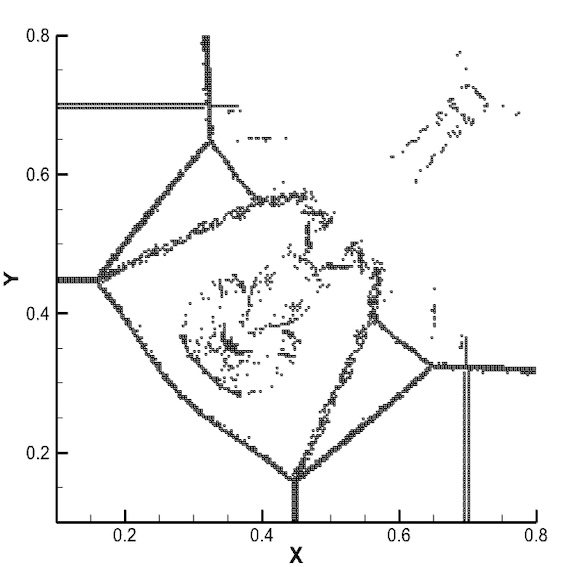}
  \end{minipage}
  \begin{minipage}[t]{0.24\textwidth}  
	\centering
	 \includegraphics[width=\textwidth]{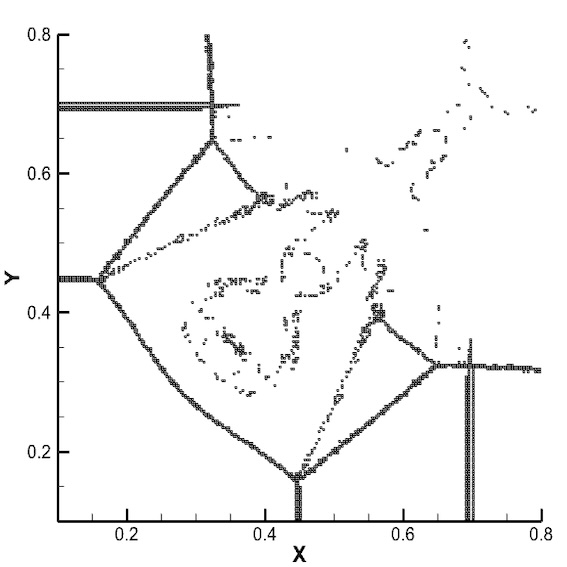}
  \end{minipage}
  \begin{minipage}[t]{0.24\textwidth}  
	\centering
	 \includegraphics[width=\textwidth]{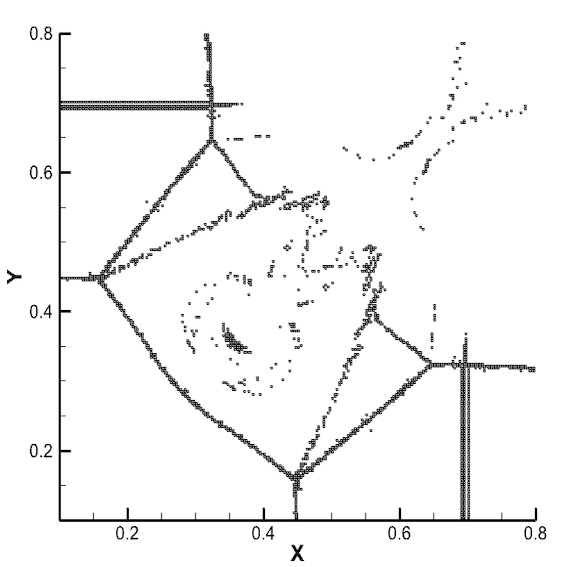}
  \end{minipage}
  }
  \caption{The troubled cells of Example 4.8 \eqref{2dr3}   at the end of time with 300$\times$300.
  From left to right: second-order($p^1$), third-order($p^2$), fourth-order($p^3$), fifth-order($p^4$).}
  \label{fig:R3KXD}
\end{figure}

%

The second problem is the shear instabilities among four initial contact discontinuities \cite{R3}.
The initial condition is given by 
\begin{eqnarray}
	\begin{aligned}
\left( {\rho, U, V, p} \right) = \left\{ \begin{matrix}
{\left( {1.0, 0.75, - 0.5, 1.0} \right) , x > 0.5, y > 0.5 ,} \\
{\left( {2.0, 0.75, 0.5, 1.0} \right) , x < 0.5, y > 0.5 ,} \\
{\left( {1.0, -0.75, 0.5, 1.0} \right) , x < 0.5, y < 0.5 ,} \\
{\left( {3.0, -0.75, -0.5, 1.0} \right) , x > 0.5, y < 0.5 .}
\end{matrix} \right.
		\end{aligned}
	\label{2dr2}
\end{eqnarray}  
The output time is $T=0.8$.
The density distributions    with 200$\times$200 cells are shown in  Fig. \ref{fig:R2KX}, 
and Fig. \ref{fig:R2KXD} is the  detection domain at the end of time $T=0.8$.  
It is evident that the DG-HGKS method effectively  captures the vortex rolling in 2D Riemann problems, 
with the flow structure being resolved more accurately as the order of the scheme increases.
\begin{figure}[h!]
	\centering
	\subfloat{
	\begin{minipage}[t]{0.5\textwidth}  
	\centering
	 \includegraphics[width=\textwidth]{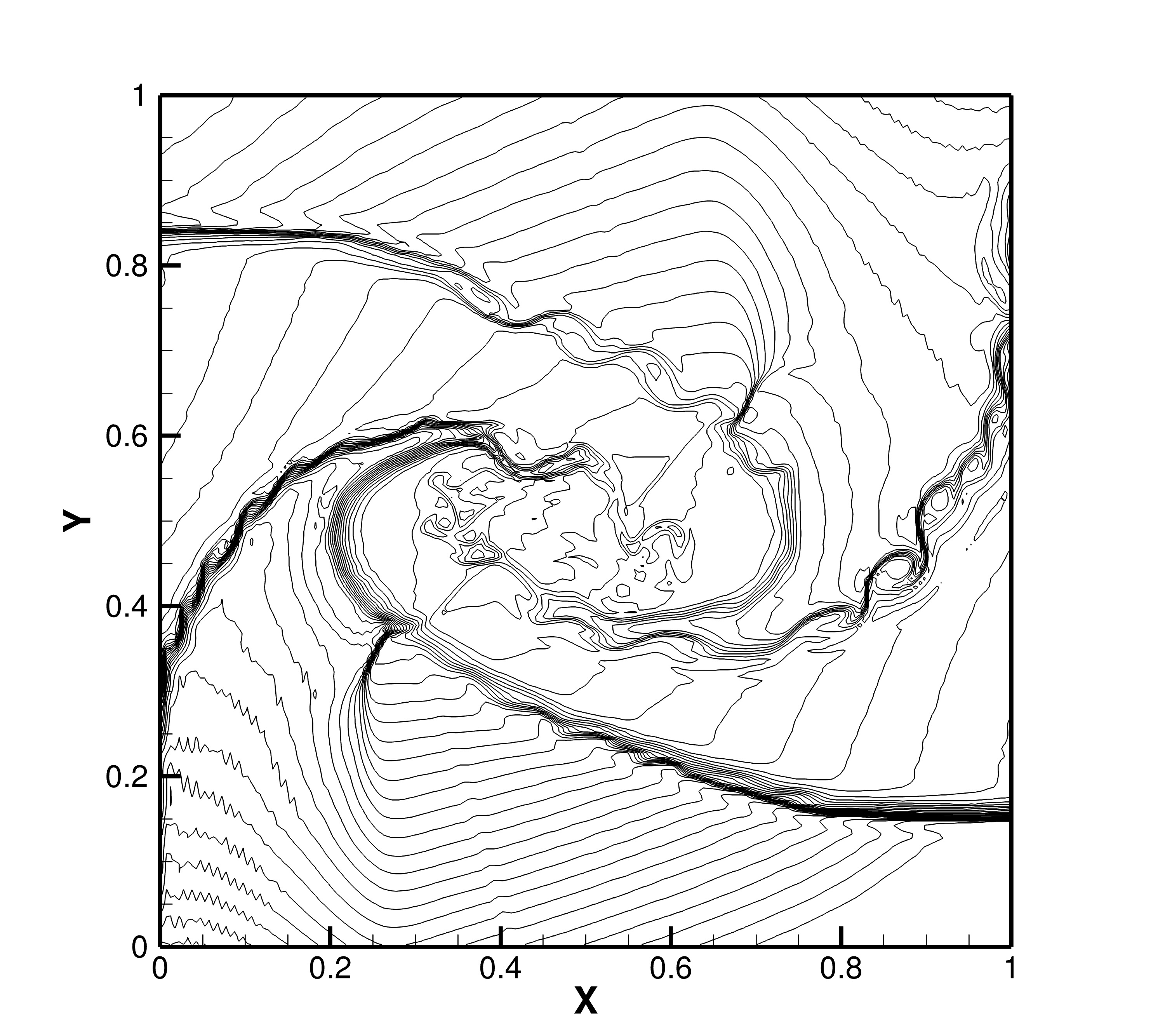}
  \end{minipage}
   \begin{minipage}[t]{0.5\textwidth}  
	\centering
	 \includegraphics[width=\textwidth]{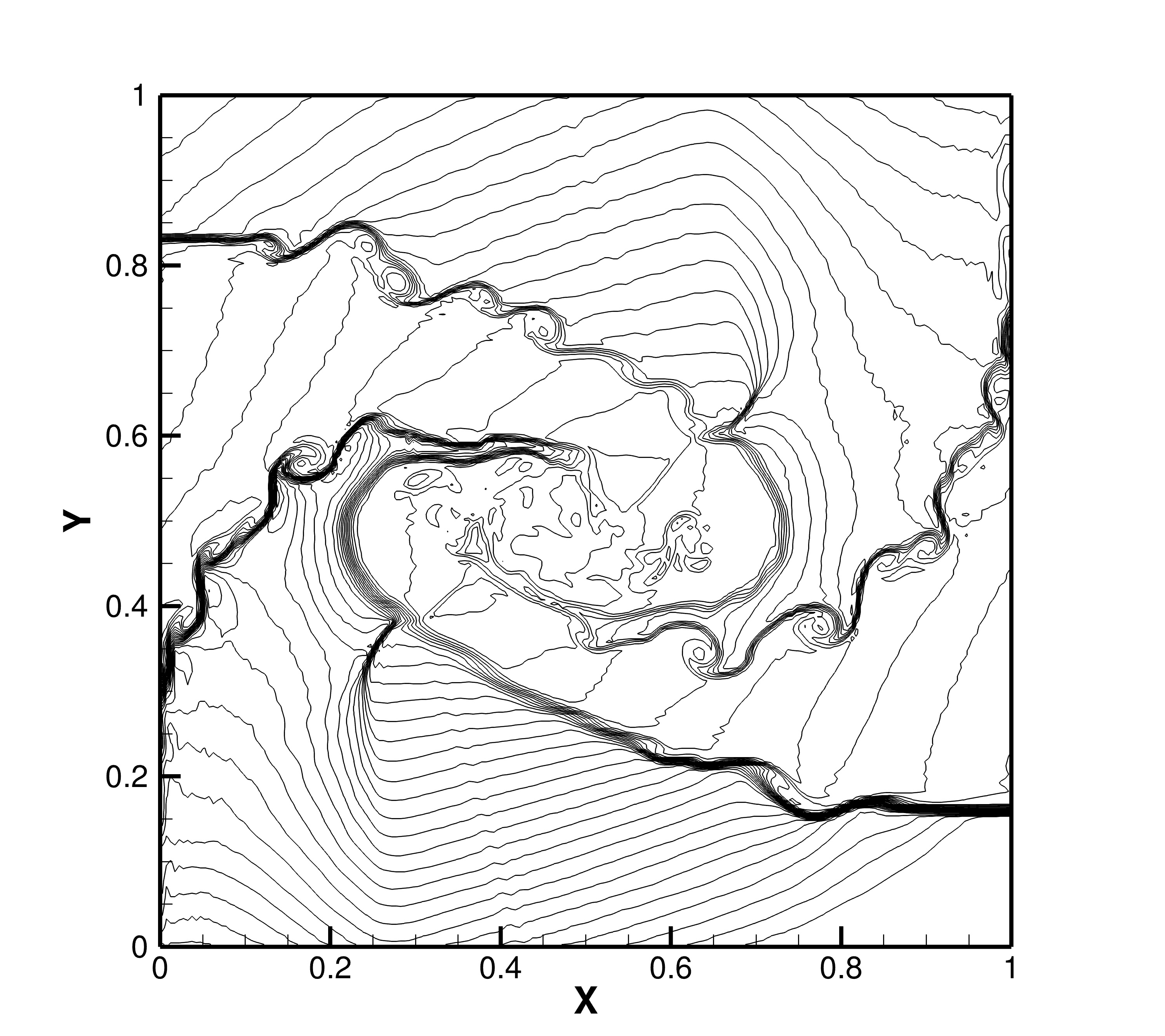}
  \end{minipage}
 }\\
  \subfloat{
  \begin{minipage}[t]{0.5\textwidth}  
	\centering
	 \includegraphics[width=\textwidth]{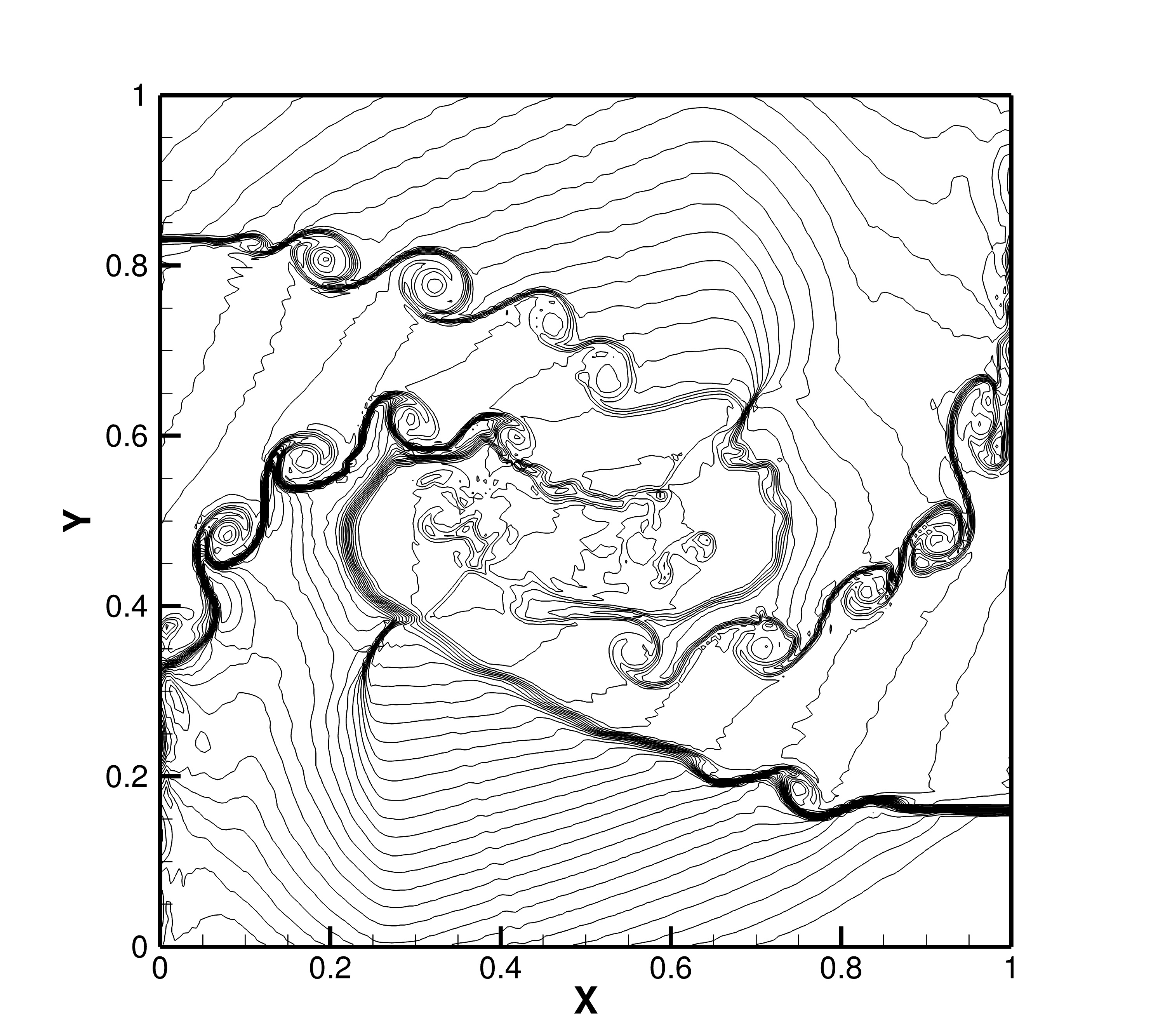}
  \end{minipage}
  \begin{minipage}[t]{0.5\textwidth}  
	\centering
	 \includegraphics[width=\textwidth]{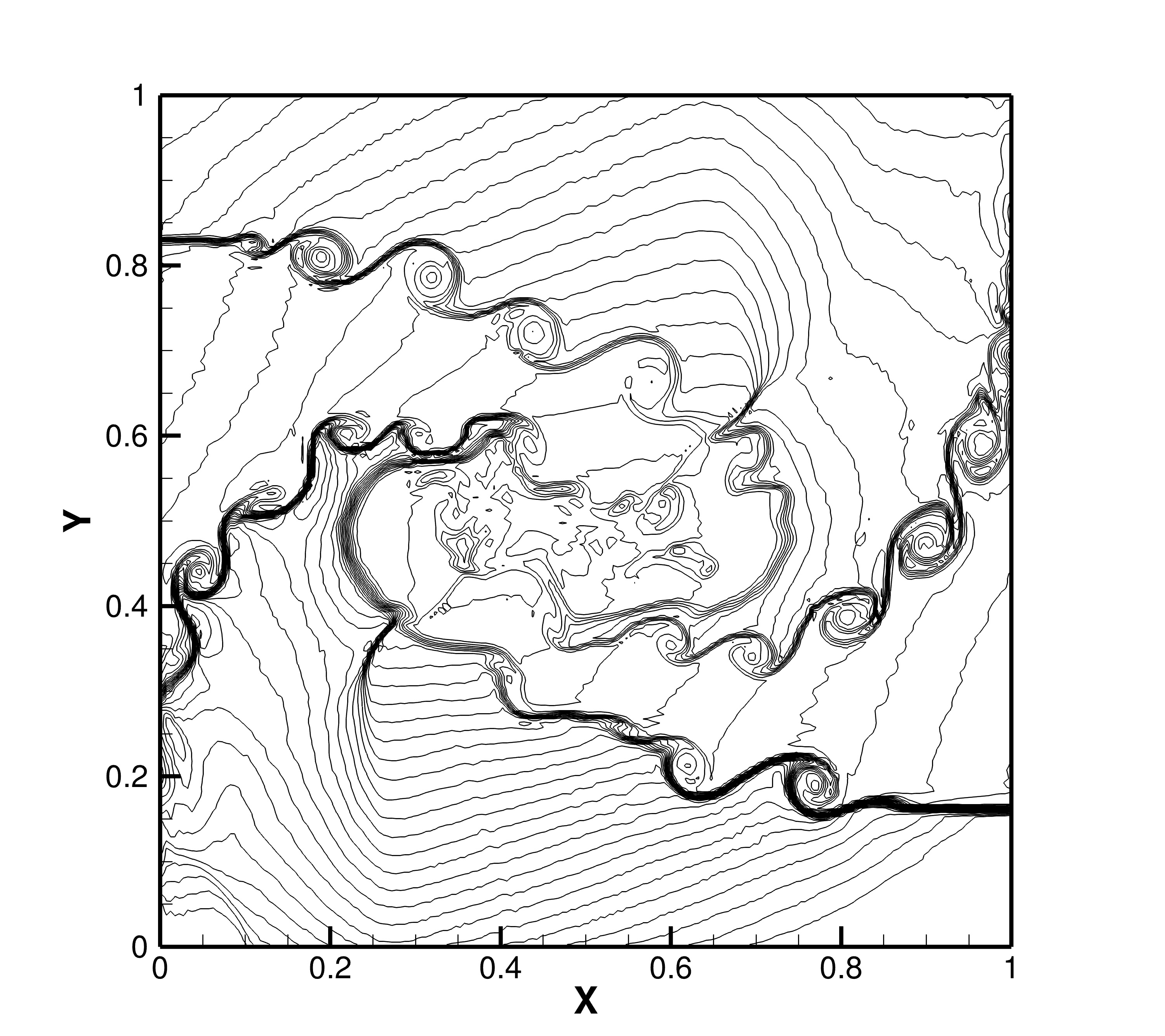}
  \end{minipage}
  }
  \caption{The density distribution  of Example 4.8 \eqref{2dr2}   with 200$\times$200.
 Left top: second-order($p^1$), right top: third-order($p^2$), left bottom: fourth-order($p^3$), right bottom: fifth-order($p^4$).
  30 contours are drawn from 0.2 to 2.2.}
  \label{fig:R2KX}
\end{figure}  
\begin{figure}[h!]
	\centering
	\subfloat{
	\begin{minipage}[t]{0.24\textwidth}  
	\centering
	 \includegraphics[width=\textwidth]{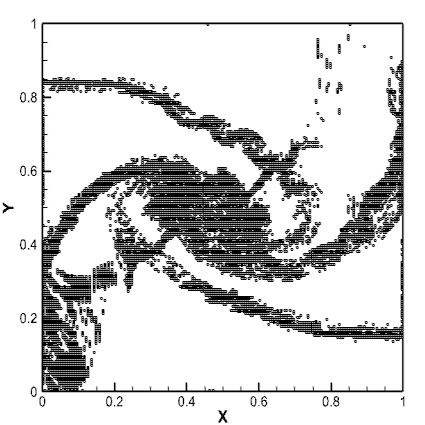}
  \end{minipage}
   \begin{minipage}[t]{0.24\textwidth}  
	\centering
	 \includegraphics[width=\textwidth]{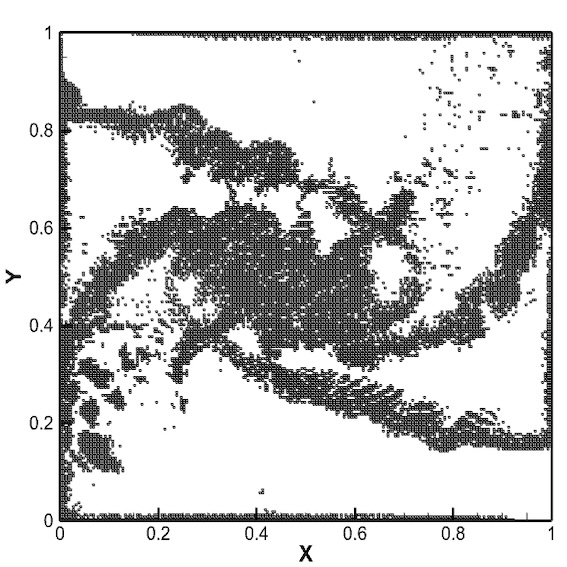}
  \end{minipage}
  \begin{minipage}[t]{0.24\textwidth}  
	\centering
	 \includegraphics[width=\textwidth]{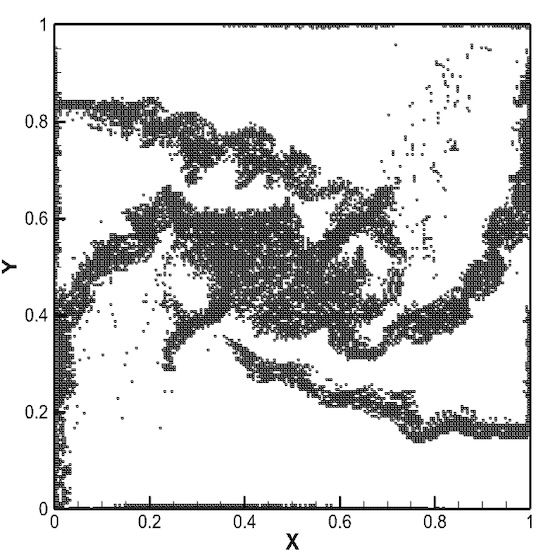}
  \end{minipage}
  \begin{minipage}[t]{0.24\textwidth}  
	\centering
	 \includegraphics[width=\textwidth]{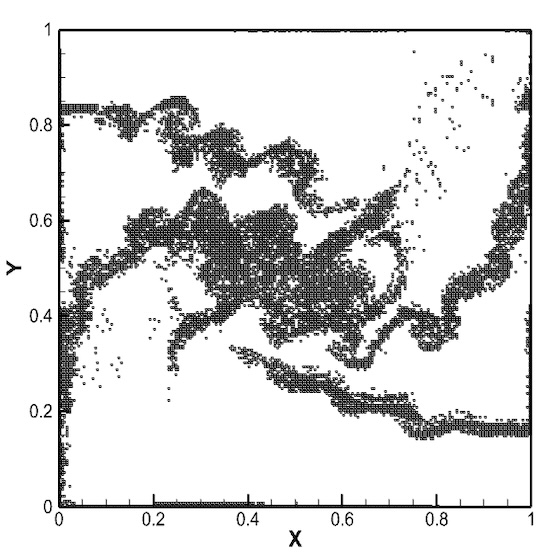}
  \end{minipage}
  }
  \caption{The troubled cells of Example 4.8 \eqref{2dr2}   at the end of time with 200$\times$200.
  From left to right: second-order($p^1$), third-order($p^2$), fourth-order($p^3$), fifth-order($p^4$).}
  \label{fig:R2KXD}
\end{figure}


{\bf Example 4.9} { The double Mach reflection problem}

This case has been extensively adopted to test the performance of numerical schemes in
the compressible flows with strong shocks \cite{DM}.
The initial condition is given by 
\begin{eqnarray}
	\begin{aligned}
\left( {\rho, U, V, p} \right) = \left\{ \begin{matrix}
{\left(8, 4.125\sqrt 3, - 4.125, 116.5 \right),~ if ~ y \ge h(x,0) ,} \\
{\left(1.4, 0, 0, 1 \right),~otherwise .} \\
\end{matrix} \right.
		\end{aligned}
	\label{dm}
\end{eqnarray}  
 Where $h(x,t)=\sqrt 3(x-\frac{1}{6})-20t$.
 A right-moving shock of  Mach $10$ is initially positioned at $(x, y) =(1/6, 0)$ with $60^{\circ}$ to the wall.
The computational domain is $[0, 3] \times [0, 0.75]$ and output time is $T=0.2$.
The reflective boundary condition is used at the wall.
The pre-shock and post-shock conditions are imposed at the rest boundaries to describe the exact motion of the shock.
The density distributions of 600$\times$150 cells with $\Delta x=\Delta y=\frac{1}{200}$   are shown in Fig. \ref{fig:DMKX} and Fig. \ref{fig:DMKXB}, 
the corresponding troubled cells detection domain at the end of time $T=0.2$ is shown in Fig. \ref{fig:DMKXD}.  
The results demonstrate   the robustness of the new DG-HGKS scheme in two-dimensional case,
and the KXRCF indicator can still accurately identify "troubled cells" in two-dimensional case.

\begin{figure}[h!]
         \centering
	  \subfloat{
	\begin{minipage}[t]{0.9\textwidth}  
	\centering
	 \includegraphics[width=\textwidth]{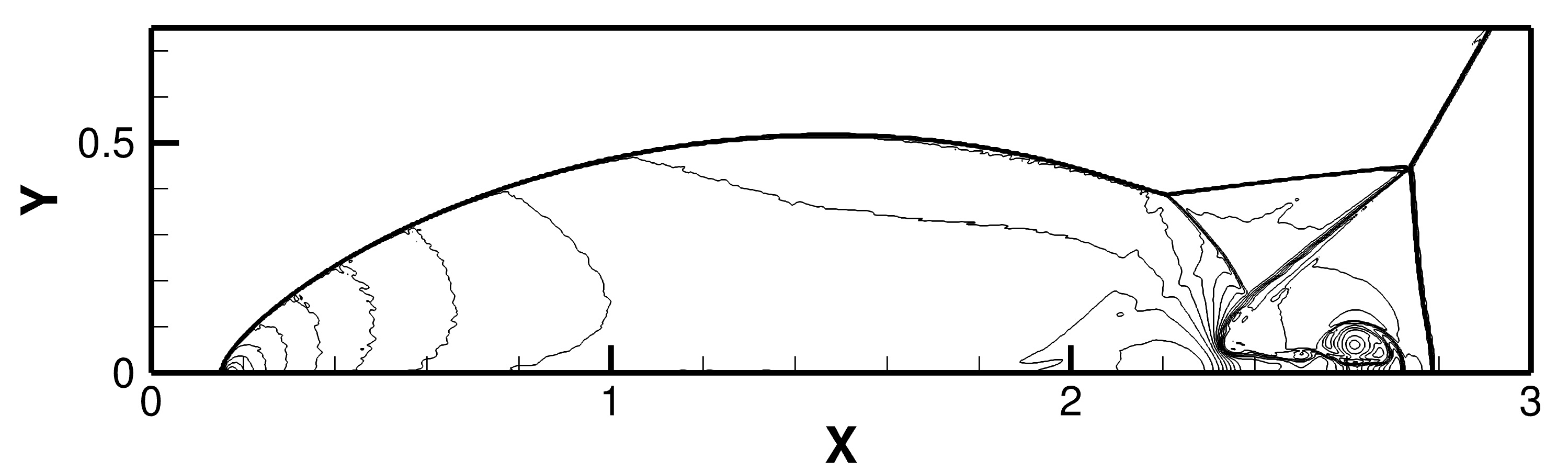}
  \end{minipage}}\\
    \subfloat{
  \begin{minipage}[t]{0.9\textwidth}  
	\centering
	 \includegraphics[width=\textwidth]{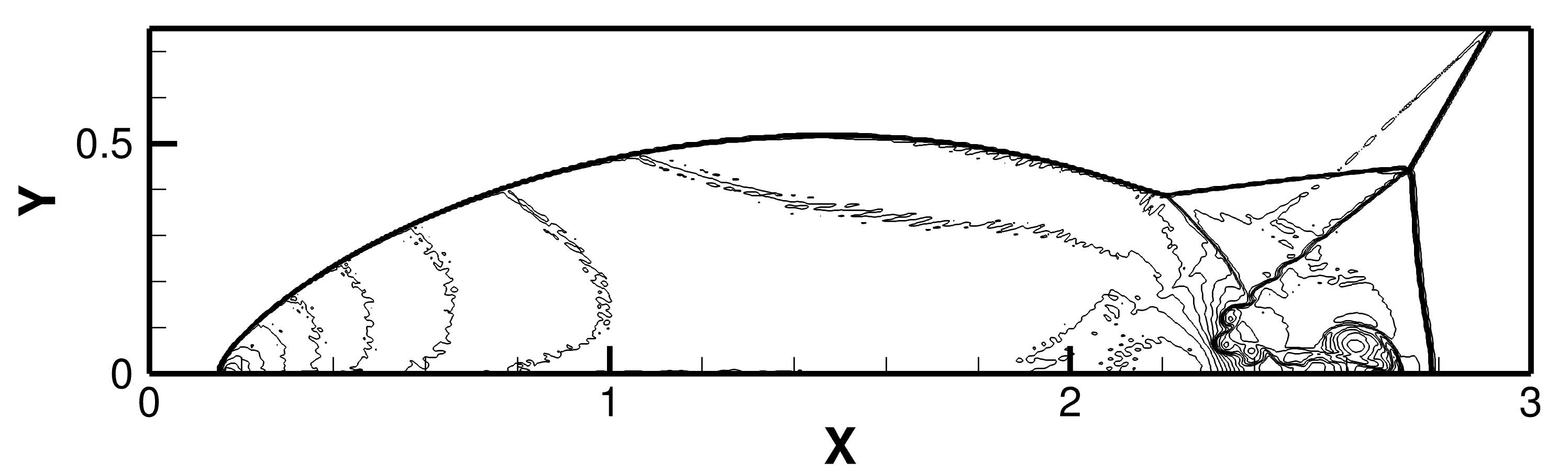}
  \end{minipage}}\\
    \subfloat{
  \begin{minipage}[t]{0.9\textwidth}  
	\centering
	 \includegraphics[width=\textwidth]{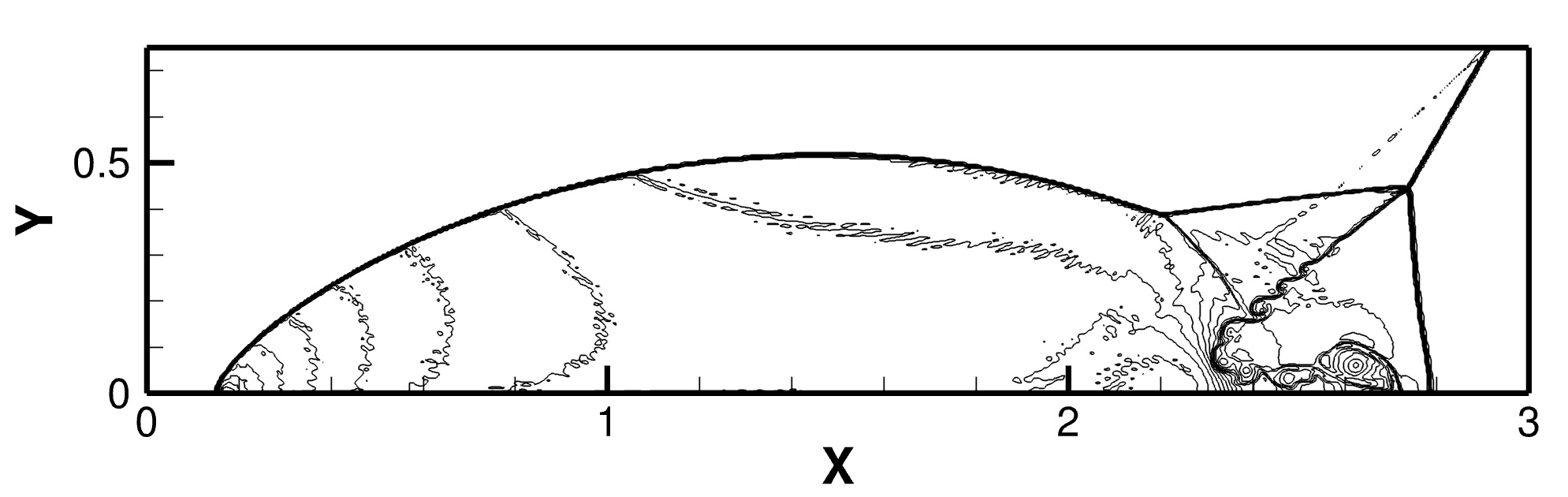}
  \end{minipage}}\\
    \subfloat{
   \begin{minipage}[t]{0.9\textwidth}  
	\centering
	 \includegraphics[width=\textwidth]{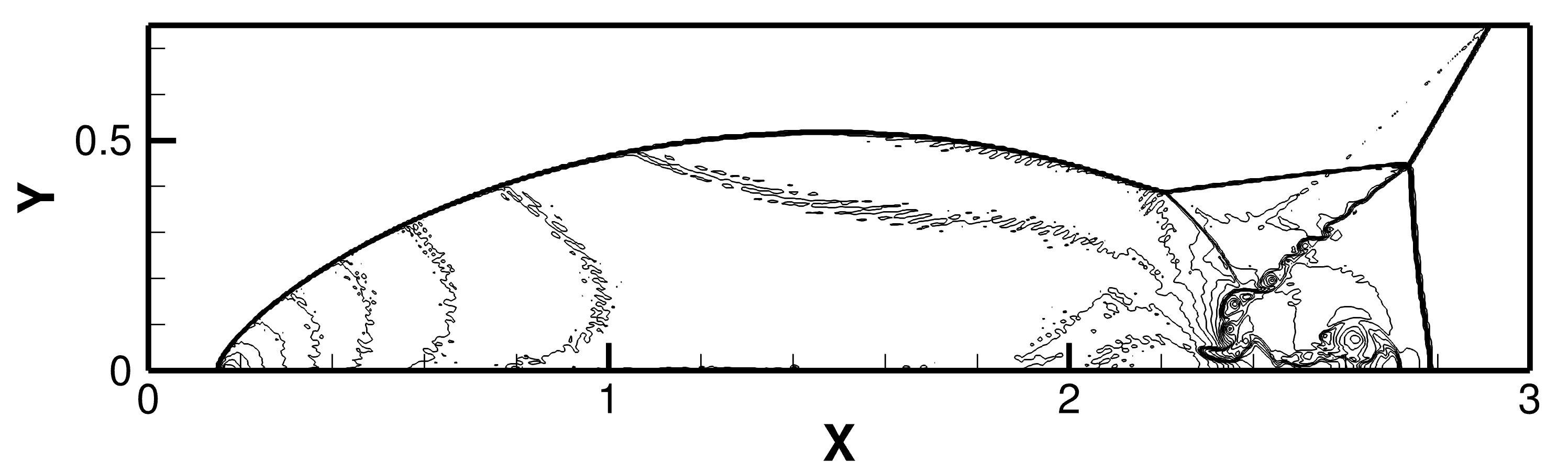}
  \end{minipage}}
 \caption{ The density distribution  of  Example 4.9 \eqref{dm} with $\Delta x=\Delta y=1/200$.
 From top to bottom: second-order($p^1$), third-order($p^2$), fourth-order($p^3$), fifth-order($p^4$).
 30 contours are drawn from 1.5 to 21.5.}
  \label{fig:DMKX}
\end{figure}
 
 \begin{figure}[h!]
         \centering
	  \subfloat{
	\begin{minipage}[t]{0.5\textwidth}  
	\centering
	 \includegraphics[width=\textwidth]{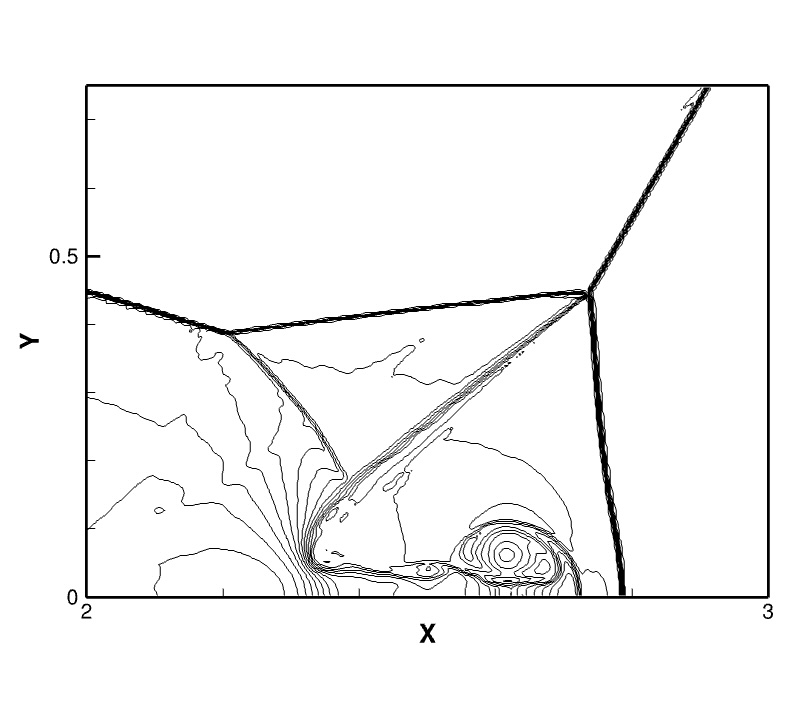}
  \end{minipage}
  \begin{minipage}[t]{0.5\textwidth}  
	\centering
	 \includegraphics[width=\textwidth]{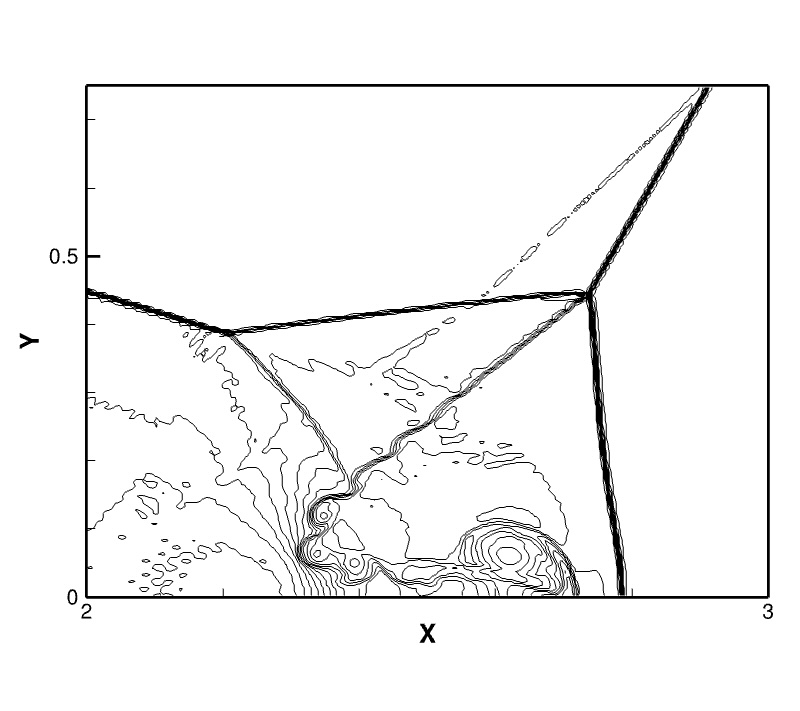}
  \end{minipage}
  }
  \\
    \subfloat{
  \begin{minipage}[t]{0.5\textwidth}  
	\centering
	 \includegraphics[width=\textwidth]{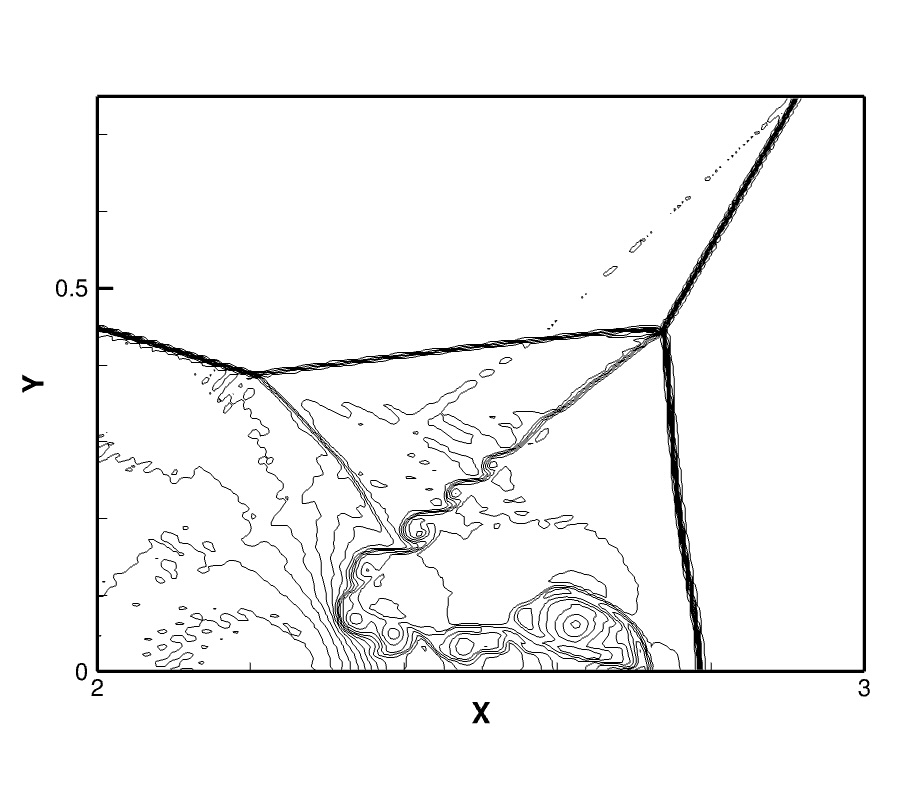}
  \end{minipage}
   \begin{minipage}[t]{0.5\textwidth}  
	\centering
	 \includegraphics[width=\textwidth]{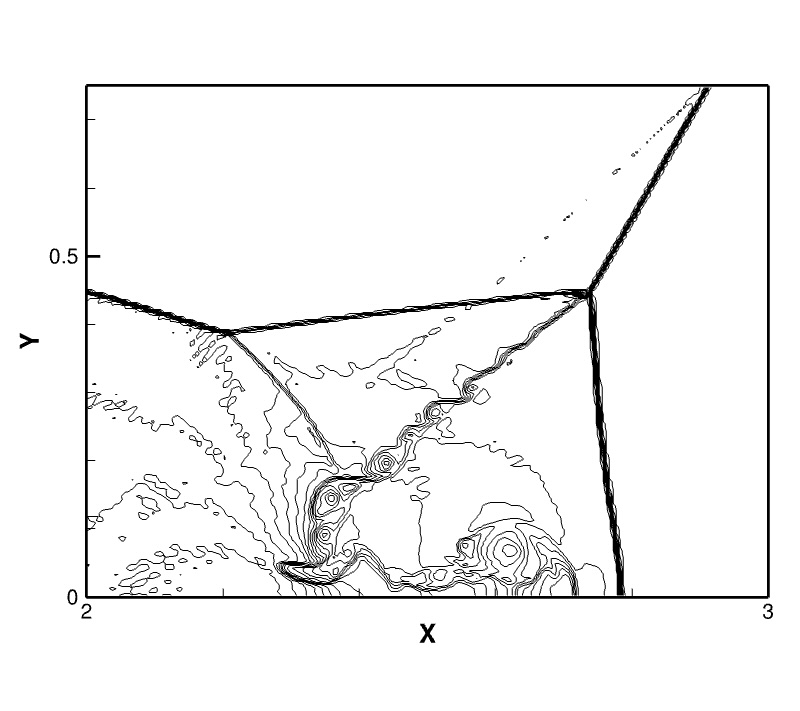}
  \end{minipage}}
 \caption{ The enlarged density distribution of Example 4.9 \eqref{dm}  with $\Delta x=\Delta y=1/200$. 
 Left top: second-order($p^1$), right top: third-order($p^2$), left bottom: fourth-order($p^3$), right bottom: fifth-order($p^4$).
 30 contours are drawn from 1.5 to 21.5.}
  \label{fig:DMKXB}
\end{figure}
 
\begin{figure}[h!]
         \centering
	  \subfloat{
	\begin{minipage}[t]{0.9\textwidth}  
	\centering
	 \includegraphics[width=\textwidth]{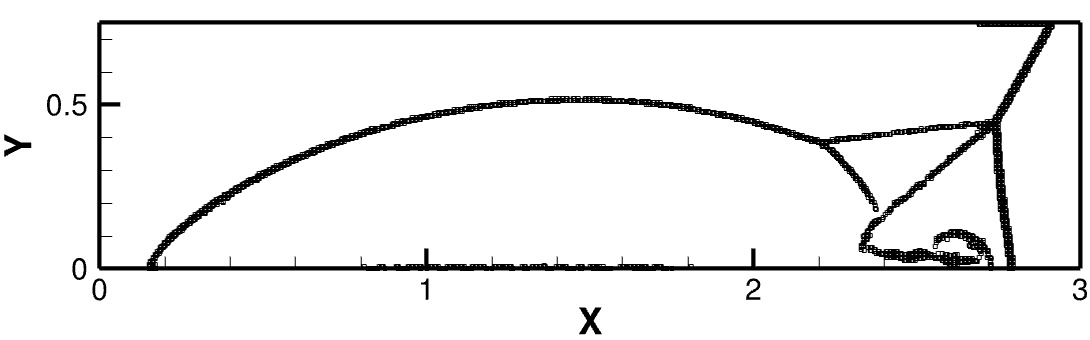}
  \end{minipage}}\\
    \subfloat{
  \begin{minipage}[t]{0.9\textwidth}  
	\centering
	 \includegraphics[width=\textwidth]{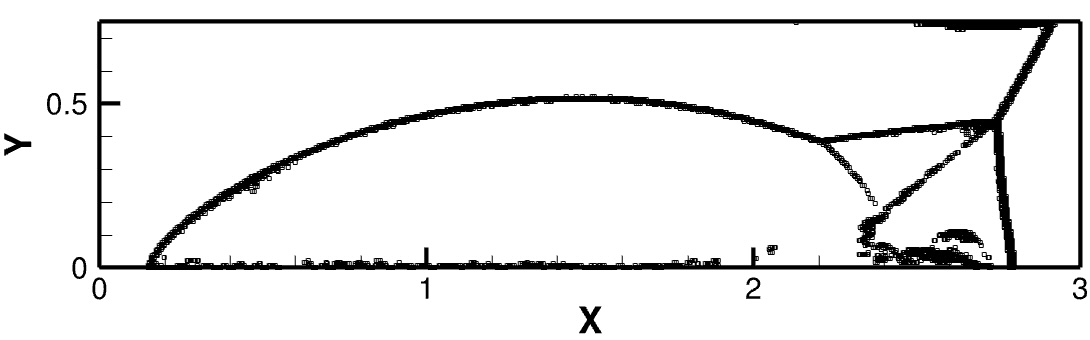}
  \end{minipage}}\\
    \subfloat{
  \begin{minipage}[t]{0.9\textwidth}  
	\centering
	 \includegraphics[width=\textwidth]{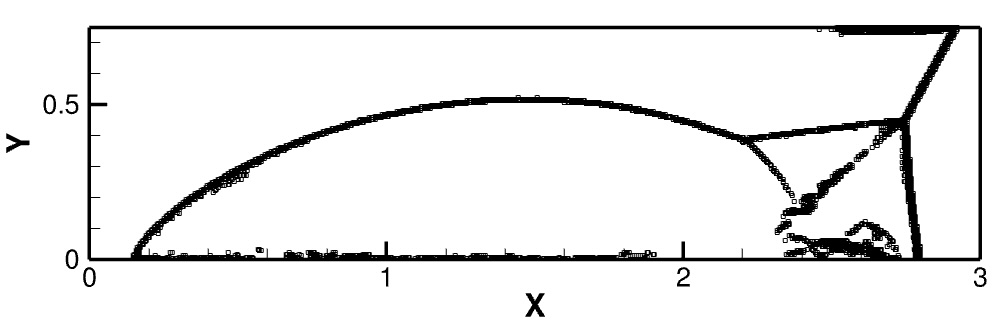}
  \end{minipage}}\\
    \subfloat{
   \begin{minipage}[t]{0.9\textwidth}  
	\centering
	 \includegraphics[width=\textwidth]{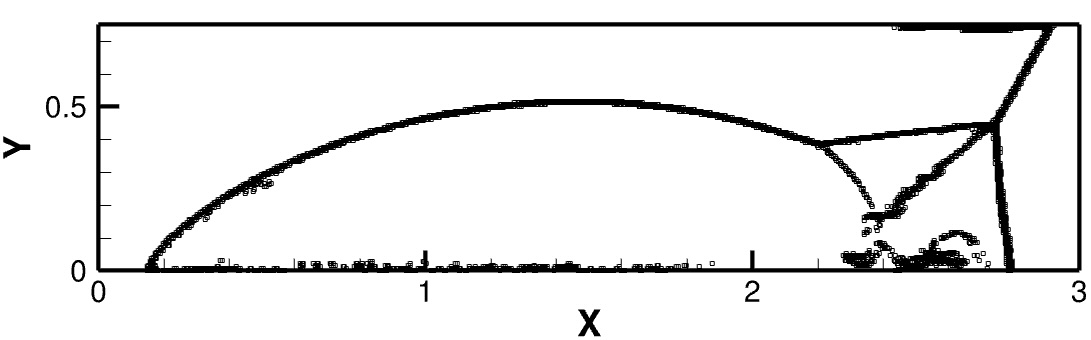}
  \end{minipage}}
 \caption{ The troubled cells of Example 4.9 \eqref{dm}   at the end of time with $\Delta x=\Delta y=1/200$.
  From top to bottom: second-order($p^1$), third-order($p^2$), fourth-order($p^3$), fifth-order($p^4$).}
  \label{fig:DMKXD}
\end{figure}

\section{Conclusions}
Inspired by the newly developed finite volume CEHGKS, a uniformly arbitrary high-order scheme, 
and SHWENO, a high-order reconstruction technique using compact stencils, 
we construct a new uniformly arbitrary high-order DG-HGKS. 
The main advantages of this scheme are as follows: 
(1) This is the first one-stage DG-HGKS scheme that can achieve higher than third-order accuracy, 
and in fact, it can achieve arbitrary high-order accuracy. 
(2) Compared with the original finite volume CEHGKS, 
the new scheme fully utilizes and maintains the compactness of the DG method, 
effectively reducing the numerical dissipation of the calculation. 
(3) Compared with the RKDG method using the same detector and limiter strategy, 
the “troubled cell” region required by the new scheme is significantly smaller, 
and the compactness of the DG method is better maintained. 
A series of 1D and 2D numerical results and comparisons with the existing CEHGKS and RKDG schemes validate the good performance of the new scheme.


\section*{Acknowledgements}

This work is supported by National Key R\&D Program of China (Grant Nos.        
2022YFA-1004500), 
NSFC (Grant Nos. 12271052, 12101063), 
and National Key Project (Grant No. GJXM92579).

\clearpage
 

%
%

\begin{appendices}
\setcounter{table}{0}   
\setcounter{figure}{0}
\setcounter{section}{0}
\setcounter{equation}{0}
\section{Comparison with  the limiter applied to all cells}\label{app:FLOPs}
\renewcommand{\thetable}{A.\arabic{table}}
\renewcommand{\thefigure}{A.\arabic{figure}}
\renewcommand{\theequation}{A.\arabic{equation}}

Here are the comparison results between applying the limiter to all cells and only to  "troubled cells".
When applying the limiter to all cells,
the CFL numbers can reach the level of RKDG method.
The numerical results include  the  Shu-Osher shock acoustic wave interaction   in 1D case,
and the  Riemann problem in 2D case.
See Fig. \ref{fig:SOA} to  Fig. \ref{fig:R2A} for details. 
It is obvious that the method of applying the limiter only to the "troubled cells" has better results.
\begin{figure}[h!]
	\centering
	\begin{minipage}[t]{0.235\textwidth}  
	\centering
	 \includegraphics[width=\textwidth]{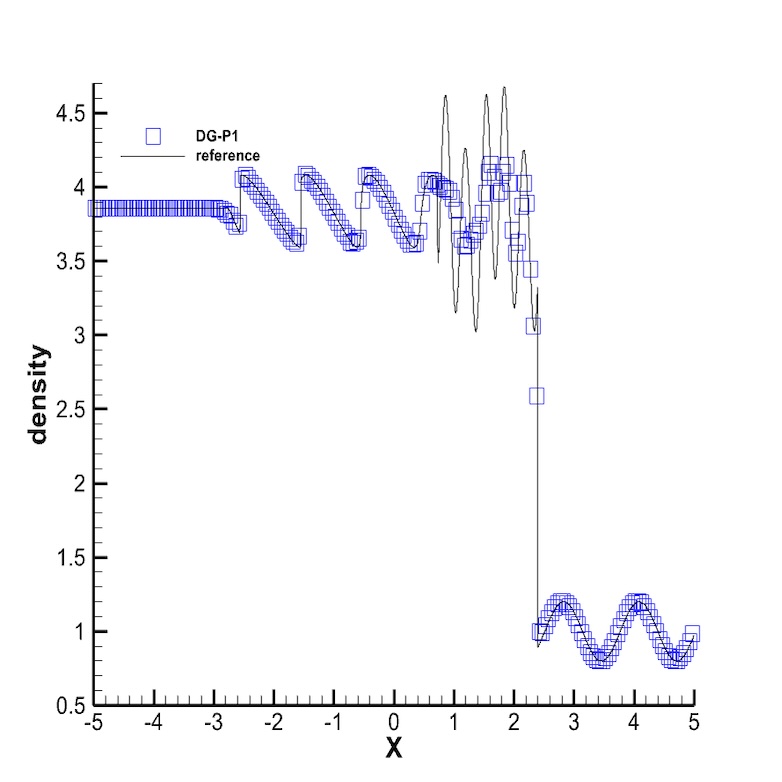}
  \end{minipage}
   \begin{minipage}[t]{0.235\textwidth}  
	\centering
	 \includegraphics[width=\textwidth]{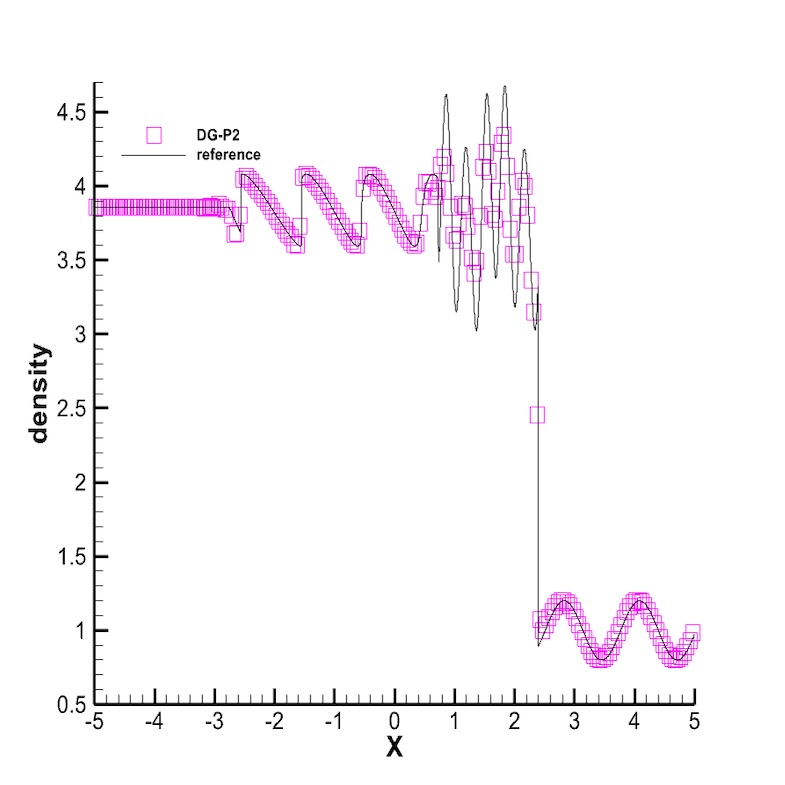}
  \end{minipage}
  \begin{minipage}[t]{0.235\textwidth}  
	\centering
	 \includegraphics[width=\textwidth]{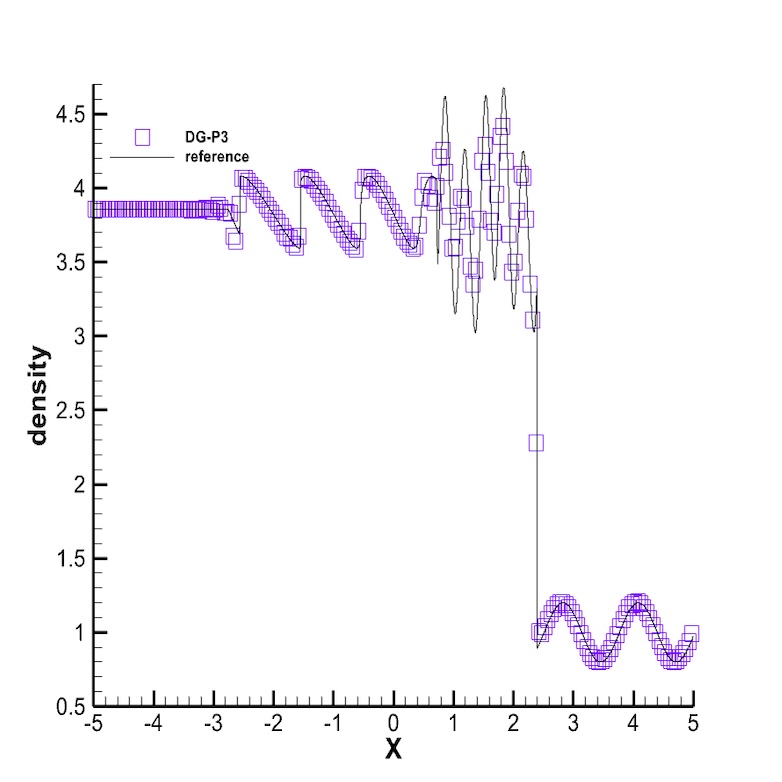}
  \end{minipage}
  \begin{minipage}[t]{0.235\textwidth}  
	\centering
	 \includegraphics[width=\textwidth]{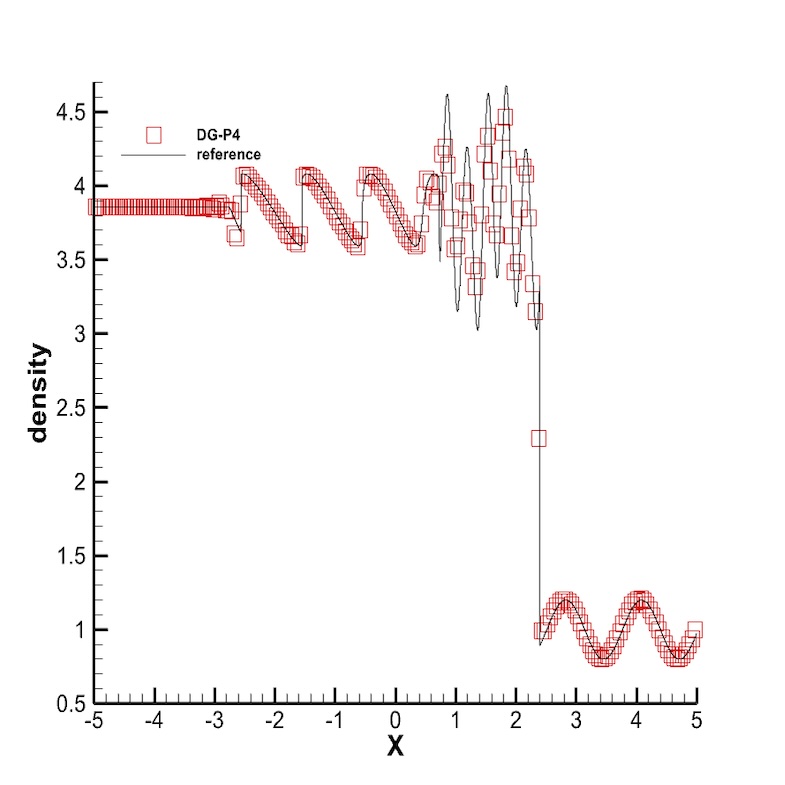}
  \end{minipage}
  \caption*{(a)  DG-HGKS, apply the limiter to all cells, 
  from left to right: second-order($p^1$), third-order($p^2$), fourth-order($p^3$), fifth-order($p^4$). }
 \begin{minipage}[t]{0.235\textwidth}  
	\centering
	 \includegraphics[width=\textwidth]{SOP1KX}
  \end{minipage}
   \begin{minipage}[t]{0.235\textwidth}  
	\centering
	 \includegraphics[width=\textwidth]{SOP2KX}
  \end{minipage}
  \begin{minipage}[t]{0.235\textwidth}  
	\centering
	 \includegraphics[width=\textwidth]{SOP3KX}
  \end{minipage}
  \begin{minipage}[t]{0.235\textwidth}  
	\centering
	 \includegraphics[width=\textwidth]{SOP4KX}
  \end{minipage}
  \caption*{(b) DG-HGKS, apply the limiter to "troubled cells", 
  from left to right: second-order($p^1$), third-order($p^2$), fourth-order($p^3$), fifth-order($p^4$).}
  \caption{The density distribution of Example 4.3 \eqref{so} with   $N=200.$}
  \label{fig:SOA}
\end{figure}

\begin{figure}[h!]
	\centering
	\begin{minipage}[t]{0.235\textwidth}  
	\centering
	 \includegraphics[width=\textwidth]{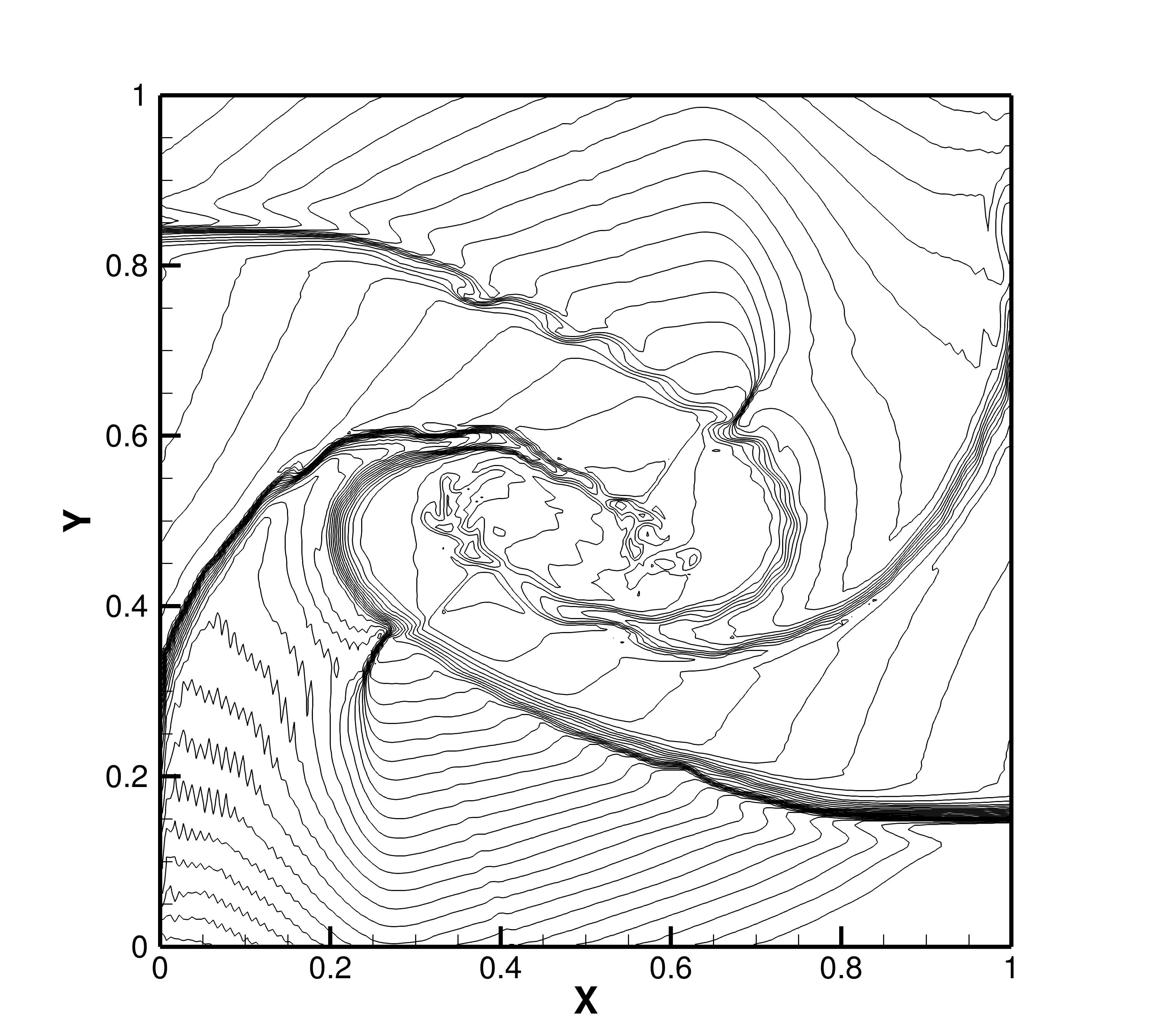}
  \end{minipage}
   \begin{minipage}[t]{0.235\textwidth}  
	\centering
	 \includegraphics[width=\textwidth]{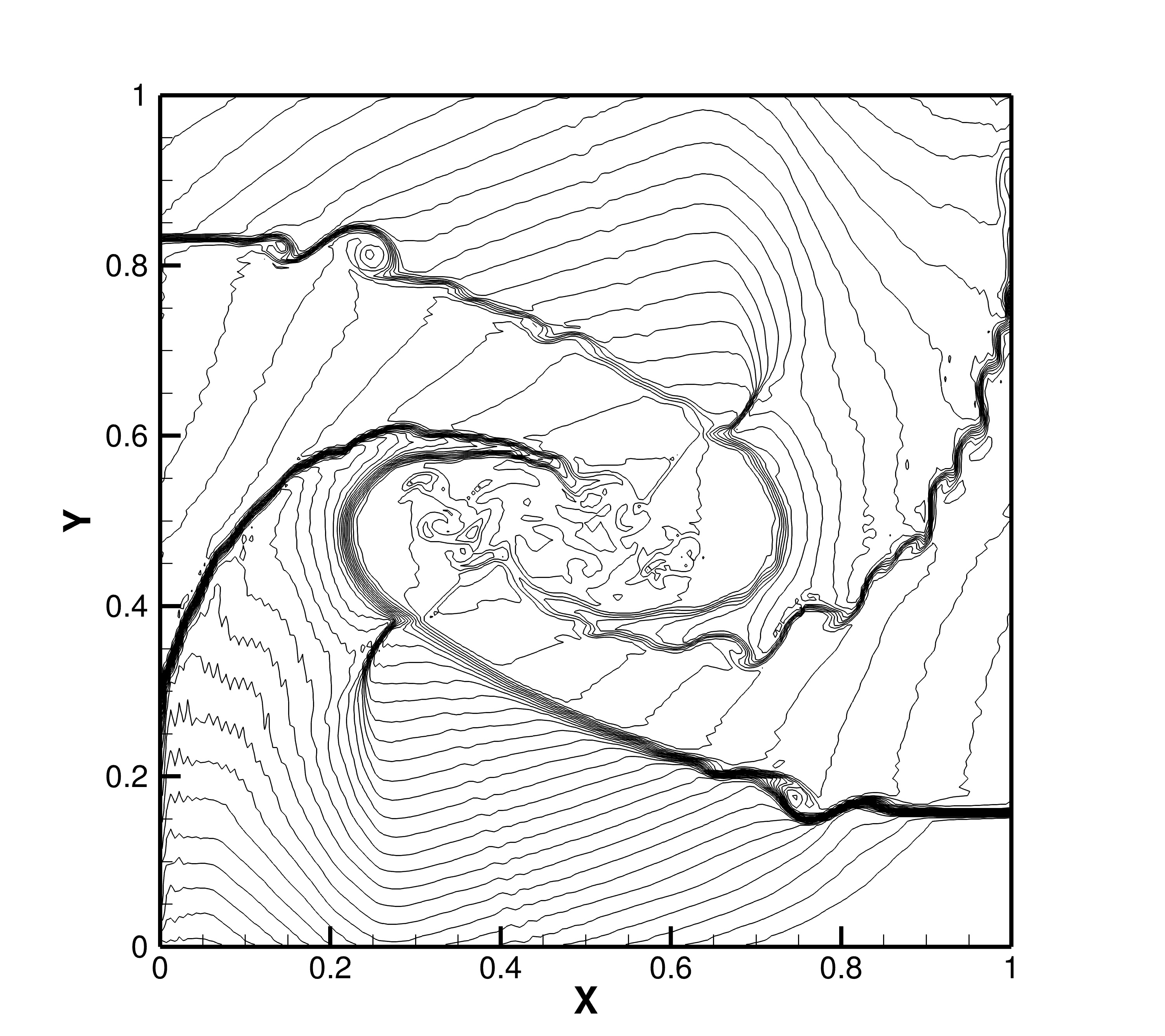}
  \end{minipage}
  \begin{minipage}[t]{0.235\textwidth}  
	\centering
	 \includegraphics[width=\textwidth]{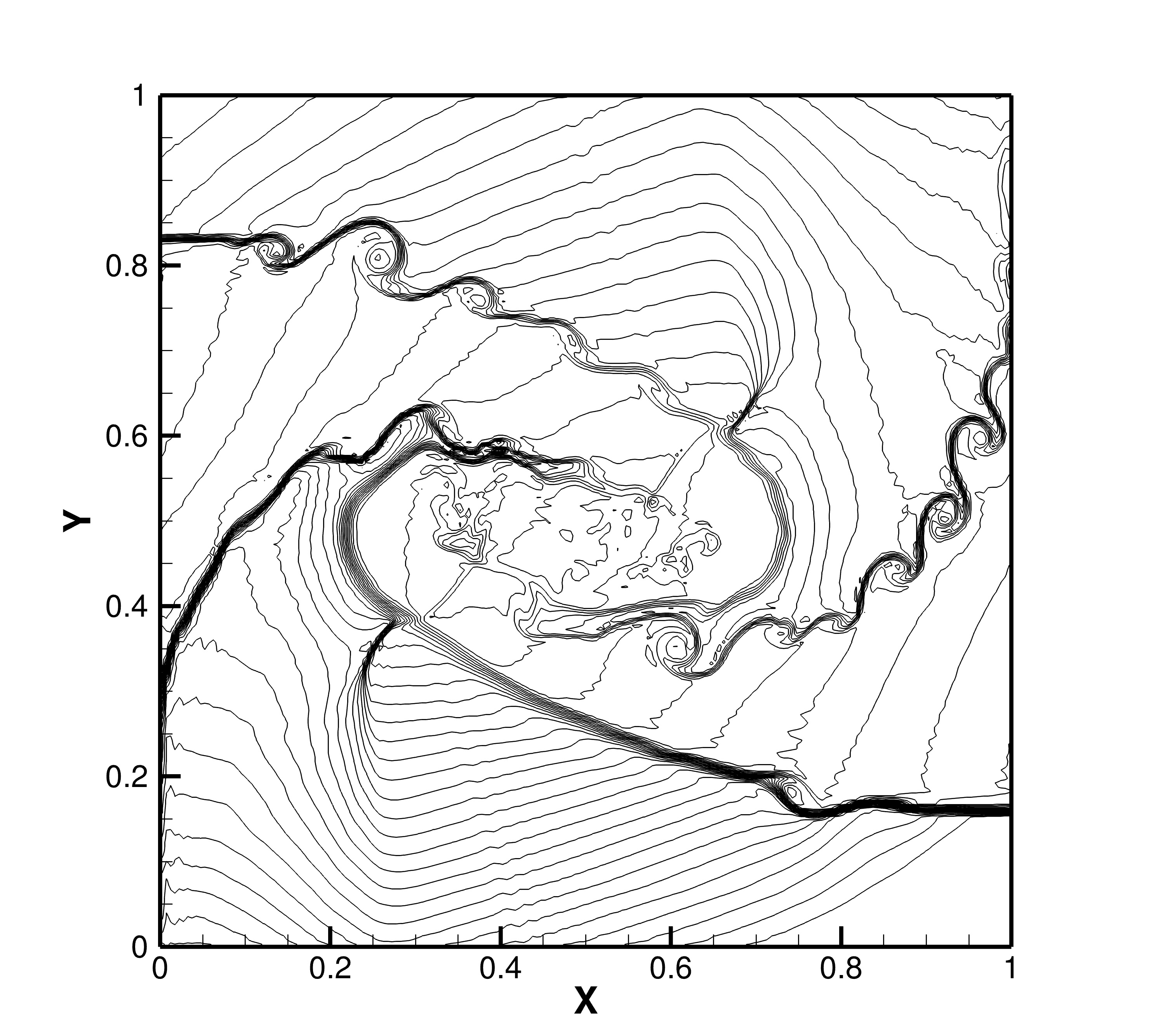}
  \end{minipage}
  \begin{minipage}[t]{0.235\textwidth}  
	\centering
	 \includegraphics[width=\textwidth]{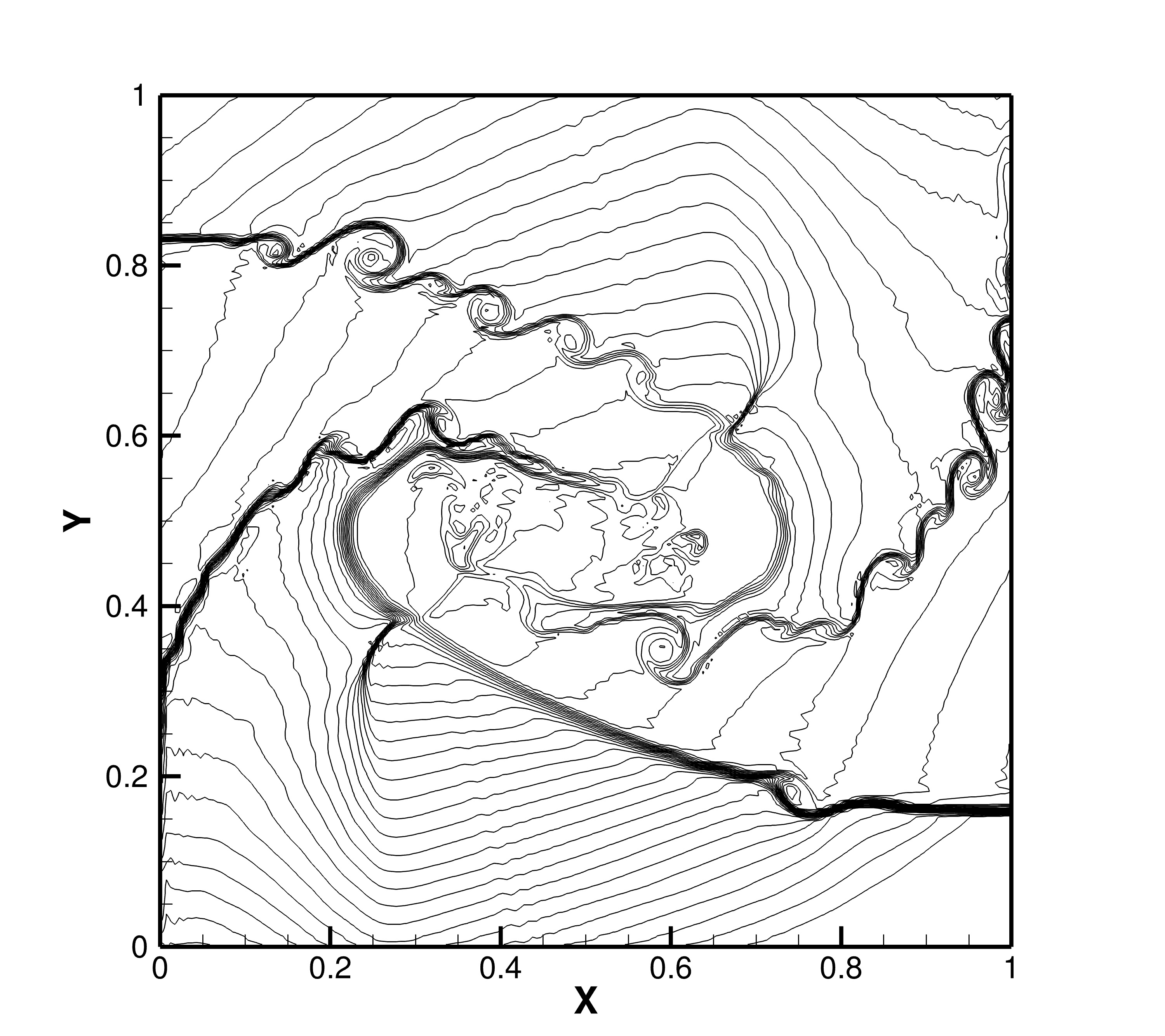}
  \end{minipage}
  \caption*{(a)  DG-HGKS, apply the limiter to all cells, 
  from left to right: second-order($p^1$), third-order($p^2$), fourth-order($p^3$), fifth-order($p^4$). }
  \begin{minipage}[t]{0.235\textwidth}  
	\centering
	 \includegraphics[width=\textwidth]{R2P1200KX}
  \end{minipage}
   \begin{minipage}[t]{0.235\textwidth}  
	\centering
	 \includegraphics[width=\textwidth]{R2P2200KX}
  \end{minipage}
  \begin{minipage}[t]{0.235\textwidth}  
	\centering
	 \includegraphics[width=\textwidth]{R2P3200KX}
  \end{minipage}
  \begin{minipage}[t]{0.235\textwidth}  
	\centering
	 \includegraphics[width=\textwidth]{R2P4200KX}
  \end{minipage}
  \caption*{(b) DG-HGKS,   apply the limiter to "troubled cells", 
   from left to right: second-order($p^1$), third-order($p^2$), fourth-order($p^3$), fifth-order($p^4$).}
  \caption{The density distribution of Example 4.8 \eqref{2dr2} with  200$\times$200.}
  \label{fig:R2A}
\end{figure}

\section{Comparison with the other method}
\renewcommand{\thetable}{B.\arabic{table}}
\renewcommand{\thefigure}{B.\arabic{figure}}
\renewcommand{\theequation}{B.\arabic{equation}}
\subsection*{B.1  Comparison with  the CEHGKS}

Here we give some comparison results about the  CEHGKS  and the fifth-order($p^4$) DG-HGKS scheme.
Including the 1D Shu-Osher shock acoustic wave interaction 
and the 2D Riemann problem.
See Fig. \ref{fig:socomp} to Fig. \ref{fig:R3comp} for details. 
It can be seen that under the same computational cells,
the DG-HGKS scheme has better numerical resolution than the CEHGKS scheme.
This suggests that the DG-HGKS  developed in this paper is more compact than the CEHGKS.

\begin{figure}[h!]
         \centering
	  \subfloat{
	\begin{minipage}[t]{0.48\textwidth}  
	\centering
	 \includegraphics[width=\textwidth]{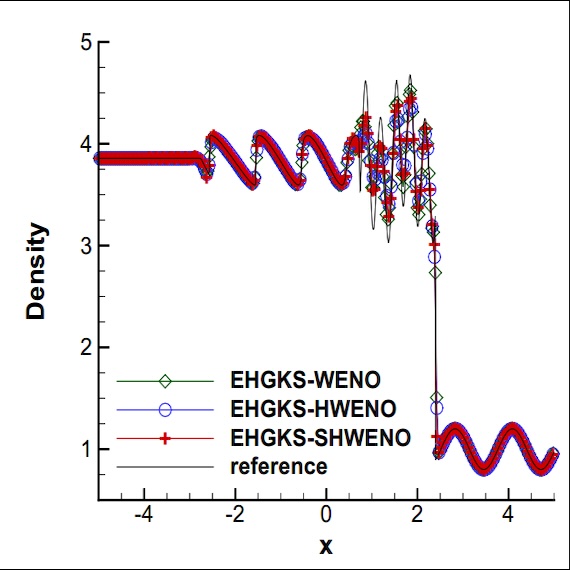}
	 \caption*{(a) CEHGKS, from\cite{LSY2}, the red line}
  \end{minipage}}
    \subfloat{
   \begin{minipage}[t]{0.48\textwidth}  
	\centering
	 \includegraphics[width=\textwidth]{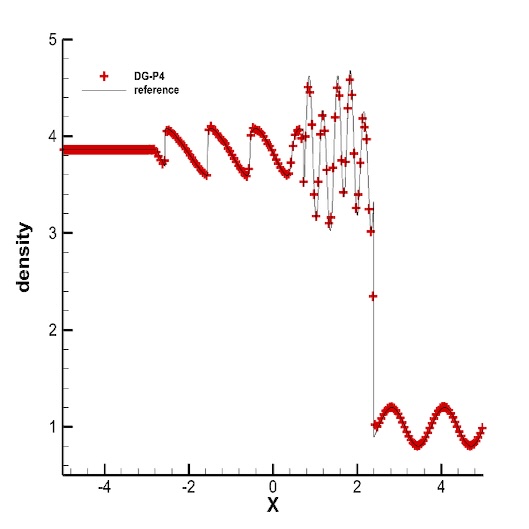}
	  \caption*{(b) DG-HGKS, fifth-order($p^4$), {\color{blue} only apply limiter to "troubled cells".} }
  \end{minipage}}
   \caption{ The density distribution  of  Example 4.3 \eqref{so} with  $N=200$.}
  \label{fig:socomp}
\end{figure}



\begin{figure}[h!]
	\centering
	\subfloat{
	\begin{minipage}[t]{0.5\textwidth}  
	\centering
	 \includegraphics[width=\textwidth]{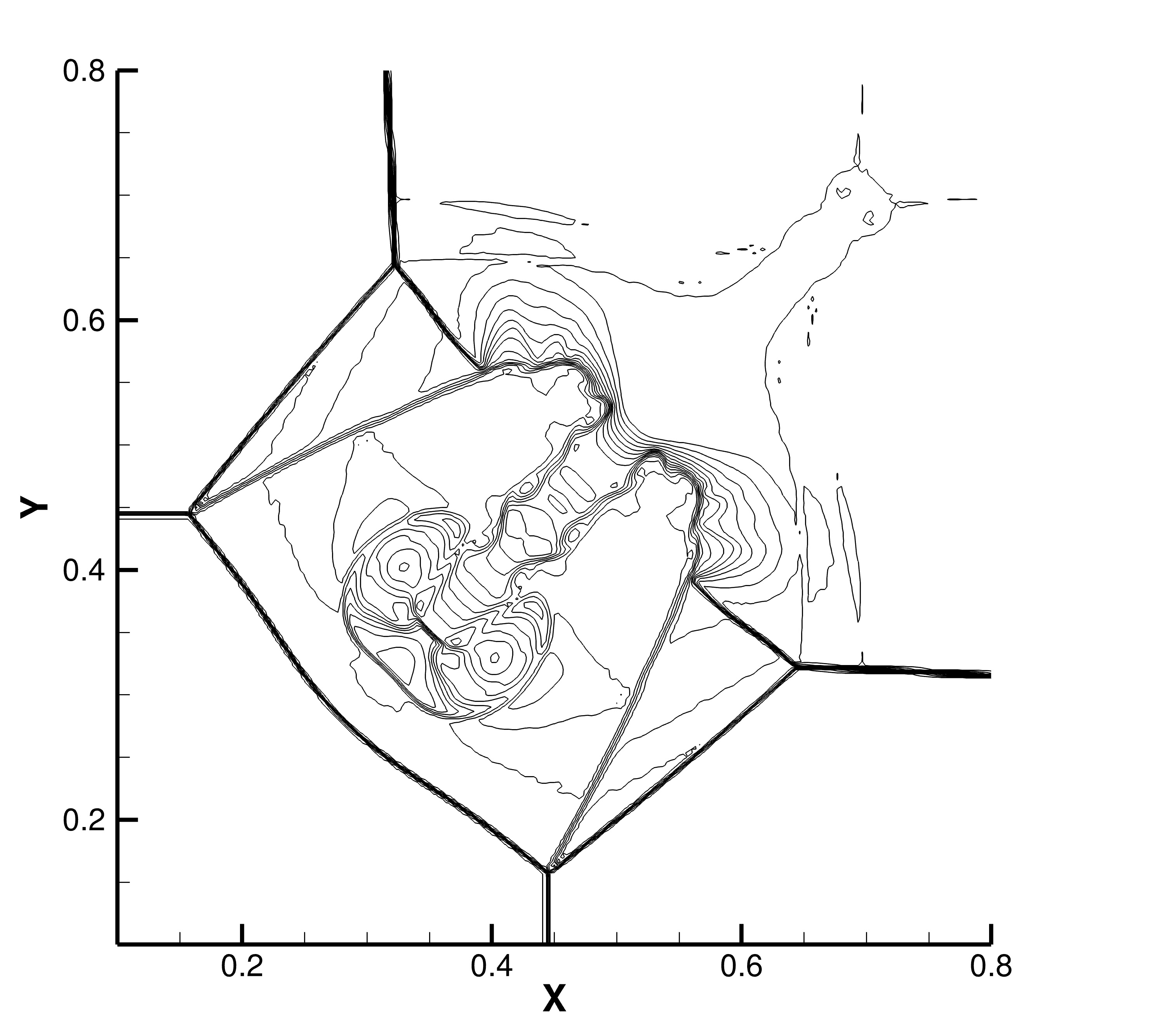}
	 \caption*{(a) CEHGKS(EHGKS-SHWENO)}
  \end{minipage}
   \begin{minipage}[t]{0.5\textwidth}  
	\centering
	 \includegraphics[width=\textwidth]{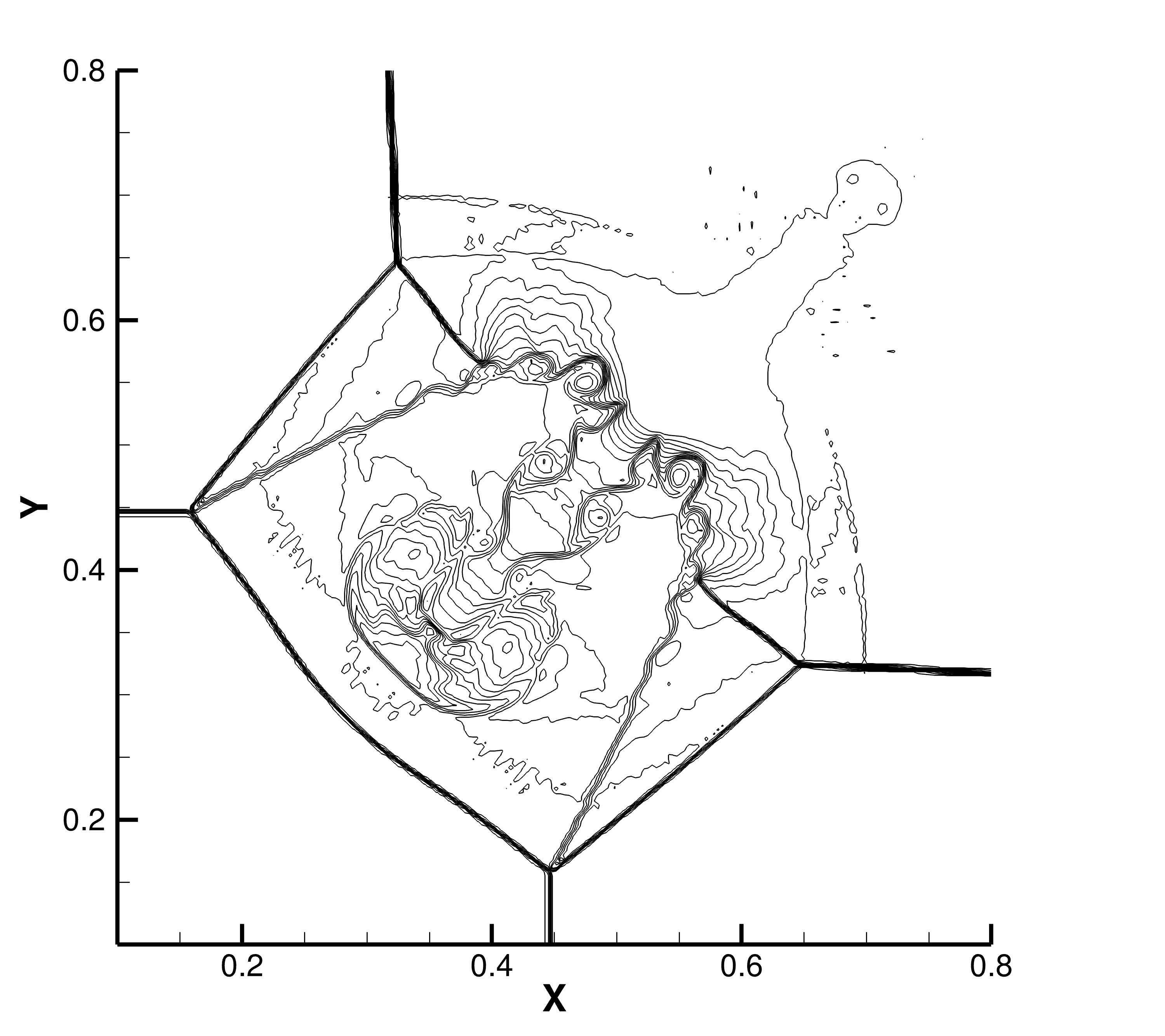}
	 \caption*{(b) DG-HGKS, fifth-order($p^4$),  only apply limiter to "troubled cells".}
  \end{minipage} }
  \caption{The density distribution of Example 4.8 \eqref{2dr3}   with 300 $\times$ 300 cells.}
  \label{fig:R3comp}
\end{figure}  

\subsection*{B.2  Comparison with  the existing RKDG}
This appendix  provides  a comparison with the RKDG  with the HWENO limiter,
which utilize the same KXRCF indicator and the same stentil \cite{SHWENOL} strategy.
The Fig. \ref{fig:dg-bw} to Fig. \ref{fig:dg-dm-kxrcf} give the results for the Woodward-Colella blast wave   and the double Mach reflection problem.
These results demonstrate that our method can detect fewer “troubled cells”,
and at the same order of accuracy, our numerical results are  more accurate.

\begin{figure}[h!]
         \centering
	\begin{minipage}[t]{0.81\textwidth}  
	\centering
	 \includegraphics[width=\textwidth]{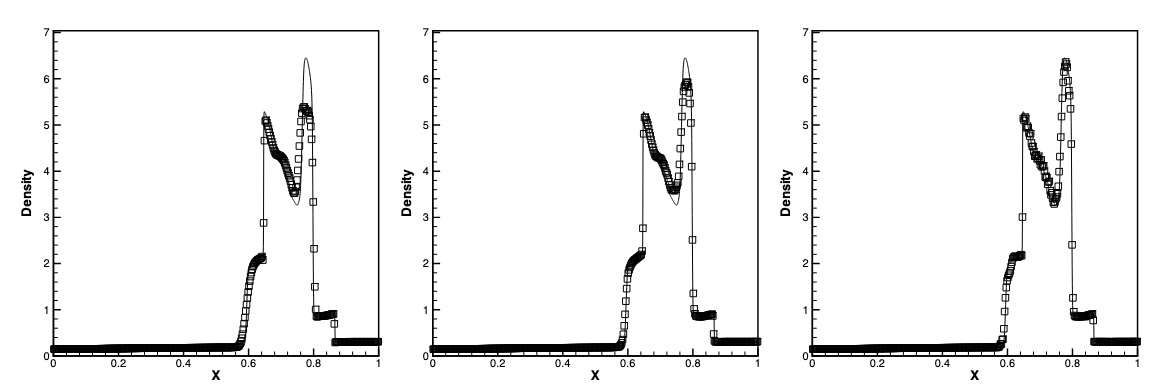} 
  \end{minipage} 
  \caption*{(a) RKDG, from\cite{SHWENOL}, 
  from left to right: second-order($p^1$), third-order($p^2$), fourth-order($p^3$).}
	\begin{minipage}[t]{0.27\textwidth}  
	\centering
	 \includegraphics[width=\textwidth]{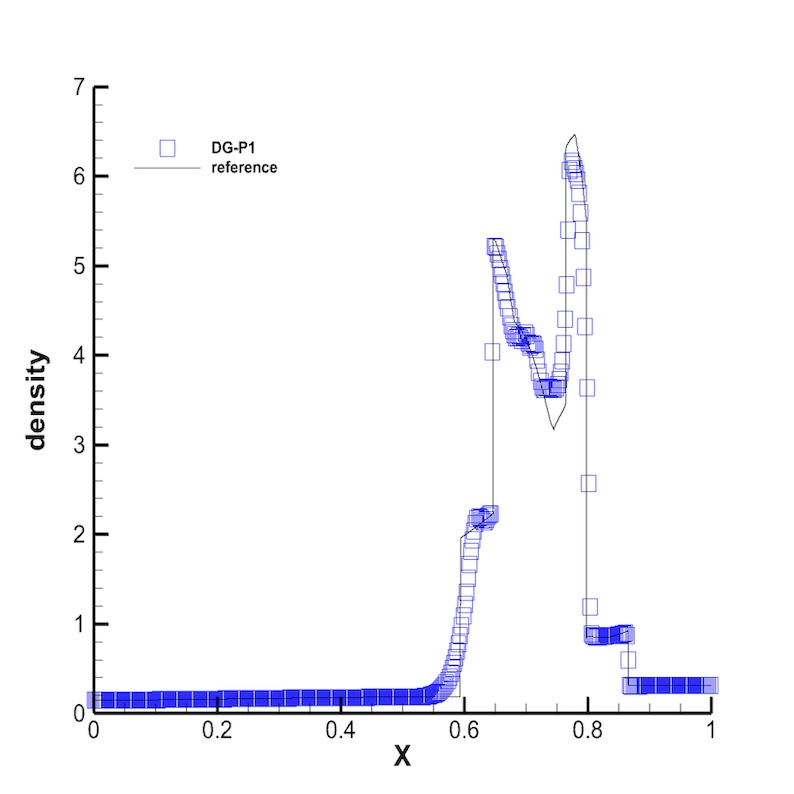}
  \end{minipage}
   \begin{minipage}[t]{0.27\textwidth}  
	\centering
	 \includegraphics[width=\textwidth]{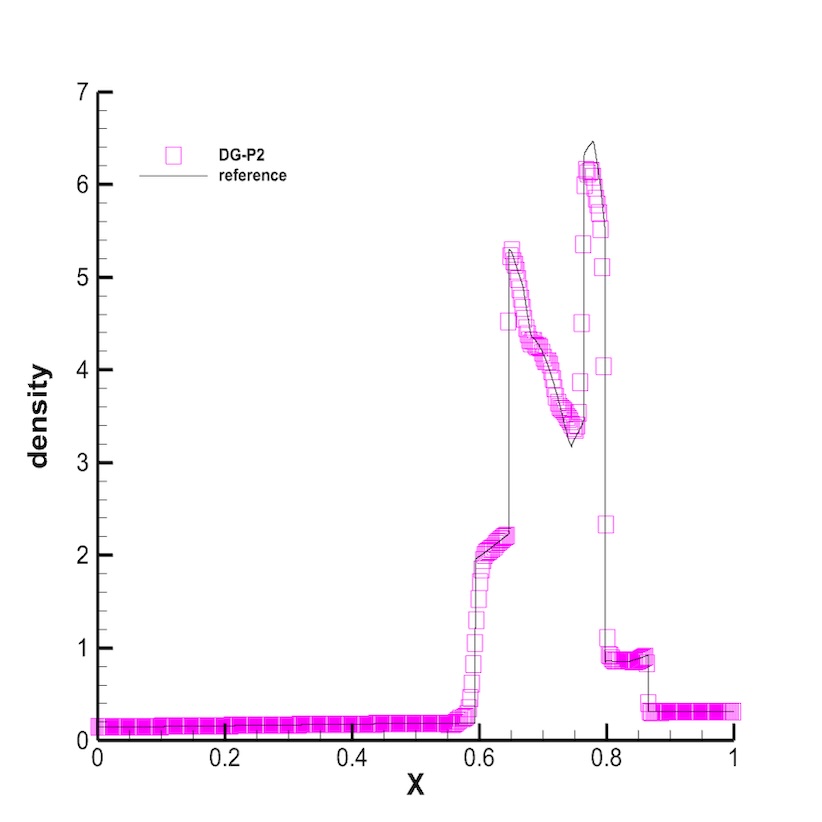}
  \end{minipage}
  \begin{minipage}[t]{0.27\textwidth}  
	\centering
	 \includegraphics[width=\textwidth]{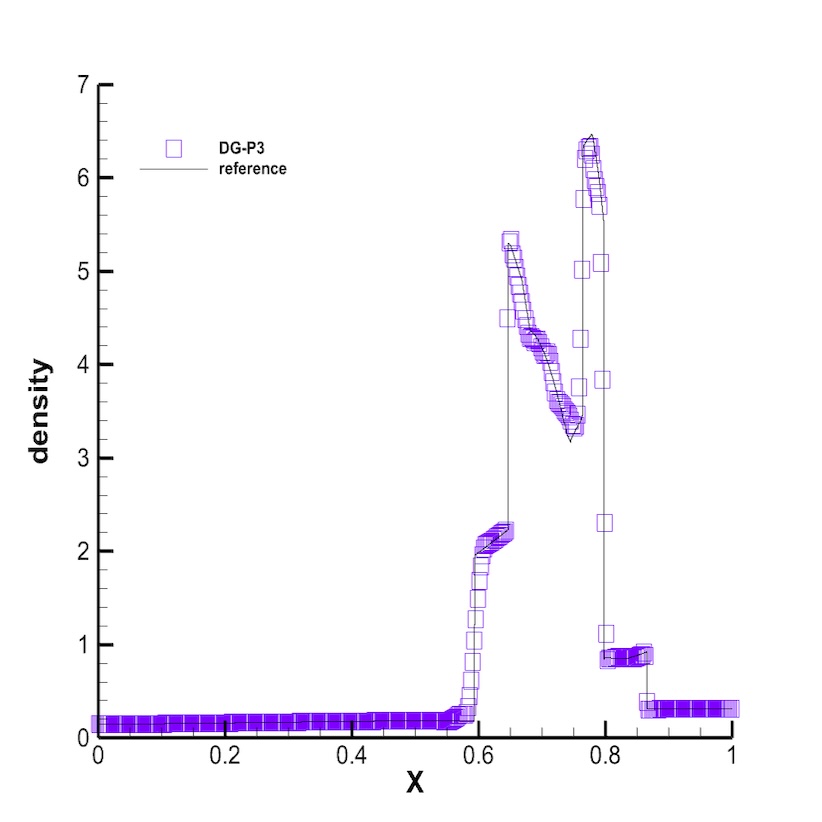}
  \end{minipage}
   \caption*{(b) DG-HGKS, from left to right: second-order($p^1$), third-order($p^2$), fourth-order($p^3$).}
     \caption{The density distribution of Example 4.4 \eqref{bw}  with   $N=400.$}
  \label{fig:dg-bw}
\end{figure}

\begin{figure}[h!]
         \centering
	\begin{minipage}[t]{0.76\textwidth}  
	\centering
	 \includegraphics[width=\textwidth]{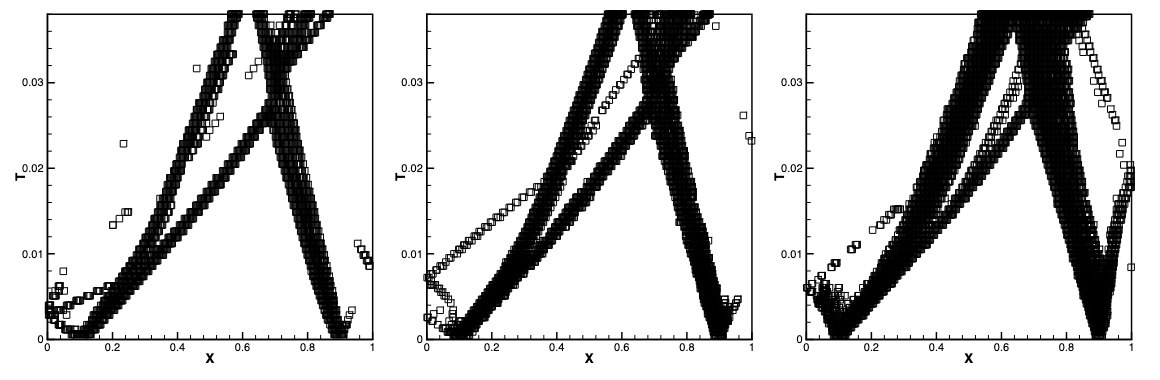} 
  \end{minipage} 
  \caption*{(a) RKDG, from\cite{SHWENOL},
  from left to right: second-order($p^1$), third-order($p^2$), fourth-order($p^3$).}
	\begin{minipage}[t]{0.25\textwidth}  
	\centering
	 \includegraphics[width=\textwidth]{BWP1KXD}
  \end{minipage}
   \begin{minipage}[t]{0.25\textwidth}  
	\centering
	 \includegraphics[width=\textwidth]{BWP2KXD}
  \end{minipage}
  \begin{minipage}[t]{0.25\textwidth}  
	\centering
	 \includegraphics[width=\textwidth]{BWP3KXD}
  \end{minipage}
   \caption*{(b) DG-HGKS,
   from left to right: second-order($p^1$), third-order($p^2$), fourth-order($p^3$).}
     \caption{The troubled cells of Example 4.4 \eqref{bw}  with   $N=400.$}
  \label{fig:dg-bw-kxrcf}
\end{figure}

\begin{figure}[h!]
         \centering
     \begin{minipage}[t]{0.32\textwidth}  
	\centering
	 \includegraphics[width=\textwidth]{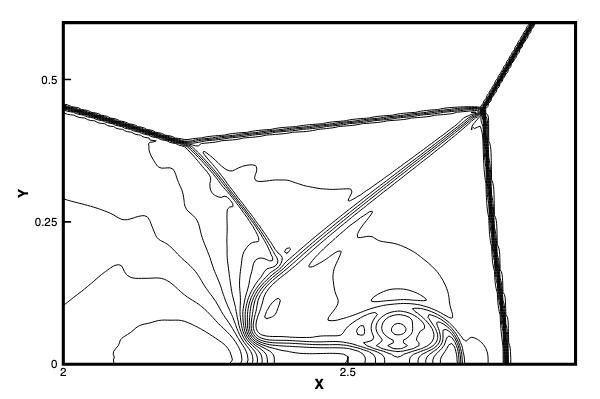}
    \end{minipage}
   \begin{minipage}[t]{0.32\textwidth}  
	\centering
	 \includegraphics[width=\textwidth]{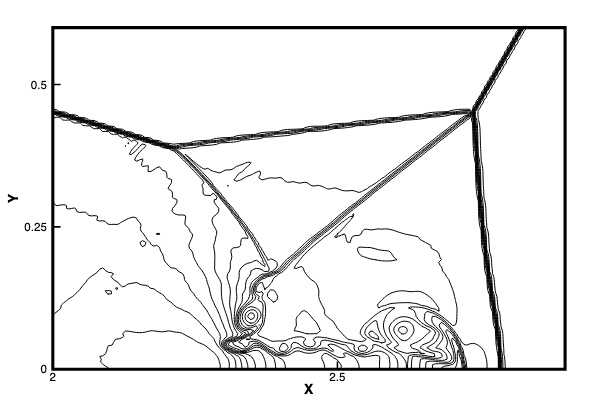}
   \end{minipage}
   \begin{minipage}[t]{0.32\textwidth}  
	\centering
	 \includegraphics[width=\textwidth]{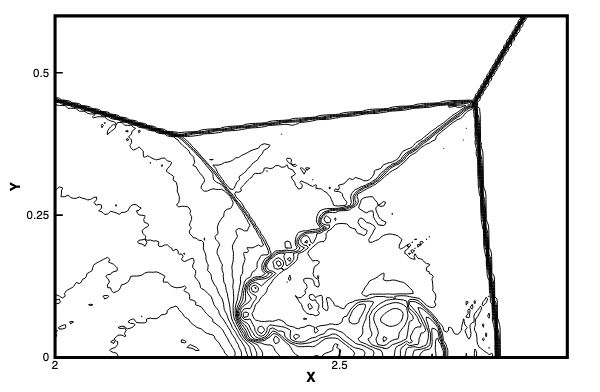}
   \end{minipage}
   \caption*{(a) RKDG, from\cite{SHWENOL},
   from left to right: second-order($p^1$), third-order($p^2$), fourth-order($p^3$).}
	\begin{minipage}[t]{0.32\textwidth}  
	\centering
	 \includegraphics[width=\textwidth]{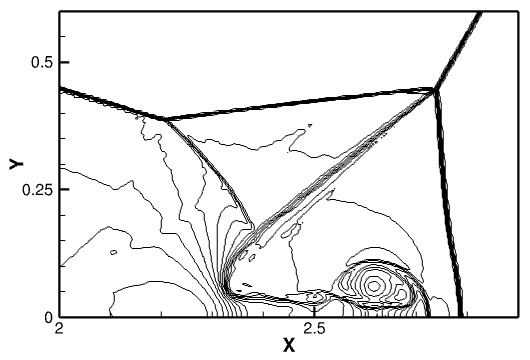}
  \end{minipage}
   \begin{minipage}[t]{0.32\textwidth}  
	\centering
	 \includegraphics[width=\textwidth]{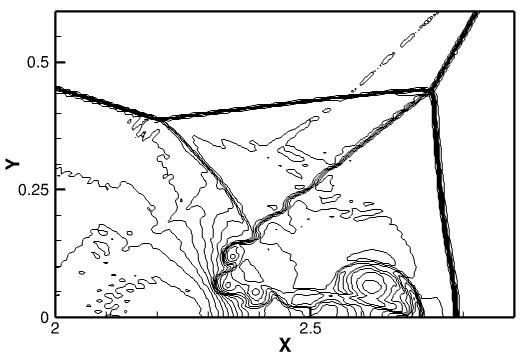}
  \end{minipage}
  \begin{minipage}[t]{0.32\textwidth}  
	\centering
	 \includegraphics[width=\textwidth]{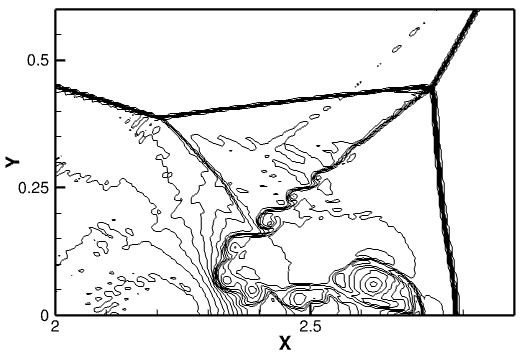}
  \end{minipage}
   \caption*{(b) DG-HGKS, from left to right: second-order($p^1$), third-order($p^2$), fourth-order($p^3$).}
     \caption{The enlarged density distribution of Example 4.9 \eqref{dm} with $\Delta x=\Delta y=1/200$.
     30 contours are drawn from 1.5 to 21.5.}
  \label{fig:dg-dm}
\end{figure}

\begin{figure}[h!]
         \centering
	  \subfloat{
	\begin{minipage}[t]{0.5\textwidth}  
	\centering
	 \includegraphics[width=\textwidth]{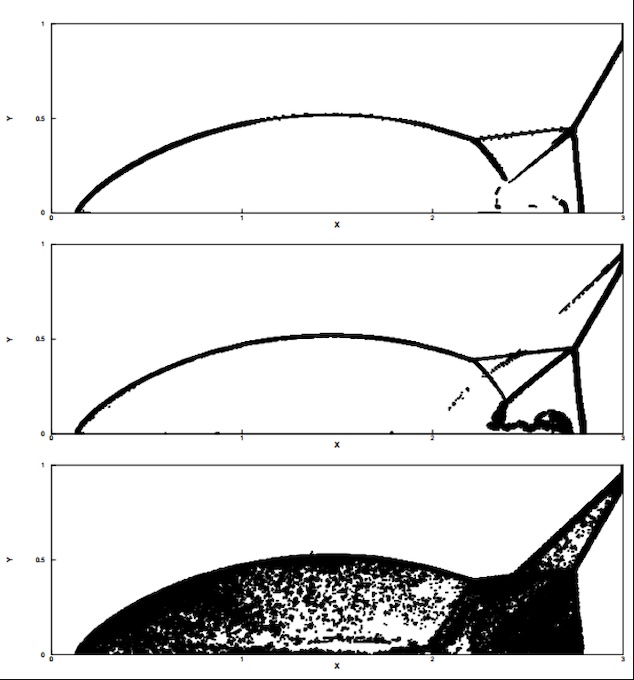}
	 \caption*{(a) RKDG, from\cite{SHWENOL},
	 from top to bottom: second-order($p^1$), third-order($p^2$), fourth-order($p^3$).}
  \end{minipage}}
    \subfloat{
   \begin{minipage}[t]{0.5\textwidth}  
	\centering
	 \includegraphics[width=\textwidth]{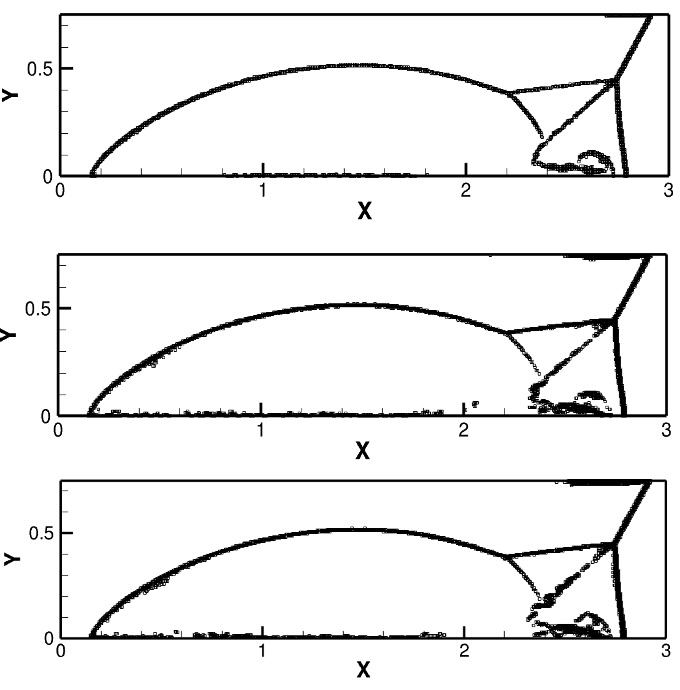}
	  \caption*{(b) DG-HGKS, 
	   from top to bottom: second-order($p^1$), third-order($p^2$), fourth-order($p^3$). }
  \end{minipage}}
   \caption{ The troubled cells  of  Example 4.9 \eqref{dm} at the end of time with $\Delta x=\Delta y=1/200$.}
  \label{fig:dg-dm-kxrcf}
\end{figure}

\subsection*{B.3  Comparison with  the RKDG-GKS}
To better demonstrate the advantages of the single-stage DG-HGKS method proposed in this paper over the multi-stage RKDG method, 
and to ensure a fair comparison, 
this appendix compares it with the RK3-DG-$p^4$ method, which is also based on the gas-kinetic flux, 
in terms of both numerical stability and computational efficiency.

Fig. \ref{fig:rkcomp} shows a comparison of numerical results for 
1D Woodward-Colella blast wave \eqref{bw}
and the modified shock/turbulence interaction \eqref{gp}
under the same computational conditions. 
The numerical results show that the single-stage method has better computational performance than the multi-stage method.

Regarding computational efficiency, 
we present the $L^2$ error as a function of CPU time for both the RKDG and DG-HGKS methods in Fig. \ref{fig:rkcpu}, which pertains to the one-dimensional advection of density perturbations. 
Additionally, the Table \ref{tab:RKDOF} provides a comparison of the computational times for both methods when applied to the same number of cells, which measured in terms of efficiency per degree of freedom (DOF). The results indicate that the single-stage method exhibits superior computational efficiency.

\begin{figure}[h!]
	\centering
	\subfloat{
	\begin{minipage}[t]{0.5\textwidth}  
	\centering
	 \includegraphics[width=\textwidth]{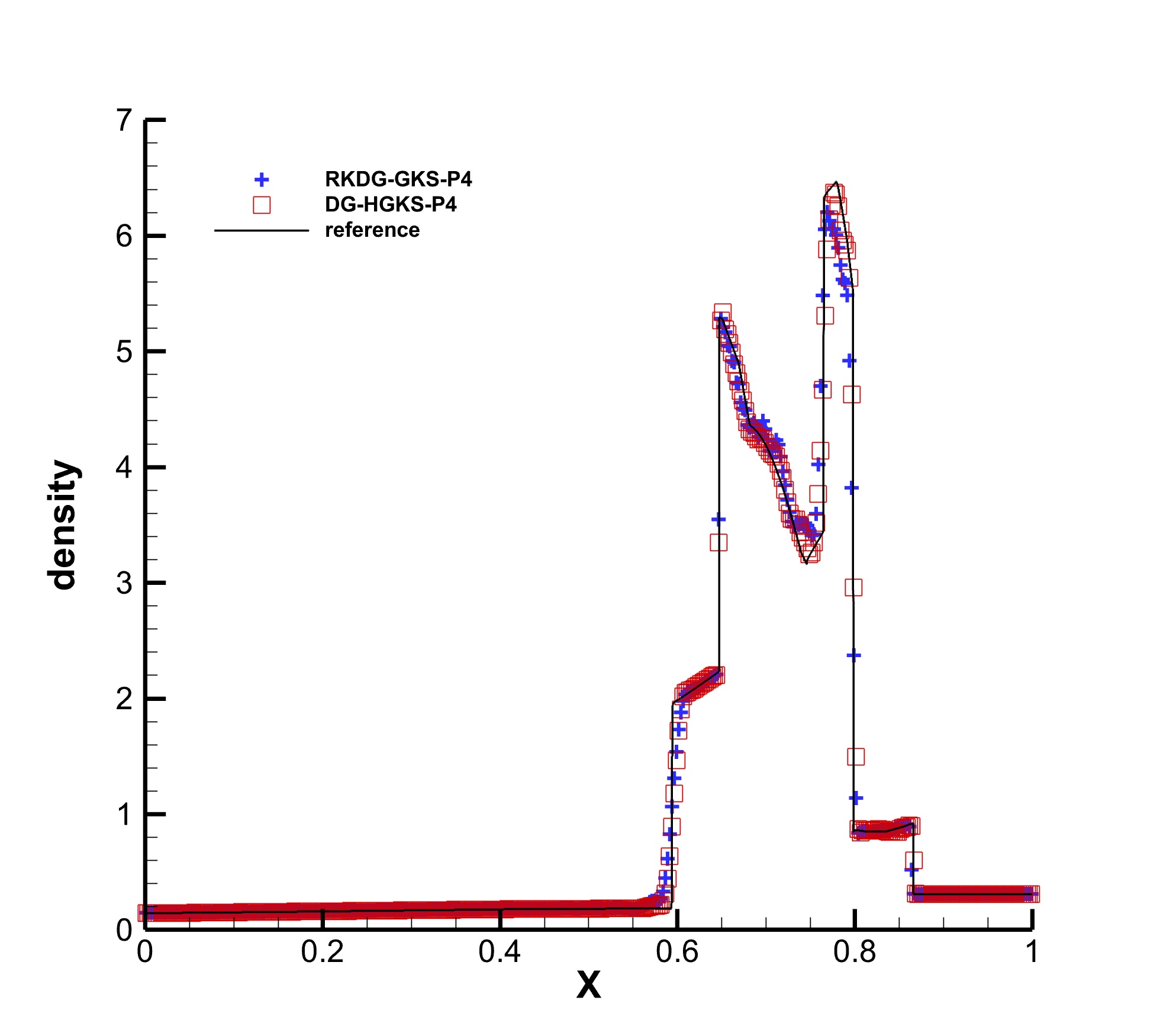}
	 \caption*{(a) {The density distribution  of Example 4.4 \eqref{bw} with $N=400.$}}
  \end{minipage}
   \begin{minipage}[t]{0.5\textwidth}  
	\centering
	 \includegraphics[width=\textwidth]{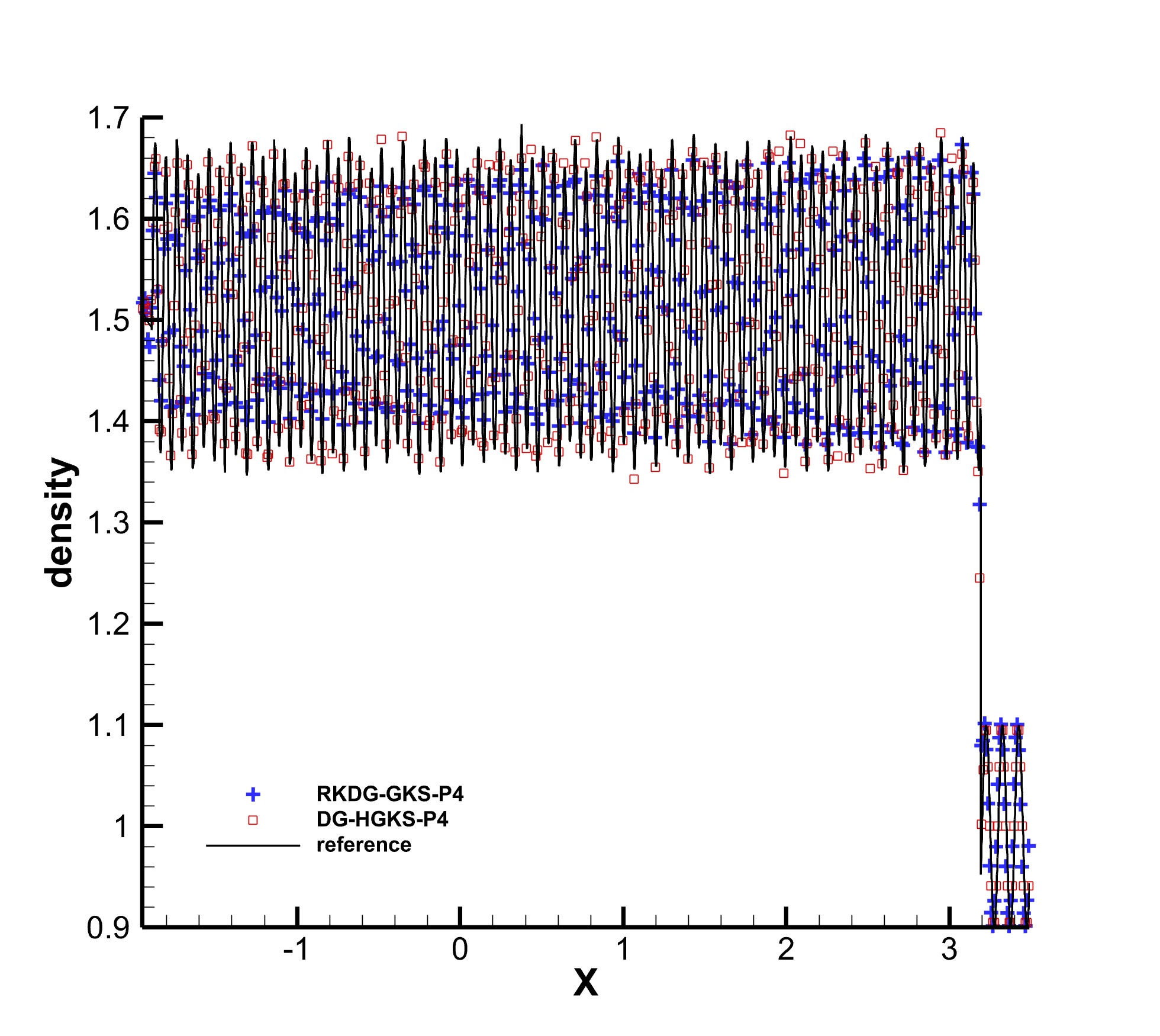}
	 \caption*{(b) {The density distribution  of Example 4.5 \eqref{gp} with $N=1000.$}}
  \end{minipage} }
\caption{Comparison results with the DG-HGKS and RKDG-GKS, fifth-order($p^4$).}
  \label{fig:rkcomp}
\end{figure}  

\begin{figure}[h!]
	\centering
	\begin{minipage}[t]{0.5\textwidth}  
	\centering
	 \includegraphics[width=\textwidth]{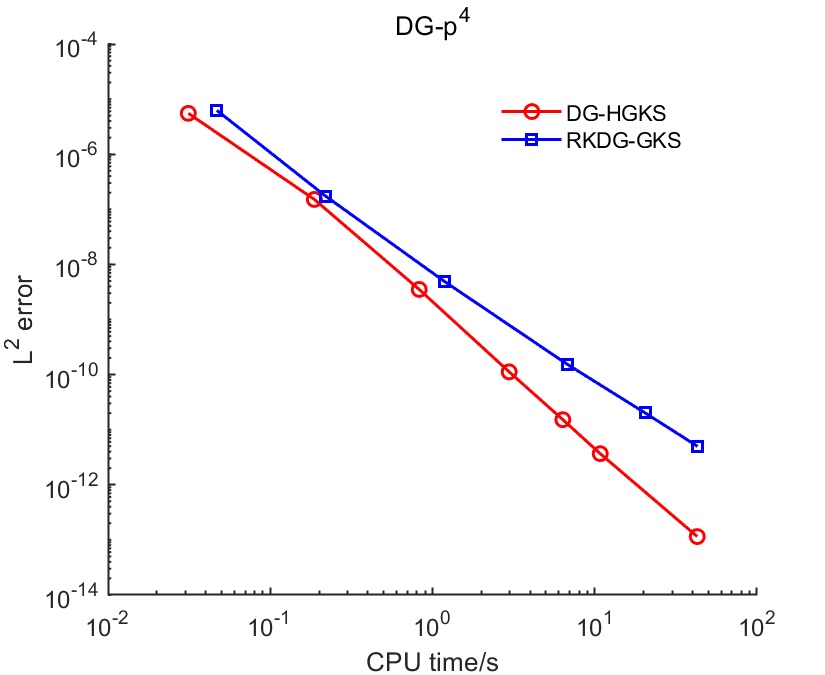}
          \end{minipage}
\caption{The CPU time vs  $L^2$ error between RKDG and DG-HGKS.}
  \label{fig:rkcpu}
\end{figure}  
\begin{table}[h]
\centering
\begin{tabular}{|l|c|c|c|}
\hline
 & $N=40$ & $N=80$ & $N=160$ \\ \hline
RKDG-GKS-$p^4$ & 1.078s & 3.750s & 14.438s \\ \hline
DG-HGKS-$p^4$ & 0.234s & 1.281s & 4.203s \\ \hline
\end{tabular}
\caption{The CPU time per DOF between RKDG and DG-HGKS.}
\label{tab:RKDOF}
\end{table}

\subsection*{B.4  Comparison with  the existing DG-GKS}
To further demonstrate the significant computational advantages of the proposed DG-HGKS method, this appendix provides a comparative analysis with the existing single-stage DG-GKS , with references to\cite{NGX}.

On one hand, as introduced in the introduction, our DG-HGKS method has successfully broken through the third-order accuracy limitation of the traditional DG-GKS method and can achieve  fifth-order accuracy in both space and time.
On the other hand, even at the same order of accuracy, our DG-HGKS method still demonstrates superior computational performance. 
Taking the Shu-Osher problem (Example 2 in \cite{NGX}) as an example for comparison,
the Fig. \ref{fig:SO2comp}  illustrates the comparison of the two methods under the third-order accuracy with  $N = 400$  cells. 
Which shown that the scheme constructed in this paper has a superior agreement with the reference solution and more precise accuracy, 
fully demonstrating its significant advantages in computational performance.
In terms of computational cost, 
the DG-GKS method takes 11.719 seconds for this example, 
while the DG-HGKS method only takes 6.891 seconds, 
demonstrating a higher computational efficiency.

\begin{figure}[h!]
	\centering
	\subfloat{
	\begin{minipage}[t]{0.48\textwidth}  
	\centering
	 \includegraphics[width=\textwidth]{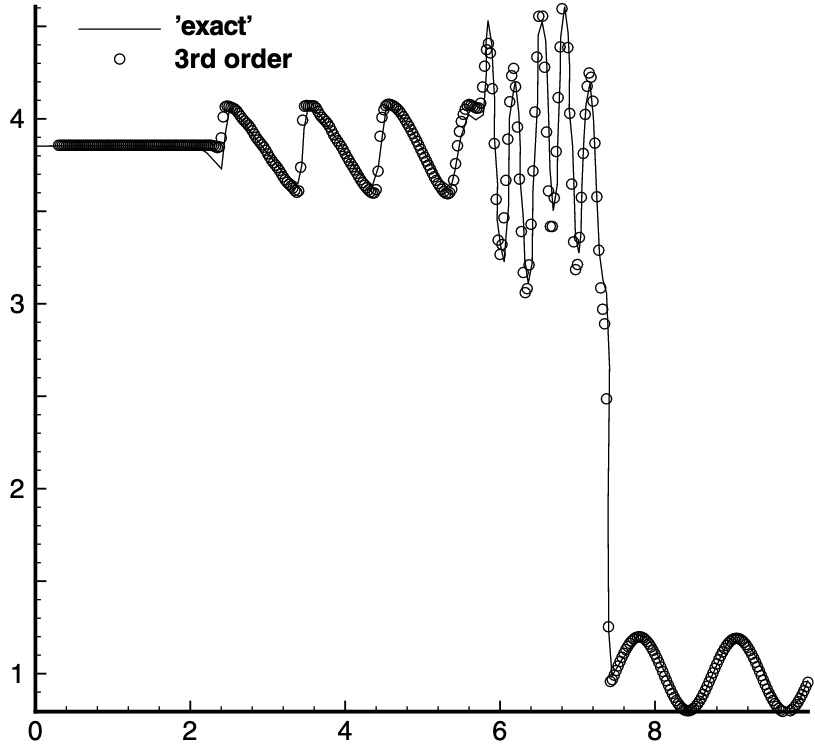}
	 \caption*{(a) DG-GKS,third-order,from\cite{NGX}.}
  \end{minipage}
   \begin{minipage}[t]{0.5\textwidth}  
	\centering
	 \includegraphics[width=\textwidth]{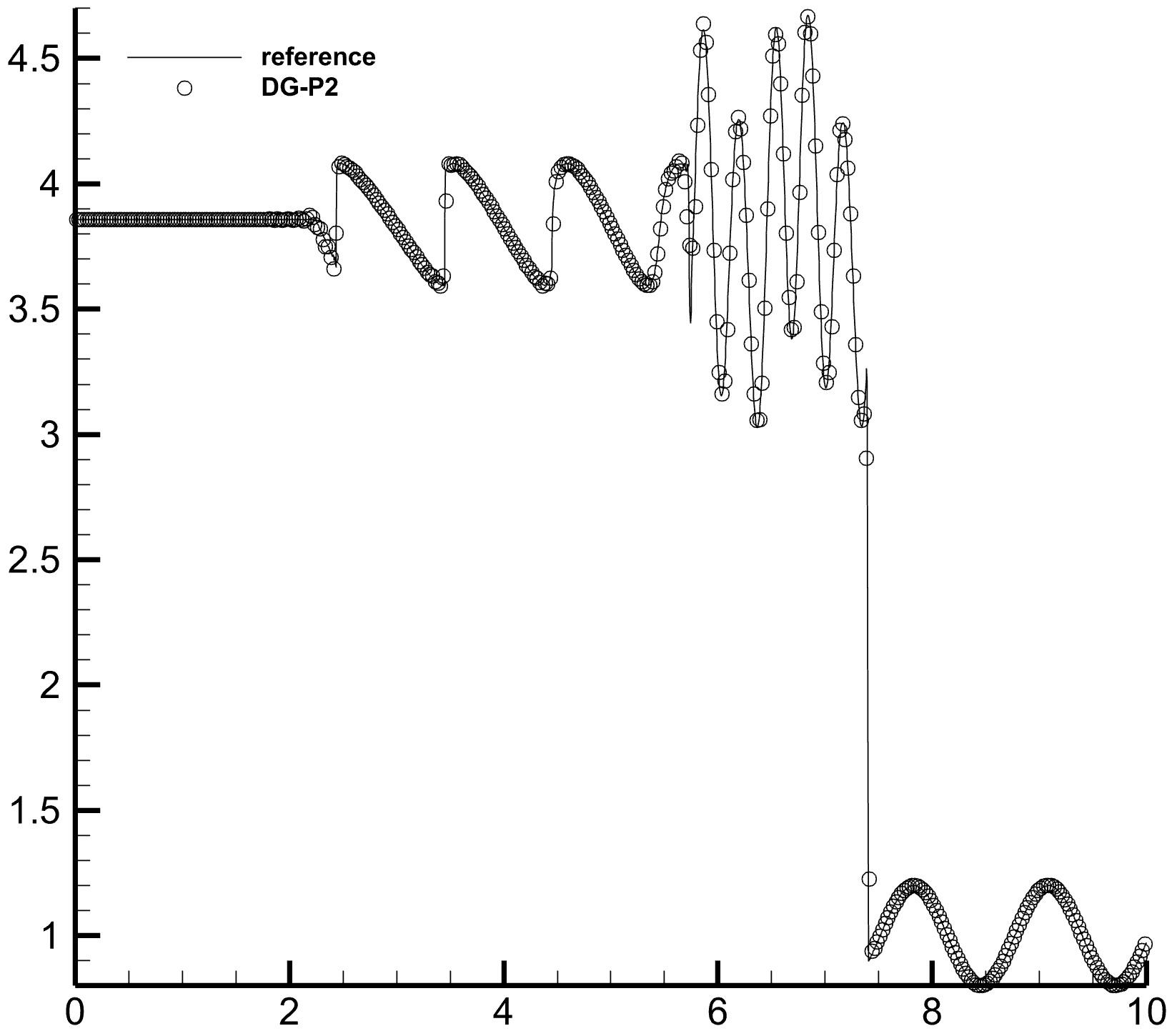}
	 \caption*{(b) DG-HGKS, third-order($p^2$).}
  \end{minipage} }
\caption{The density distribution  of  shock and sound wave interaction \cite{NGX} with  $N=400$.}
  \label{fig:SO2comp}
\end{figure}  


\subsection*{B.5  Comparison with  the DG-ADER}
To highlight the advantages of the DG-HGKS method constructed in this paper in single-stage computational methods, we compared it with the single-stage DG-ADER method.
The same limiter strategy as that in this paper was adopted for DG-ADER   to ensure the consistency.
Fig. \ref{fig:BW2comp} shows a comparison  for 
1D Woodward-Colella blast wave \eqref{bw}
and the modified shock/turbulence interaction \eqref{gp}. 
It is evident that the DG-HGKS method in this paper significantly outperforms the DG-ADER method in terms of numerical resolution.

In terms of computational cost, 
this paper compares the computational time of the two methods.
For the Woodward-Colella blast wave \eqref{bw},
the DG-ADER method takes 147.078 seconds, and the DG-HGKS method  takes 130.891 seconds.
For the modified shock/turbulence interaction \eqref{gp},
the DG-ADER method takes 526.064 seconds, while the DG-HGKS method  takes 469.936 seconds. 
This clearly demonstrates that the DG-HGKS method has a higher computational efficiency.

\begin{figure}[h!]
	\centering
	\subfloat{
	\begin{minipage}[t]{0.5\textwidth}  
	\centering
	 \includegraphics[width=\textwidth]{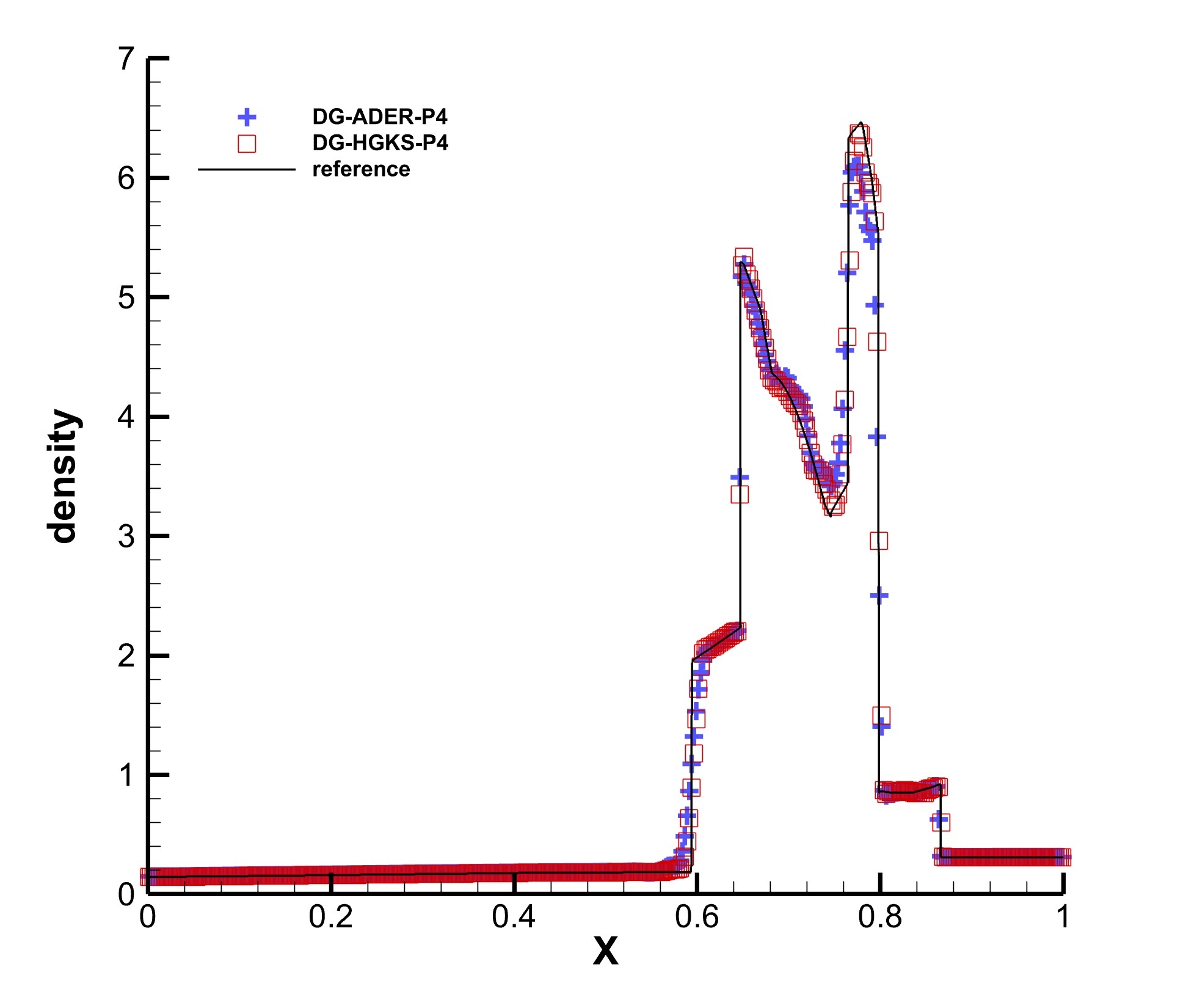}
	 \caption*{(a) {The density distribution  of Example 4.4 \eqref{bw} with $N=400.$}}
  \end{minipage}
   \begin{minipage}[t]{0.5\textwidth}  
	\centering
	 \includegraphics[width=\textwidth]{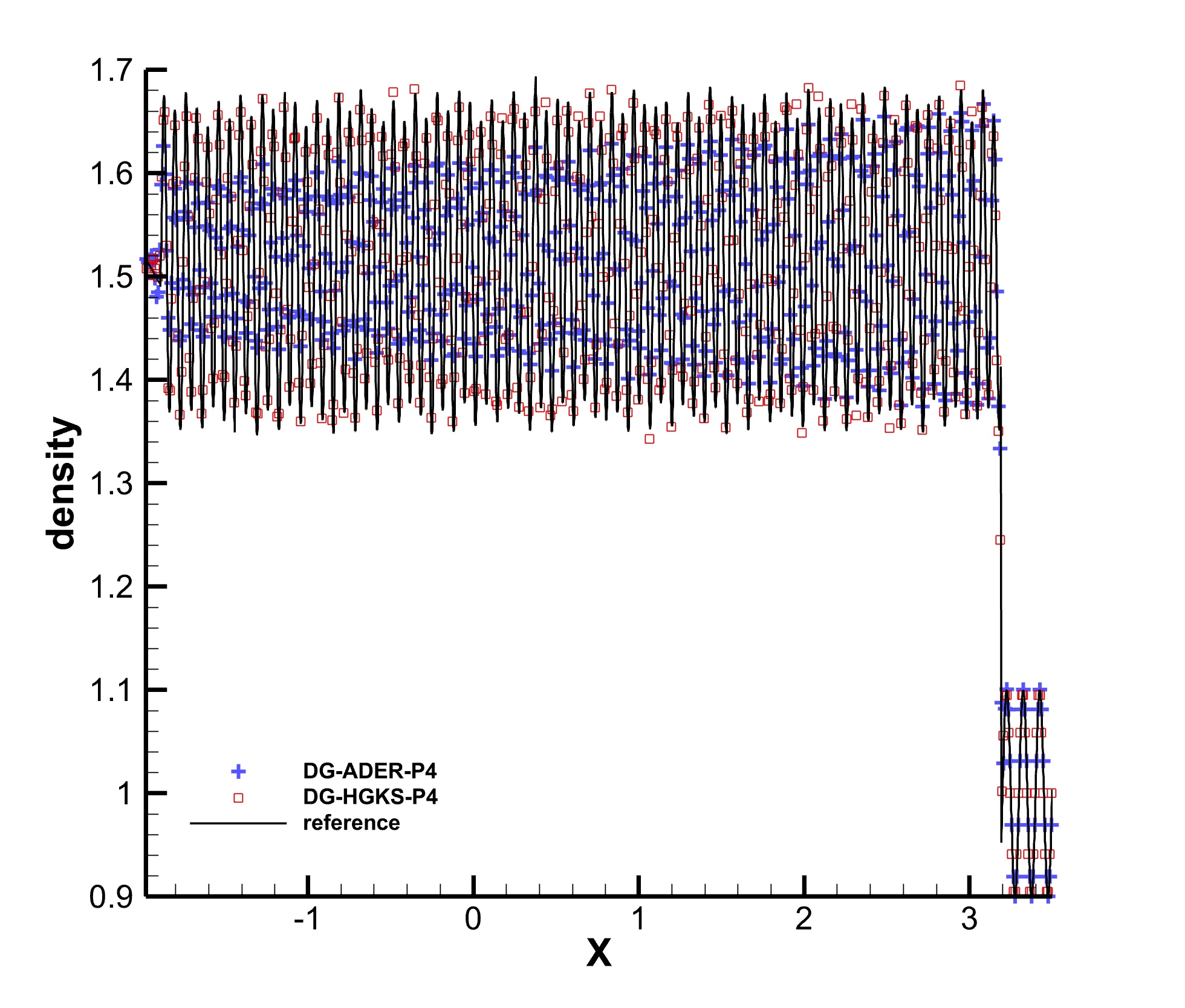}
	 \caption*{(b) {The density distribution  of Example 4.5 \eqref{gp} with $N=1000.$}}
  \end{minipage} }
\caption{Comparison results with the DG-HGKS and DG-ADER, fifth-order($p^4$).}
  \label{fig:BW2comp}
\end{figure}  

%


{\section{The extension of the viscous Navier-Stokes equations}}
\renewcommand{\thetable}{C.\arabic{table}}
\renewcommand{\thefigure}{C.\arabic{figure}}
\renewcommand{\theequation}{C.\arabic{equation}}
In this paper, the DG-HGKS scheme primarily focuses on the inviscid Euler equations. 
There are two distinct approaches to extend the DG-HGKS method to the viscous Navier-Stokes equations. 
The first method, following the classical GKS framework\cite{13,15,RXD}, which offers a more natural treatment of both viscous and inviscid flows. 
However, extending this method to higher-order formulations leads to a significant increase in the complexity of space-time transformations. 
Alternatively, one may adopt Luo's ADER method\cite{LDM}, which effectively reduces the complexity of space-time transformation in viscous flow equations.

In this appendix, we introduce the approach for solving the one-dimensional viscous Navier-Stokes equations using the DG-HGKS method. 
The viscous Navier-Stokes equations in 1D is given by\cite{NS} 
$$
\begin{cases}
W_t + F(W)_x = G(W, W_x)_x, \\
W(x,0) = W_0(x),
\end{cases}
$$
where the conservative variables $W$ and the flux $F(W)$ are the same as the Euler equations \eqref{eq:euler}.
And
$$
G(W, W_x) = \begin{pmatrix}
0 \\
\frac{4}{3} \mu U_x \\
\frac{4}{3} \mu U U_x + \frac{\mu}{Pr(\gamma - 1)} T_x
\end{pmatrix}
$$
where $T$ is the temperature, $Pr$ is the Prandtl number and $\mu$ is the molecular viscosity computed by the Sutherland’s law $\mu = \frac{1 + c }{T + c}T^{\frac{3}{2}}$,
here $c = \frac{110.4}{T_{\infty}}$ and $T_{\infty}=288K$ is the reference temperature.

The DG scheme is 
$$
\begin{aligned}
w_m(t_{n+1}) &= w_m(t_n) + \frac{2m+1}{\Delta x} \left[ F_{\text{volume}} - G_{\text{volume}} \right]\\
&- \frac{2m+1}{\Delta x} \left[ \hat{F}_{i+\frac{1}{2}} - \hat{F}_{i-\frac{1}{2}} \right] + \frac{2m+1}{\Delta x} \left[ \hat{G}_{i+\frac{1}{2}} - \hat{G}_{i-\frac{1}{2}} \right] 
\end{aligned}
$$
the objective is to extend the high-order DG-HGKS scheme presented in this paper to viscous flows. Therefore, here only focus on the computation of the high-order terms,
employing the Gaussian quadrature rule to approximate the flux terms
$$
\hat{F}_{i+\frac{1}{2}} = \int_{t_n}^{t_{n+1}} F(W)_{i+\frac{1}{2}} v_m(x_{i+\frac{1}{2}}^{-}) dt 
\approx  v_m(x_{i+\frac{1}{2}}^{-}) \sum_g w_g F(\mathcal{W}^e(x_{i+\frac{1}{2}}, t_g))
$$
$$
\hat{G}_{i+\frac{1}{2}} =  \int_{t_n}^{t_{n+1}} G(W, W_x)_{i+\frac{1}{2}} v_m(x_{i+\frac{1}{2}}^{-}) \, dt \\
\approx  v_m(x_{i+\frac{1}{2}}^{-}) \sum_g w_g G(\mathcal{W}^e(x_{i+\frac{1}{2}}, t_g), \mathcal{W}_x^e(x_{i+\frac{1}{2}}, t_g))
$$
where $t_g$ and $w_g$ are the corresponding Gaussian points and weights. 
For the equilibrium state, the approach is consistent with the previously introduced method, 
employing a Taylor expansion in time for approximation
$$
\begin{aligned}
\mathcal{W}^e(x_{i+\frac{1}{2}}, \tau) & \approx W^e(x_{i+\frac{1}{2}}, 0^+) + \sum_{k=1}^{4} \left[ \frac{\partial_t^{(k)} W^e(x_{i+\frac{1}{2}}, 0^+)}{k!} \right] \frac{\tau^k}{k!}\\
\mathcal{W}_x^e(x_{i+\frac{1}{2}}, \tau) & \approx W_x^e(x_{i+\frac{1}{2}}, 0^+) + \sum_{k=1}^{3} \left[ \frac{\partial_t^{(k)} W_x^e(x_{i+\frac{1}{2}}, 0^+)}{k!} \right] \frac{\tau^k}{k!}
\end{aligned}
$$
where $\tau=t_g-t_n$. 
For each order of time derivatives ${\partial_t^{(k)} W^e(x_{i+\frac{1}{2}}, 0^+)}$, 
the Lax-Wendroff procedure is used to convert them into a combination of spatial derivatives ${\partial_x^{(k)} W^e(x_{i+\frac{1}{2}}, 0^+)}$. 
As for each order of spatial derivatives can also be calculated by \eqref{eq12}
$$
\begin{aligned}
W^e(x_{i+\frac{1}{2}}, 0^+) &= \omega_e W(x^{-}_{i+\frac{1}{2}}, 0^+) + (1 - \omega_e) W(x^{+}_{i+\frac{1}{2}}, 0^+)\\
{\partial_x^{(k)} W^e(x_{i+\frac{1}{2}}, 0^+)} &= \omega_e {\partial_x^{(k)} W(x^{-}_{i+\frac{1}{2}}, 0^+)} + (1 - \omega_e) {\partial_x^{(k)} W(x^{+}_{i+\frac{1}{2}}, 0^+)}
\end{aligned}
$$
in the fifth-order scheme, $ k = 1, \ldots, 4$.

As for the volume term 
$$
\begin{aligned}
F_{\text{volume}} &= \int_{t_n}^{t_{n+1}} \int_{x_{i-\frac{1}{2}}}^{x_{i+\frac{1}{2}}} F(W) v'_m(x) \, dx \, dt,
\\
G_{\text{volume}} &= \int_{t_n}^{t_{n+1}} \int_{x_{i-\frac{1}{2}}}^{x_{i+\frac{1}{2}}} G(W, W_x) v'_m(x) \, dx \, dt.
\end{aligned}
$$
Similar to the volume term of the inviscous part \eqref{3-x},\eqref{vx} and \eqref{wht}.
$$
\begin{aligned}
&\int_{t_n}^{t_{n+1}} \int_{x_{i-\frac{1}{2}}}^{x_{i+\frac{1}{2}}} G(W, W_x) v'_m(x) \, dx \, dt \\
&\approx \sum_{g_x} \sum_{g_t} w_{g_x} w_{g_t} P_m'(s_{g_x}) G \left( \mathcal{W}(x_{g_x}, t_{g_t}), \mathcal{W}_x(x_{g_x}, t_{g_t}) \right)
\end{aligned}
$$
here $x_{g_x} = \frac{\Delta x}{2} s_{g_x} + x_i, \, t_{g_t} = \frac{\Delta t}{2} s_{g_t} + t_n, \,$ 
where $s_{g_x}, s_{g_t}$  are the standard Gauss points and  $w_{g_x}, w_{g_t}$ are the weights,
similar to the high-order term of the flux
$$
\begin{aligned}
\mathcal{W}(x_{g_x}, \tau) &\approx W(x_{g_x}, 0^+) + \sum_{k=1}^{4} \left[ \partial_t^{(k)} W(x_{g_x}, 0^+) \right] \frac{\tau^k}{k!}\\
\mathcal{W}_x(x_{g_x}, \tau) &\approx W_x(x_{g_x}, 0^+) + \sum_{k=1}^{3} \left[ \partial_t^{(k)} W_x(x_{g_x}, 0^+) \right] \frac{\tau^k}{k!}
\end{aligned}
$$
The temporal derivatives of various orders are still converted into combinations of spatial derivatives through the Lax-Wendroff procedure,
and the spatial derivatives of various orders are directly obtained by numerical solution \eqref{so1d}.
Regarding the Lax-Wendroff procedure of the Navier-Stokes equations, 
please refer to the literature\cite{LDM} for details.

We test the fifth-order DG-HGKS scheme on several one-dimensional Navier-Stokes examples to evaluate its accuracy, effectiveness, and robustness.
The reference solution with a small viscosity coefficient $\mu$ is consistent with that of the Euler equations,
while the reference solution with a large viscosity coefficient $\mu$ is obtained using a highly refined mesh.

\textbf{Example 1.} First we test the accuracy of the DG-HGKS method for Navier-Stokes equation,
the example 4.1 is extended to the viscous Navier-Stokes equations with the same initial condition \eqref{1dorder} and periodic boundary condition.
The parameters $Pr$ is set as $\frac{2}{3}$.
The errors and numerical orders with different parameters $\mu=0.0001$ and $\mu=0.1$ are shown in Table \ref {tab:acvs1d},
which demonstrate that the DG-HGKS scheme can achieve the fifth-order
when the advection is dominated, i.e., $\mu = 0.0001$, 
and the order is decreased while the diffusion is dominated, i.e., $\mu = 0.1$.
Specifically, the numerical convergence rates are estimated using the asymptotic convergence error proposed in \cite{[36]}, as there are no exact solutions for this problem.

 \begin{table}[h]
\centering
\begin{tabular}{|c|c|cccc|cccc|}
\hline
 &   & \multicolumn{4}{c|}{ \textbf{DG-HGKS with $\mu=0.0001$}} & \multicolumn{4}{c|} {\textbf{DG-HGKS with $\mu=0.1$}} \\ \hline
 &{ { N }}   & { {$L^1$ error}} & { {order}} & { {$L^\infty$ error}} & { {order}} & { {$L^1$ error}} & { {order}} & { {$L^\infty$ error }}& { {order}} \\ \hline
\multirow{4}{*}{{ {$\boldsymbol{p^4}$}}} 
& 20   & 1.93e-06 & --      & 1.26e-06 & --     & 2.67e-06 & --      &1.29e-06 & --\\
& 40   & 5.57e-08 & 5.11 & 3.50e-08 & 5.17 & 8.61e-08 & 4.96 & 3.41e-08 & 5.25 \\
& 80   & 1.95e-09 & 4.84 & 1.23e-09 & 4.83 & 2.93e-09 & 4.88 & 1.15e-09 & 4.88\\
& 160 & 5.88e-11 & 5.05 & 3.71e-11 & 5.05 & 1.08e-10 & 4.76 & 4.38e-11 & 4.72\\ \hline
\end{tabular}
\caption{{\label{tab:acvs1d} Accuracy test  for 1D of Example 1 with initial condition \eqref{1dorder}, $t=1.0$.}}
\end{table}

%
%

\textbf{Example 2.} The Shu-Osher  problem in example 4.3 is extended to the viscous Navier-Stokes equation,
which with the same initial condition\eqref{so} and boundary condition.
The parameters $Pr$ is set as $0.72$.
The density distributions for fifth-order with   $N=400$ under different parameters $\mu=0.0001$ and $\mu=0.01$ are shown in Fig. \ref{fig:SO-vs}.  
It can be seen that the DG-HGKS scheme effectively captures discontinuities in viscous flows. 
Additionally, the numerical dissipation increases with higher viscosity.

\begin{figure}[h!]
	\centering
	 \subfloat{
	\begin{minipage}[t]{0.45\textwidth}  
	\centering
	 \includegraphics[width=\textwidth]{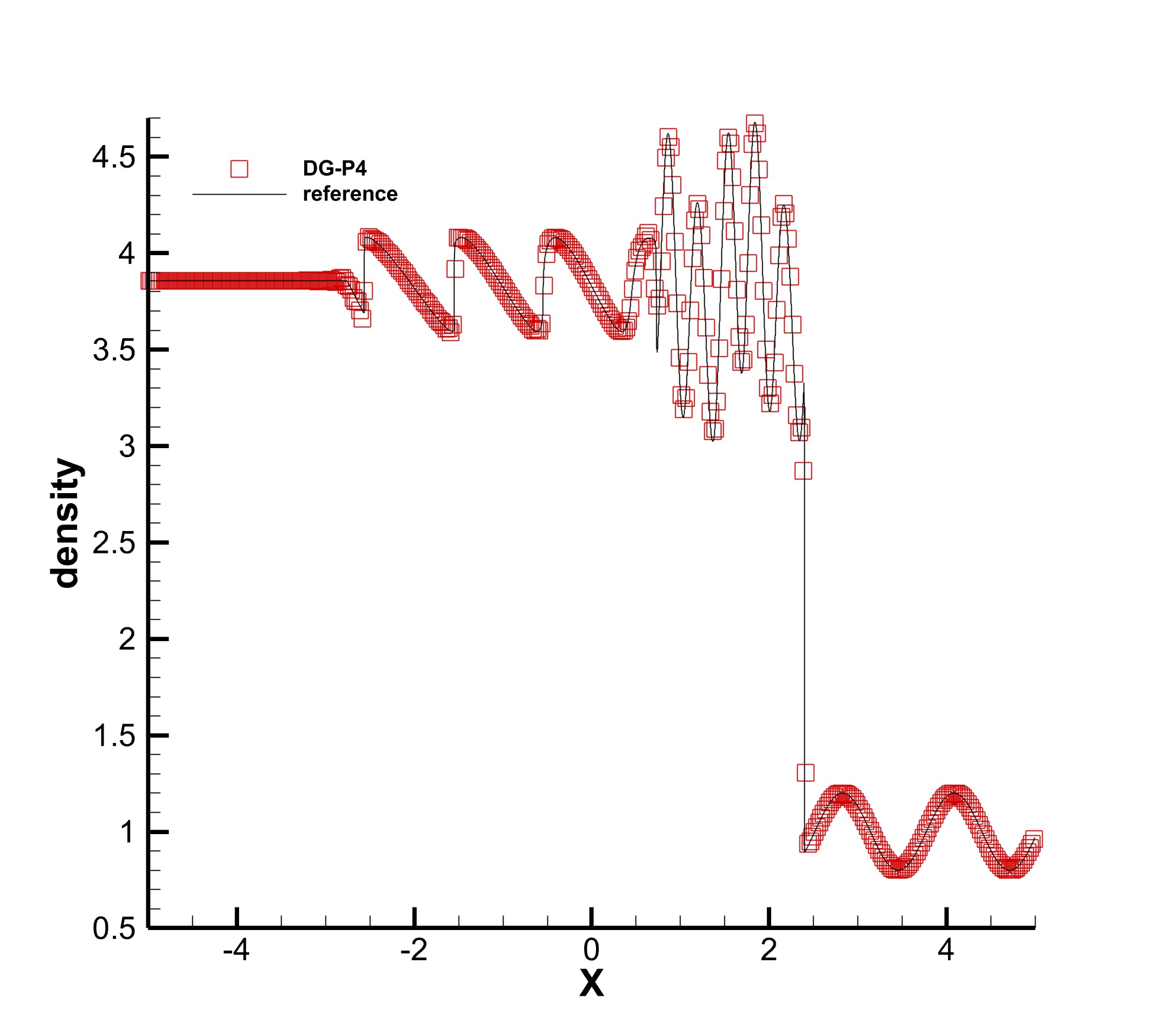}
  \end{minipage}
   \begin{minipage}[t]{0.45\textwidth}  
	\centering
	 \includegraphics[width=\textwidth]{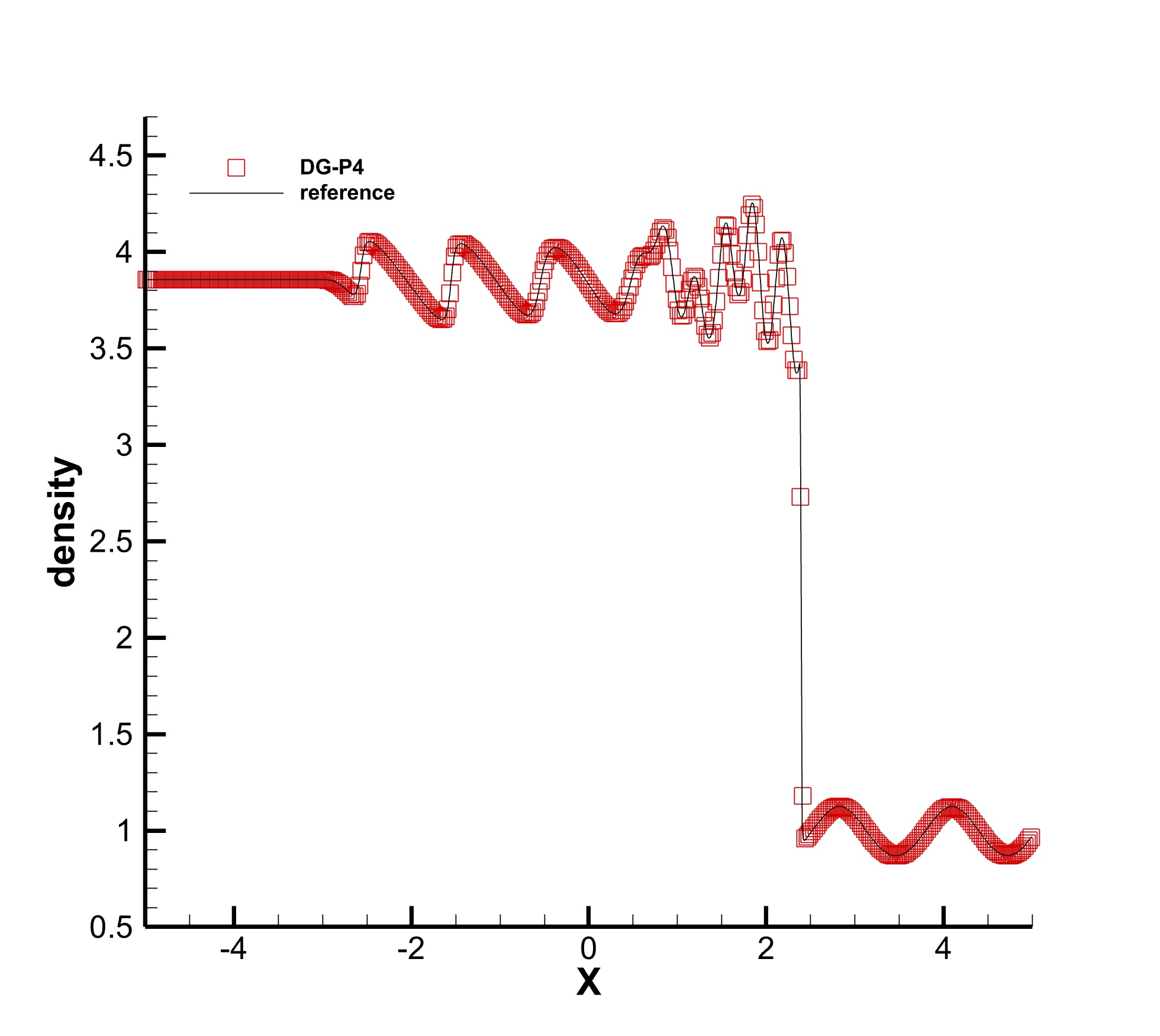}
  \end{minipage}
 }
  \caption{The density distribution of  Example 4.3 \eqref{so} in fifth-order($p^4$), with $N=400, t=1.8$.
Left: $\mu=0.0001$, right: $\mu=0.01$.}
  \label{fig:SO-vs}
\end{figure}

\textbf{Example 3.} The Woodward-Colella blast wave in example 4.4 is extended to the viscous Navier-Stokes equation,
which with the same initial condition\eqref{bw} and boundary condition.
The parameters $Pr$ is set as $0.72$.
The density distributions for fifth-order with  $N=800$ under different parameters $\mu=0.001$ and $\mu=0.01$ are shown in Fig. \ref{fig:BW-vs}. 
The numerical results demonstrate that the DG-HGKS scheme presented in this paper remains highly effective for viscous flows featuring strong shock waves.

\begin{figure}[h!]
	\centering
	\subfloat{
	\begin{minipage}[t]{0.45\textwidth}  
	\centering
	 \includegraphics[width=\textwidth]{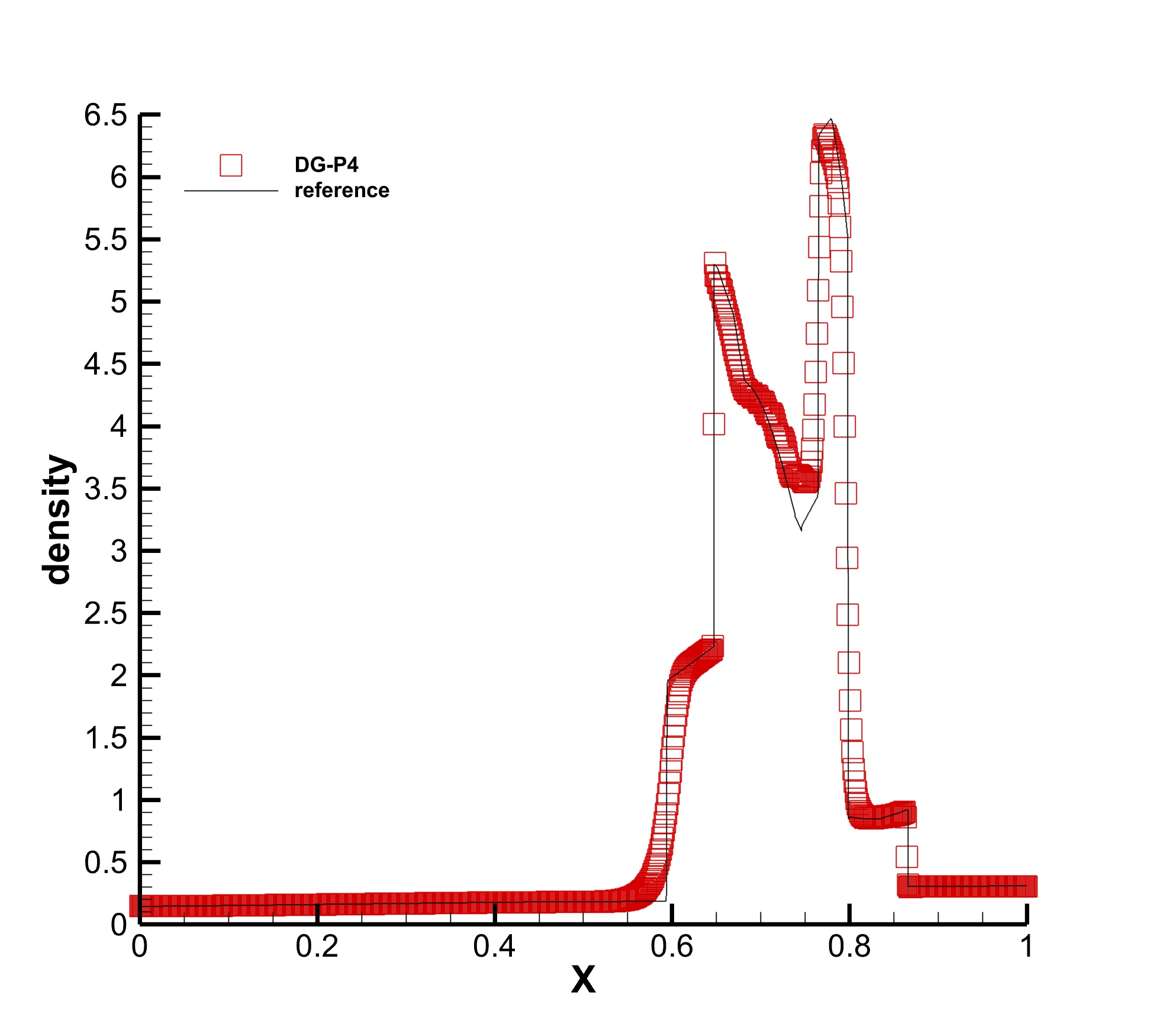}
  \end{minipage}
   \begin{minipage}[t]{0.45\textwidth}  
	\centering
	 \includegraphics[width=\textwidth]{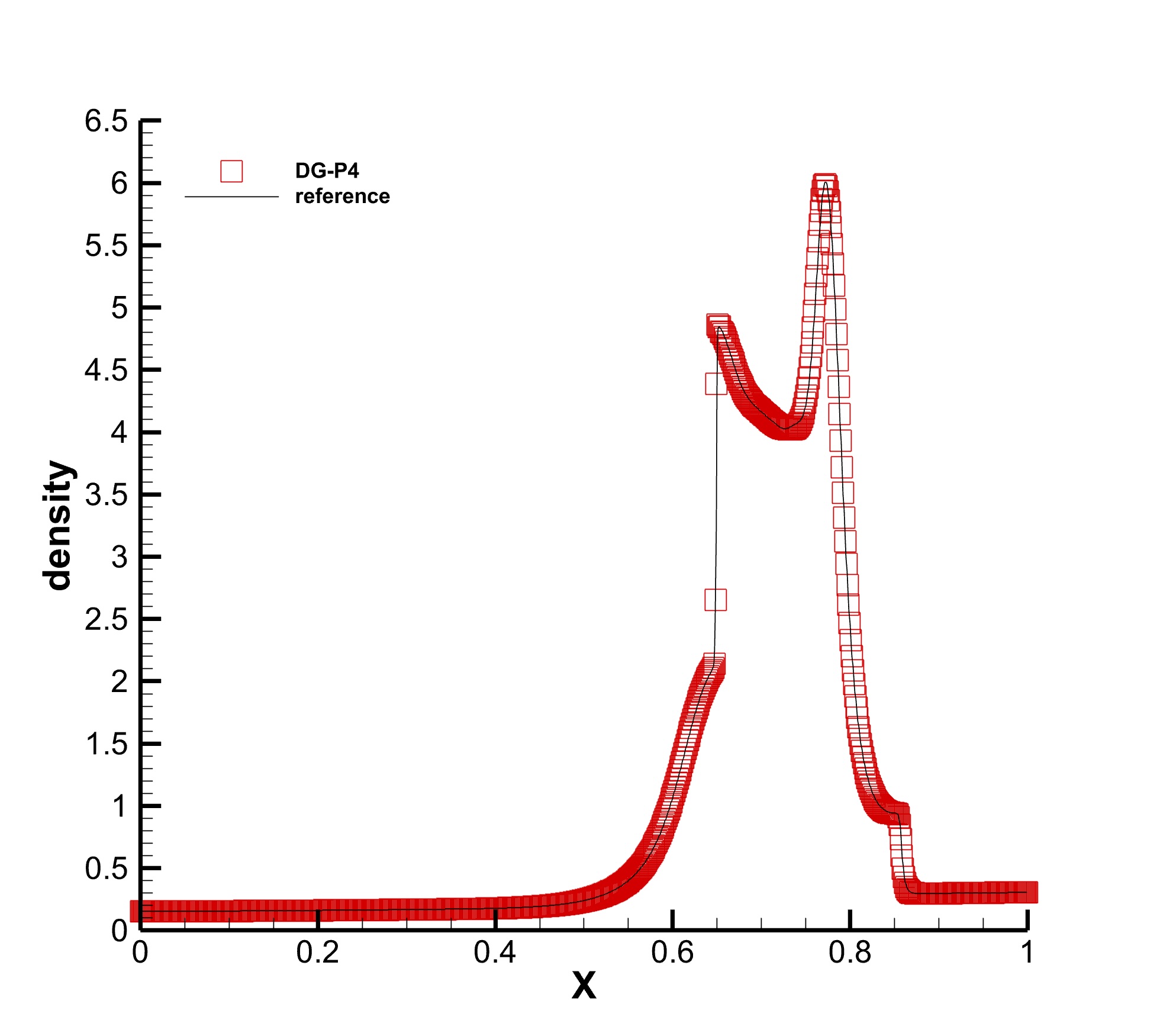}
  \end{minipage}
 }
  \caption{The density distribution of  Example 4.4 \eqref{bw} in fifth-order($p^4$), with  $N=800, t=0.038.$
  Left: $\mu=0.001$, right: $\mu=0.01$. }
  \label{fig:BW-vs}
\end{figure}

\end{appendices}

\end{document}